\documentclass[notitlepage,leqno,10pt]{article}
\usepackage{amsmath, amsfonts, amssymb, amscd,
amsthm,enumitem,extarrows}
\usepackage[all]{xy}
\usepackage{latexsym}
\usepackage{amsmath}
\usepackage{amsfonts}  
\usepackage{amssymb}
\usepackage{amscd}
\usepackage{mathtools, amsthm,enumitem}
\usepackage{hyperref}  
\newcommand{\eps}{\varepsilon}
\newcommand*\rfrac[2]{{}^{#1}\!/_{#2}}

\renewcommand{\d}{\delta }
\newcommand{\D }{\Delta }

\renewcommand{\l }{\lambda }

\newcommand{\n }{\nabla }

\newcommand{\ov}{\overline}
\newcommand{\intbar}{\mathop{\int\makebox(-13.5,0){\rule[4pt]{.7em}{0.3pt}}%
\kern-6pt}\nolimits}
\newcommand{\wtilde }{\widetilde}

\newcommand{\be}{\begin{equation}}
\newcommand{\ee}{\end{equation}}
\newcommand{\bes}{\begin{equation*}}
\newcommand{\ees}{\end{equation*}}
\newcommand{\ba}{\begin{eqnarray}}
\newcommand{\ea}{\end{eqnarray}}
\newcommand{\bas}{\begin{eqnarray*}}
\newcommand{\eas}{\end{eqnarray*}}
\newenvironment{pf}{\noindent{\sc Proof}.\enspace}{\rule{2mm}{2mm}\medskip}
\newenvironment{pfn}{\noindent{\sc Proof}}{\rule{2mm}{2mm}\medskip}

\newcommand{\R}{\mathbb{R}}

\newcommand{\Z}{\mathbb{Z}}

\newcommand{\N}{\mathbb{N}}

\DeclareMathAlphabet{\mymathbb}{U}{BOONDOX-ds}{m}{n}

\author{Martin MAYER$^{a}$,\;\; Cheikh Birahim NDIAYE$^{b}$}
\date{}

\title{\bf Fractional Yamabe problem on locally flat conformal infinities of Poincar\'e-Einstein manifolds}

\begin{document}

\newtheorem{lem}{Lemma}[section]
\newtheorem{pro}[lem]{Proposition}
\newtheorem{thm}[lem]{Theorem}
\newtheorem{rem}[lem]{Remark}
\newtheorem{cor}[lem]{Corollary}
\newtheorem{df}[lem]{Definition}

\maketitle

\begin{center}

{\small

\noindent  $^{a}$ Scuola Superiore Meridionale, \\Via Mezzocannone 4, Naples, ITALY.

}
\
\
\
\

{\small

\noindent  $^{b}$ Department of Mathematics  of Howard University \\  Annex 3, Graduate School of Arts and Sciences, \# 217 \\ DC 20059 Washington, USA.

}

\
\
{\small

\noindent

}

\end{center}

\footnotetext[1]
{
E-mail: m.mayer@ssmeridionale.it, \;
cheikh.ndiaye@howard.edu.\\
\thanks
{
\\C. B. Ndiaye was partially supported by NSF grant DMS--2000164.
\\  M.Mayer has been supported by the Italian MIUR Department of Excellence grant CUP E83C18000100006. 
}
}

\

\begin{center}
{\bf Abstract}

\end{center}
We study the fractional Yamabe problem first considered by Gonzalez-Qing\cite{gq} on the conformal infinity \;$(M^n, \;[h])$\; of a Poincar\'e-Einstein manifold \;$(X^{n+1}, \;g^+)$\; with either \;$n=2$\; or $n\geq 3$\; and \;$(M^n, \;[h])$\;  locally flat - namely\;$(M, h)$\;is  locally conformally flat.
However, as for the classical Yamabe problem, because of the involved quantization phenomena, the variational analysis of the fractional one exhibits a local situation and also a global one. The latter global situation 
includes the case of  conformal infinities of Poincar\'e-Einstein manifolds of dimension either \;$n=2$\; or of dimension  \;$n\geq 3$ and  which are locally flat, and hence the minimizing technique of Aubin\cite{au}-Schoen\cite{sc} in that case clearly requires an analogue of the positive mass theorem of Schoen-Yau\cite{sy1},
which is not known to hold. Using the algebraic topological argument of Bahri-Coron\cite{bc}, we bypass the latter positive mass issue and show that any conformal infinity of a Poincar\'e-Einstein manifold of dimension either \;$n=2$\; or of dimension \;$n\geq 3$\; and  which is locally flat admits a Riemannian metric of constant fractional scalar curvature.

\begin{center}
 \
 \


\bigskip\bigskip

\bigskip

\centerline{\bf AMS subject classification:  53C21, 35C60, 58J60, 55N10.}

\bigskip
\end{center}

\begin{center}
$_{}$
\vspace{100pt}
\tableofcontents
\end{center}

\section{Introduction and statement of the results}
In recent years there has been a lot of study about fractional order operators in the context of elliptic theory with nonlocal operators, nonlinear diffusion involving nonlocal operators, nonlocal aggregations and balance between nonlinear diffusions and nonlocal attractions. The elliptic theory of  fractional order operators is well understood in many directions like semi-linear equations, free-boundary value problems and non-local minimal surfaces (see \cite{cafroq}, \cite{cafsyl}, \cite{cafsoug}, \cite{cafval}, \cite{cafvas}, \cite{fabkenser}). Furthermore a connection to classical conformally covariant operators arising in conformal geometry and their associated conformally invariant geometric variational problems is established (see \cite{cg}, \cite{fg}, \cite{qr}, \cite{gz}, \cite{guil}, \cite{gms}, \cite{gq}). In this paper we are interested in the latter aspect of fractional order operators, precisely in the fractional Yamabe problem first considered by Gonzalez-Qing\cite{gq}. 
\vspace{6pt}

\noindent
To discuss the fractional Yamabe problem, we first recall some definitions in the theory of asymptotically hyperbolic metrics. Given \;$X=X^{n+1}$\; a smooth manifold with boundary $M=M^n$\; and \;$n\geq 2$, we say that \;$\varrho$\; is a defining function of the boundary \;$M$\; in \;
$X$,\; if 
$$
\varrho>0\;\;\text{ in }\;\;X,\;\;\varrho=0\;\;\text{ on }\;\;M\;\;\text{ and }\;\;d\varrho\neq 0\;\;\text{ on }\;\; M.
$$
A Riemannian metric \;$g^+$\; on \;$X$\; is said to be conformally compact, if for some defining function \;$\varrho$\; the Riemannian metric 
\begin{equation}\label{eq:cmetric}
 g:=\varrho^2g^+
\end{equation}
 extends to \;$\ov{X}:=X\cup M$\; so that \;$(\ov{X}, \, g)$\; is a compact Riemannian manifold with boundary \;$M$ and interior \;$X$. Clearly this induces a conformal class of Riemannian metrics 
$$ [h]=[ g|_{TM}]$$ on \;$M$, where $TM$ denotes the tangent bundle of \;$M$, when the defining functions \;$\varrho$\; vary. The resulting conformal manifold \;$(M, [h])$\; is called conformal infinity of \;$(X, \;g^+)$. Moreover  a Riemannian metric \;$g^+$\; in \;$X$\; is said to be asymptotically hyperbolic, if it is conformally compact and its sectional curvature tends to \;$-1$\; as one approaches the conformal infinity of \;$(X, \;g^+)$,
 which is equivalent to 
\begin{equation*}
|d\varrho|_{\bar g}=1
\end{equation*}
 on\;$M$, see \cite{mazzeo1}, and in such a case \;$(X, \;g^+)$\; is called an asymptotically hyperbolic manifold. Furthermore  a Riemannian metric \;$g^+$\; on \;$X$\; is said to be conformally compact Einstein or Poincar\'e-Einstein (PE), if it is asymptotically hyperbolic and satisfies the Einstein equation 
\begin{equation*}
Ric_{g^+}=-ng^+,
\end{equation*}
where $Ric_{g^+}$\; denotes the Ricci tensor of \;$(X, \;g^+)$.
\vspace{6pt}

\noindent
On one hand for every asymptotically hyperbolic manifold \;$(X, \;g^+)$\; and every choice of the representative \;$h$\; of its conformal infinity \;$(M, [h])$\; 
there exists a unique geodesic defining function \;$ y $\; of \;$M$\; in \;$X$\; 
such that  in a tubular neighborhood of \;$M$\; in \;$X$  the Riemannian metric \;$g^+$\; takes the following normal form
\begin{equation}\label{eq:uniqdef}
g^+=\frac{d y ^2+h_{ y }}{ y ^2},
\end{equation}
where \;$h_{ y }$\; is a family of Riemannian metrics on \;$M$\; satisfying \;$h_0=h$.  
We say that the conformal infinity \;$(M, \;[\hat h])$\; of an asymptotically hyperbolic manifold \;$(X, \;g^+)$\; is locally flat, 
if \;$h$\; is locally conformally flat, and clearly this is independent of the representative \;$h$\; of \;$[h]$. 
Moreover  we say that \;$(M, [h])$\; is umbilic, if \;$(M, h)$\; is umbilic in \;$(X, \; g)$ \;where \;$g$\; is given by \eqref{eq:cmetric} 
and \;$ y $\; is the unique geodesic defining function given by \eqref{eq:uniqdef}, 
and this is again independent of the representative \;$h$\; of \;$[h]$, as easily seen from the uniqueness of the normal form \eqref{eq:uniqdef} or Lemma 2.3 in \cite{gq}. 
Similarly we say that \;$(M, [h])$\; is minimal if \;$H_{g}=0$\; with \;$H_{g}$\; denoting the mean curvature of \;$(M,  \;h)$\; in \;$(\ov X, \;g)$ with respect to the inward direction, and this is again clearly independent of the representative of \;$h$\; of \;$[h]$, as easily seen from Lemma 2.3 in \cite{gq}. Finally we say that \;$(M, [h])$\; is totally geodesic, if \;$(M, [h])$\; is umbilic and minimal.                                                                                                                                         
\begin{rem}\label{eq:minimal} 
We remark that in the conformally compact Einstein case \;$h_{ y }$\; as in \eqref{eq:uniqdef} has an asymptotic expansion which contains only even powers of \;$ y $, at least up to order \;$n$, see \cite{cg}.
In particular the conformal infinity \;$(M, [h])$\; of  any Poincar\'e-Einstein manifold \;$(X, g^+)$\; is totally geodesic.                                                                                                                                                                                                                                                                                                                                                                                                                                                                                                                                                                                                                                                                                                                                                                                                                                                                                                                                                                                                                                                                                                                                                                                                                                                                                                                                                                                                                                                                                                                                                                                                                                                                                                                                                                                                                                                                                                                                                                                                                                                                                                                                                                                                                                                                                                                                                                                                                                                                                                                                                                                                                                                                                                                                                                                                                                                                                                                                                                                                                                                                          \end{rem}
\begin{rem}
As every \;$2$-dimensional Riemannian manifold is locally conformally flat, we will say {\em locally flat conformal infinity of a Poincar\'e-Einstein manifold}  to mean just the conformal infinity of a Poincar\'e-Einstein manifold when \;$n=2$\; or which is furthermore locally flat, when \;$n>2$. 
\end{rem}

\noindent
On the other hand to any asymptotically hyperbolic manifold \;$(X, g^+)$\; with conformal  infinity \;$(M, [h])$  
Graham-Zworsky\cite{gz}\; have attached a family of scattering operators \;$S(s)$, 
which is a meromorphic family of pseudo-differential operators on \;$M$\; defined on \;$\mathbb{C}$,  
by considering Dirichlet-to-Neumann operators for the scattering problem for \;$(X, \;g^+)$\; and a meromorphic continuation argument. 
Indeed it follows from  \cite{gz} and \cite{mazmel} 
that for every \;$f\in C^{\infty}(M)$ and for every \;$s\in \mathbb{C}$\; such that \;$Re(s)>\frac{n}{2}$\; and \;$s(n-s)$\; is not an \;$L^2$-eigenvalue of \;$-\D_{g^+}$ the following generalized eigenvalue problem
\begin{equation}\label{eq:geneig}
 -\D_{g^+}u-s(n-s)u=0\,\;\;\text{ in }\;\;X
\end{equation}
has a solution of the form
$$
u=F y ^{n-s}+G y ^s, \;\;\; F,\;G\in C^{\infty}(\ov X), \;\;\;F|_{ y =0}=f,
$$
where \;$ y $\; is given by  \eqref{eq:uniqdef} and for those values of \;$s$\; the scattering operator \;$S(s)$\; on \;$M$\; is defined  as
\begin{equation}\label{eq:dtnscat}
S(s)f=G|_{M}.
\end{equation}
Furthermore,  using a meromorphic continuation argument, Graham-Zworsky\cite{gz} extend \;$S(s)$\; defined by \eqref{eq:dtnscat} to a meromorphic family of pseudo-differential operators on \;$M$\; defined on all \;$\mathbb{C}$\; and still denoted by \;$S(s)$\; with only a discrete set of poles including the trivial ones \;$s=\frac{n}{2}, \frac{n}{2}+1, \cdots, $\; which are simple poles of finite rank, and possibly some others corresponding to the \;$L^2$-eigenvalues of \;$-\D_{g^+}$. Using the regular part of the scattering operators \;$S(s)$, to any \;$\gamma = s-\frac{n}{2}\in (0, 1)$\; such that 
\begin{equation*}
\left(\frac{n}{2}\right)^2-\gamma^2<\l_1(-\D_{g^+})
\end{equation*}
with \;$\l_1(-\D_{g^+})$\: denoting the first eigenvalue of \;$-\D_{g_+}$,\; Chang-Gonzalez\cite{cg} have associated the following fractional order pseudo-differential operators,
referred to as fractional conformal Laplacians or fractional Paneitz operators
\begin{equation}\label{eq:fracopd}
P^{\gamma}[g^+, \;h]:=-d_\gamma S\left(\frac{n}{2}+\gamma\right),
\end{equation}
where \;$d_\gamma$\; is a positive constant depending only on \;$\gamma$\; and chosen  such that the principal symbol of \;$P^{\gamma}[g^+,  h]$\; is exactly the same as the one of the fractional Laplacian \;$(-\D_{h})^{\gamma}$, when
\begin{equation*}
X=\R^{n+1}_{+}, \;M=\R^{n},\;h=g_{\R^{n}} \;\;\text{ and }\;\; g^{+}=g_{\mathbb{H}^{n+1}}.
\end{equation*}
When there is no possible confusion with the  metric \;$g^+$, we just use the simple notation
\begin{equation*}
P^{\gamma}_{h}:=P^{\gamma}[g^+,  \;h].
\end{equation*}
Similarly to the other well studied conformally covariant differential operators Chang-Gonzalez\cite{cg} associate to each \;$P^{\gamma}_{h}$\; the curvature quantity
\begin{equation*}
 Q^{\gamma}_{h}:=P^{\gamma}_{h}(1).
\end{equation*}
The functions \;$Q^\gamma_{h}$\; are referred to as fractional scalar curvatures, fractional \;$Q$-curvatures or simply \;$Q^\gamma$-curvatures. 
Of particular importance to conformal geometry is the covariance property 
\begin{equation}\label{eq:confinv}
P^{\gamma}_{h_u}(v)=v^{-\frac{n+2\gamma}{n-2\gamma}}P^{\gamma}_{h}(uv)
\;\;\text{ for }\;\; h_v=v^{\frac{4}{n-2\gamma}} h \;\;\text{ and }\;\; 0<v\in C^{\infty}(M),
\end{equation}
verified by \;$P^{\gamma}_{h}$, see \cite{cg} or Subsection 3.2 in \cite{martndia3}. 
As for the classical scalar curvature, for the \;$Q^\gamma$-curvature Gonzalez-Qing\cite{gq} have introduced 
the fractional \;$\gamma$-Yamabe problem which asks for conformal metrics of constant fractional scalar curvature \;$Q^{\gamma}$. Moreover,  for an asymptotically hyperbolic manifold \;$(X, g^+)$\; with conformal infinity \;$(M, [h])$ being minimal in case \;$\gamma\in (\frac{1}{2},1)$, Chang-Gonzalez\cite{cg} showed the equivalence between the Dirichlet-to-Neumann operators of the scattering problem \eqref{eq:geneig} and the ones of some uniformly degenerate elliptic boundary value problems defined on \;$(X, g)$, which coincide with the extension problem of Caffarelli-Silvestre\cite{cafsyl} when 
\begin{equation*}
(X, g^{+},  h)=(\mathbb{H}^{n+1}, g_{\mathbb{H}^{n+1}}, g_{\R^n}),
\end{equation*} 
and hence 
\begin{equation*}
(-\D_{g_{\R^n}})^{\gamma}=P^{\gamma}[g_{\mathbb{H}^{n+1}}, g_{\R^n}].
\end{equation*}
\vspace{6pt}

\noindent
The latter established relation allows Gonzalez-Qing\cite{gq} to derive a Hopf-type maximum principle 
by taking inspiration from the work of Cabre-Sire\cite{cabsir}, 
which deals with the Euclidean half-space, see Theorem 3.5 and Corollary 3.6 in \cite{gq}. 
Clearly the latter Hopf-type maximum principle of Gonzalez-Qing\cite{gq} opens the door to variational arguments for existence for the \;$\gamma$-Yamabe problem, 
as explored by Gonzalez-Qing\cite{gq}, Gonzalez-Wang\cite{gw} and Kim-Musso-Wei\cite{kmw1}. 
\vspace{6pt}

\noindent
In terms of geometric differential equations the fractional Yamabe problem is equivalent to finding a positive smooth  solution to the semi-linear pseudo-differential equation with critical Sobolev nonlinearity
\begin{equation}\label{eq:fracyam}
P^{\gamma}_{h}u=cu^{\frac{n+2\gamma}{n-2\gamma}}\;\;\text{ on }\;\;M
\end{equation}
for some constant \;$c$. The non-local equation \eqref{eq:fracyam} has a variational structure and  thanks to the regularity theory for uniformly degenerate elliptic boundary value problems (see \cite {cabsir}, \cite{fabkenser}, \cite{gq}), the above cited local interpretation of \;$P^\gamma_h$\; of Chang-Gonzalez\cite{cg} and the Hope-type maximum principle of Gonzalez-Qing\cite{gq}, we have that positive smooth solutions to \eqref{eq:fracyam} can be found by looking at critical points of the following fractional Yamabe functional 
\begin{equation}\label{eq:fracfunc}
\mathcal{E}_{h}^\gamma(u):=\frac{\langle u, u\rangle_{P^{\gamma}_{ h}}}{\left(\int_M u^{\frac{2n}{n-2\gamma}}dV_{ h}\right)^{\frac{n-2\gamma}{n}}}, \;\; u\in W^{\gamma, 2}_+(M,  h):=\{v\in W^{\gamma, 2}(M, h) \; : \;v\geq 0, \;v\not \equiv 0\}, 
\end{equation}
where \;$W^{\gamma, 2}(M,  h)$\; denotes the usual fractional Sobolev space on \;$M$\;with respect to the Riemannian metric \;$h$
and 
\begin{equation}\label{fracscal}
\langle u, u\rangle_{P^{\gamma}_{ h}}=\langle P^{\gamma}_{ h}u, u\rangle_{L^{2}_{h}(M)},
\end{equation} 
with \;$L^{2}_{h}(M)$\; denoting the usual $L^2$-space on \;$M$\; with respect to \;$h$\; and $\langle \cdot, \cdot\rangle_{L^{2}_{h}(M)}$\; 
denoting the scalar product on \;$L^{2}_{h}(M)$.  
For more informations see  \cite{dinpalval}, \cite{gt} and \cite{s}. 
\vspace{6pt}

\noindent
However, as for the classical Yamabe problem and for the same reasons, 
the variational analysis of \;$\mathcal{E}_{h}^\gamma$\; has a local regime, 
namely the situation where the local geometry can be used to ensure a solution (even a minimizer), 
and a global one, where the local geometry cannot be used to find a solution and just a global one, usually called {\em mass}, 
can be used to apply the Aubin-Schoen's minimizing technique. 
We refer to the introduction of \cite{nss} for a precise definition of local and global regimes. 
\vspace{6pt}

\noindent 
Furthermore, still as for the classical Yamabe problem, there is a natural invariant called the \,$\gamma$-Yamabe invariant of \;$(M, [h])$, 
denoted by $\mathcal{Y}^{\gamma}(M, [h])$\; and defined by the formula
\begin{equation}\label{eq:gaminv}
\mathcal{Y}^{\gamma}(M, [h]):=\inf_{u\in W^{\gamma, 2}_+(M)}\frac{\langle u, u\rangle_{P^{\gamma}_{h}}}{(\int_{M}u^{\frac{2n}{n-2\gamma}}dV_{h})^{\frac{n-2\gamma}{2}}}.
\end{equation} 
From the work of Gonzalez-Qing\cite{gq} it is known that \;$\mathcal{Y}^{\gamma}(M, [h])$\; satisfies the  rigidity estimate
\begin{equation*}
\mathcal{Y}^{\gamma}(M, [h])\leq \mathcal{Y}_\gamma(\mathbb{S}^n):= \mathcal{Y}^{\gamma}(S^n, [g_{S^n}]),
\end{equation*}
provided \;$(M,[h])$\; is minimal in case \;$\gamma\in (\frac{1}{2},1)$, where \;$(S^n, [g_{S^n}])$\; is the conformal infinity of the Poincar\'e ball model of the hyperbolic space.
Moreover, as mentioned in the abstract, the global situation clearly includes the case of a locally flat conformal infinity of a Poincar\'e-Einstein manifold, 
as observed by Kim-Musso-Wei\cite{kmw1}, 
while the existence results of Gonzalez-Qing\cite{gq}, Gonzalez-Wang\cite{gw} and some part of the existence results in Kim-Musso-Wei\cite{kmw1} 
deal with situations which clearly belong to the local regime. 
Moreover in the global situation, as for the classical Yamabe problem, to run the minimizing technique of Aubin\cite{au}-Schoen\cite{sc} 
one needs an analogue of the positive mass theorem of Schoen-Yau\cite{sy1}, which is not known to hold, see Conjecture 1.6 in Kim-Musso-Wei\cite{kmw1}. 

\vspace{6pt}

\noindent
We would like to point out that in case of a locally flat conformal infinity of a Poincar\'e-Einstein manifold, using geometric flow techniques, 
in \cite{dsv} solvability of \eqref{eq:fracyam} is proved under the extra assumption of positivity of the \textit{classical} Yamabe invariant. 
There are also works related to the issue of compactness of \eqref{eq:fracyam}, see \cite{kmw} and \cite{kmw2}, and on the singular fractional Yamabe problem, see \cite{acdfgw}, \cite{adgw}, \cite{ddgw}, \cite{gms}. Finally we refer the reader to the survey by Gonzalez\cite{gonz1} for more informations about the fractional Yamabe problem. 
\vspace{6pt}

\noindent
Our main goal in this work is to show that  with the analytical results provided by the works of 
Caffarelli-Silvestre\cite{cafsyl}, Chang-Gonzalez\cite{cg}, Gonzalez-Qing\cite{gq} and Mayer-Ndiaye\cite{martndia3} 
at hand
we can perform a variational argument for existence for the fractional Yamabe problem for a locally flat conformal infinity of a Poincar\'e-Einstein manifold, 
bypassing the lack of knowledge of a fractional analogue of the positive mass theorem of Schoen-Yau\cite{sc} 
and a positivity assumption on the classical Yamabe invariant. 
Indeed, using a suitable scheme of the algebraic topological argument, also called barycenter technique of Bahri-Coron\cite{bc}, 
as implemented in our previous work \cite{martndia2}, we show the following existence theorem.
\begin{thm}\label{eq:thm}
Let \;$n\geq 2$\; be a positive integer, \;$(X^{n+1}, \; g^+)$\; be a Poincar\'e-Einstein manifold with conformal infinity \;$(M^{n}, \;[h])$, \;$\gamma \in(0, 1)$, and \;$\left(\frac{n}{2}\right)^2-\gamma^2<\l_1(-\D_{g^+})$. Assuming that either\;$n=2$\; or \;$n\ge 3$\; and \;$(M, [h])$\; is locally flat, then \;$(M, [h])$\; carries a Riemannian metric of constant \;$Q^{\gamma}$-curvature. 
\end{thm}

\begin{rem}\label{eq:half}
We point out that, as observed by Gonzalez-Qing\cite{gq},  the \;$\frac{1}{2}$-Yamabe problem for Poincar\'e-Einstein manifolds is equivalent to the Riemann mapping problem of Cherrier\cite{cherrier}-Escobar\cite{es2}, which is completely solved after the series of works of \cite{almaraz1}, \cite{chen}, \cite{es2}, \cite{marques}, \cite{martndia2}, \cite{ould1}.
\end{rem}

\noindent
As already mentioned, to prove Theorem \ref{eq:thm} we use variational arguments by applying a suitable scheme of the Barycenter Technique of Bahri-Coron\cite{bc}. 
Indeed, exploiting that the conformal infinity is locally flat, the ambient space is Poincar\'e-Einstein, the conformal covariance property \eqref{eq:confinv}, 
the works of Caffarelli-Silvestre\cite{cafsyl} and Chang-Gonzalez\cite{cg}, the Hopf-type maximum principle of Gonzalez-Qing\cite{gq} 
and the standard bubbles attached to the related optimal trace Sobolev inequality, 
we define some bubbles  and show that they can be used to run a suitable scheme of the barycenter technique of Bahri-Coron\cite{bc} for existence, 
which among others has been used in the works  \cite{gam1}, \cite{gam2}, \cite{ould1} and \cite{nss}. 
We give below a brief discussion of the main ideas behind the argument and refer to the introduction of our paper\cite{martndia2} for a more geometric description
as well as to \cite{gam2} for a detailed and concise exposition. 
\vspace{6pt}
 
\noindent
The barycenter technique of Bahri-Coron\cite{bc} is an argument by contradiction, thus we assume the problem has no solution.
Then, denoting for $1\leq p\in \N$ the limiting energy of $p$-many \textit{non collapsing} bubbles, see \eqref{dfbua}, by 
$$
l_{p}=p^{\frac{2\gamma}{n}}\mathcal{Y}_{\gamma}(S^{n}),  
$$
putting $L_{0}=\emptyset$ and considering for some $\epsilon>0$ the sublevels
$$
L_p:=\{u:  \mathcal{E}_h^\gamma[u]\leq l_{p}+\epsilon\},
$$
on one hand we construct recursively singular chains   $X_{p}$ in   $L_{p}$, 
which generate non zero classes in the relevant $\mathbb{Z}_{2}$-homologies of the topological pairs $(L_{p},L_{p-1})$, 
precisely
\begin{equation}\label{nzfi}
0\neq [X_{p}] \in H_{np+p-1}(L_{p},L_{p-1},\mathbb{Z}_{2}),
\end{equation}
as follows.  
The starting point is the existence of  $X_1$ and non triviality of 
$$[X_{1}]\in H_{n}(L_{1},L_{0},\mathbb{Z}_{2})=H_{n}(L_{1},\mathbb{Z}_{2}),$$ 
which follow from
$H_n(M, \mathbb{Z}_2)\neq 0$,  
embedding $M$ into   $L_1$ via bubbling 
$$
M\longrightarrow L_{1}:a\longrightarrow v_{a,\lambda}
$$ 
and, that based on the quantization phenomenon, which $\mathcal{E}_h^\gamma $ enjoys,
$M$ survives via the deformation Lemma \ref{eq:classicdeform} and selection map \eqref{eq:select} topologically in $L_1$, see Lemma \ref{eq:nontrivialf1}.  
We then start \textit{piling up masses}  $v_{a,\l} $ over $X_1$, 
thereby  iteratively moving from the level   
$l_{p}$   
to the level    
$l_{p+1}$.
At each step one constructs a singular chain $X_{p+1}$ with a non zero class  $[X_{p+1}]$,    
which reads  $$(1-t)u + tv_{a,\l},u\in X_p,t\in  [0,1],$$ 
see Lemma \ref{eq:nontrivialrecursive}. Indeed, denoting by  $\delta_{a}$ for $a\in M$ the Dirac measure at $a$ and recalling 
the space of formal barycenter of $M$,  defined as 
 \begin{equation}\label{dfbp}
B_{p}(M)=\{\sum_{i=1}^{p}\alpha_i\d_{a_i}\;:\;a_i\in M, \;\alpha_i\geq 0,\;\; i=1,\cdots, p,\;\,\sum_{i=1}^{p}\alpha_i=1\},\;
\;B_0(M)=\emptyset,
\end{equation}
the set $B_{p+1}(M)$ as a cone over $B_p(M)$ with top $M$ survives as a non trivial cone in  $(L_{p+1},L_p )$, 
when embedding  $B_{p+1}(M)$ into   $L_{p+1}$ via $(p+1)$-convex combinations of the bubbles $v_{a, \l}$.
Again the latter survival is based on the quantization phenomenon, which $\mathcal{E}_h^\gamma $ enjoys, 
via the deformation Lemma \ref{eq:classicdeform} and the selection map \eqref{eq:select}. 
Since for all $q\geq 1$ we have the existence of  
\begin{equation}\label{orientation_classes} 
0\neq w_q\in H_{nq+q-1}(B_{q}(M), B_{q-1}(M),\mathbb{Z}_{2}),
\end{equation}
see \cite{bc}, we then obtain $[X_{p+1}]\neq 0$ for some $X_{p+1}$, which is the image of a representative of $w_{p+1}\neq 0$,
and \eqref{nzfi} is established.
On the other hand, 
because of the strong interaction phenomenon, 
for  some $p_0$  large we are actually passing from the level    
$l_{p_{0}}+\epsilon$    
to the level    
$l_{p_{0}+1}-\epsilon_{0}$   
for some $\epsilon_0>0$, that is   
\begin{equation*}
X_{p_{0}+1} \;\; \text{ is a chain in } \;\; 
\tilde L_{p_{0+1}}=\{u:\mathcal{E}^\gamma_h \leq l_{p_{0}+1} - \epsilon_0\}
\subset
L_{p_{0}+1}.
\end{equation*}
Moreover, since in absence of solutions and due to the quantization phenomenon
the Palais-Smale condition holds on the sets
$\{u:  l_{p}+\varepsilon<\mathcal{E}_h^\gamma <l_{p+1}-\varepsilon\}$ 
for all $p\geq 1$ and $0<\varepsilon \ll 1$,
the pair $(\tilde L_{p_{0}+1},L_{p_{0}})$ retracts by deformation onto the pair $({L}_{p_{0}},L_{p_{0}})$
and we conclude 
$$
[X_{p_0+1}]
\in
H_{n(p_{0}+1)+p_{0}}(\tilde L_{p_{0}+1},L_{p_{0}},\mathbb{Z}_{2})
= 
\mymathbb{0}.
$$ 
In particular $[X_{p_{0}+1}]=0$ in contradiction  to \eqref{nzfi} for $p=p_{0}+1$
and so a solution must exist.

\begin{rem}\label{rem_limiting_flatness}
Clearly the Poincar\'e-Einstein structure reflects 
the flatness of the conformal infinity \;$(M,[h])$\; into the interior of \;$X$, 
namely, as observed by Kim-Musso-Wei\cite{kmw1}, 
the metric \;$g$\; on \;$X$\; takes locally the form \;$g=\delta+O(|y|^{n})$, 
i.e. $g$ is flat to order \;$n$. 
However, since we base the calculation of the fractional Yamabe energy of a bubble on comparison via maximum principle, 
the appearance of a logarithmic term in the construction of a suitable barrier solution, cf. \eqref{appearance_of_the_log},
in case $g=\delta+O(|y|^{n})$ highly suggests the limitation of our argument to the latter order of flatness.
\end{rem}

\noindent
We would like to make some comments about the application of the barycenter technique of Bahri-Coron\cite{bc} in our situation  
and in the case of the classical Yamabe problem for locally conformally flat closed Riemannian manifolds, as studied by Bahri\cite{bah1}. 
The latter situations are counterpart to each other, however the nonlocal aspect of our situation creates an additional difficulty, 
that is not in its counterpart for the classical Yamabe problem, which is of local nature. 
In fact, even if both problems are conformally invariant and after a conformal change we are locally in the corresponding model space of singularity 
- truly in the classical case and up to a critical, but handleable lower order term in the fractional scenario, cf. Remark \ref{rem_limiting_flatness} - 
we have that the lower order term, i.e. the scalar curvature, {\em vanishes} locally for the classical Yamabe problem, 
because of its  local nature, while for the fractional Yamabe problem that does not necessarily imply that the \;$Q^{\gamma}$-curvature vanishes locally, 
because of the nonlocal aspect of the problem. 
This is the source of the difficulty we mentioned before, which is similar to the one  
encountered by Bahri-Brezis\cite{BB} and in a different framework Brendle\cite{bre2}. 
To overcome the latter issue, we use the works of Caffarelli-Silvestre\cite{cafsyl}, Chang-Gonzalez\cite{cg} and Gonzalez-Qing\cite{gq} 
to reduce ourselves to a local situation. 
Having done that, we then encounter the problem of not having an explicit knowledge of the standard bubble corresponding to the reduced local situation, 
and clearly such an explicit knowledge plays an important role in the corresponding situation of the classical Yamabe problem. 
To deal with this lack of explicit knowledge of the standard bubbles corresponding to the reduced local situation on the \;$1$-dimension augmented half space, 
we observe that its integral representation  given in Caffarelli-Silvestre\cite{cafsyl} 
can be interpreted as a suitable interaction of standard bubbles on the boundary of the latter augmented half space with different points and scales of concentration.
This interpretation provides the required estimates of the latter argument, which one gets for free from an explicit knowledge, and this is made rigorous in this work, see Lemma \ref{lem_standard interaction on Rn} and Corollary \ref{cor_sharp_estimates}. 
We point out that the role of interaction of bubbles in the existence mechanism of Yamabe type problems has been observed for the first time by Bahri-Coron\cite{bc}. 
The idea behind is, that interaction pushes the energy down from the expected critical value at infinity for multiple bubbles. 
This has been successfully used in the study of other Yamabe type problems (see \cite{BB}, \cite{gam1}, \cite{gam2}, \cite{martndia2}, \cite{ould1}, \cite{nss})
and in the study of Yamabe flow (see \cite{bre1}, \cite{bre2}).

%

\noindent
\begin{rem}
We would like to emphasize the analogy between our function \;$\mathcal{M}_{\gamma}$\; given by Definition \ref{fracmass} 
and the notion of mass appearing in the context of the classical Yamabe problem. 
Moreover we point out that  our work in \cite{martndia3} answers the first part of the Conjecture of Kim-Musso-Wei\cite{kmw1} 
about the structure of the Green's function, see Theorem 1.4 in \cite{martndia3} and gives rise to the definition of $\mathcal{M}_{\gamma}$.  
We add that under the global assumption
$M_{\gamma}>0$	
our sharp estimate in Lemma \ref{sharpenergy} gives directly the existence of a fractional Yamabe minimizer.
\end{rem}

\noindent
The structure of the paper is as follows. In Section \ref{notation_and_prelimiaries} we recall the standard bubbles of the variational problem, 
give some preliminaries and fix notations. 
In Section \ref{bubinter} we analyse the standard bubbles on \;$\R^{n+1}_{+}$, 
introduce the relevant Schoen's and Projective bubbles in the curved scenario and prove some interaction estimates of the latter. 
In Section \ref{pceinsten} we establish  sharp \;$L^{\infty}$-estimates for the difference between these bubbles 
and derive a sharp selfaction estimate for the Projective one. 
In Section \ref{varagtop} we present the variational and algebraic topological argument to prove Theorem \ref{eq:thm}. 
It is divided into two subsections.  
In Subsection \ref{var} we present a variational principle 
which extends the classical variational principle by taking into account the non-compactness phenomena via our Projective bubbles. 
In Subsection \ref{agtop}, we present the barycenter technique or algebraic topological argument for existence.  
Finally in Section \ref{appendix} we collect the proofs of some technical lemmas and estimates.

\vspace{6pt}

\noindent
{
\begin{center}
{\bf Acknowledgements}
\end{center}
The authors worked on this project when they were visiting the department of Mathematics of the University of Ulm in Germany. Parts of this paper were written when the authors were visiting the Mathematical Institute of Oberwolfach in Germany as \textit{Research in Pairs} and the 
Institut des Hautes \'Etudes Scientifiques in Paris.  We are very grateful to all these institutions for their kind hospitality.
}

\section{Preliminaries and notations}\label{notation_and_prelimiaries} 
In this section we give some preliminaries and fix notations. 
We start with the standard bubbles.

\vspace{6pt}

\noindent
For \;$a\in \R^n$\; and \;$\l>0$\; we define the standard bubbles on \;$\R^n$\; as 
\begin{equation*}
\delta_{a, \l}(x)=\left(\frac{\l}{1+\l^2|x-a|^2}\right)^{\frac{n-2\gamma}{2}}, \;\;x\in \R^n.
\end{equation*}
They are solutions of the  pseudo-differential equation
\begin{equation}\label{eq:yambaf}
(-\Delta_{\R^n})^{\gamma}\delta_{a, \l}=c_{n, \gamma}\delta_{a, \l}^\frac{n+2\gamma}{n-2\gamma}\;\;\text{on}\;\; \R^n,
\end{equation}
where \;$c_{n, \gamma}$\; is a positive constant depending only on \;$n$\; and \;$\gamma$. 
From \;$\delta_{a, \l}$\; we then define \;$\hat \d_{a, \l}$\; via
\begin{equation}\begin{split}\label{equation_for_delta_hat}
\begin{cases}
div(y^{1-2\gamma}\nabla \hat \delta_{a,\lambda})=0\;\; \text{ in} \;\; \R^{n+1}_+ \\
\hat \delta_{a,\lambda}=\delta_{a, \l} \;\; \text{on}\;\;\R^{n}
\end{cases}
\end{split}\end{equation}
and the quantities
\begin{equation*}
c_{n, 1}^{\gamma}=\int_{\R^{n}}\left(\frac{1}{1+|x|^2}\right)^{n}dx 
\;\;\text{ and }\;\; 
c_{n, 2}^{\gamma}=d_\gamma^*\int_{\R^{n+1}_+}y^{1-2\gamma}\left|\n \hat \d_{0, 1}(y,x)\right|^2dxdy, 
\end{equation*}
where
\begin{equation}\label{dsgamma}
d_{\gamma}^*=\frac{d_\gamma}{2\gamma}, 
\end{equation}
cf. \eqref{eq:fracopd}, which relate via
\;$
-d_{\gamma}^{*}\lim_{y\to 0} (\partial_{y} y^{1-2\gamma}\partial_{y}\hat \delta_{a,\lambda})
= 
c_{n,\gamma}\delta_{a,\lambda}^{\frac{n+2\gamma}{n-2\gamma}}
$\;
and 
\begin{equation*}
\begin{split}
c_{n,1}^{\gamma}
= &
\int_{\R^{n}}\delta_{a,\lambda}^{\frac{2n}{n-2\gamma}}dx
=
\frac{1}{c_{n,\gamma}}\int_{\R^{n}}(-\Delta_{\R^{n}})^{\gamma}\delta_{a,\lambda}\delta_{a,\lambda}dx
=
-\frac{d^{*}_{\gamma}}{c_{n,\gamma}}
\lim_{y\to 0} \int_{\R^{n+1}_{+}} \partial_{y}(y^{1-2\gamma}\partial_{y}\hat \delta_{a,\lambda})\delta_{a,\lambda} dxdy \\
= & \;
\frac{d^{*}_{\gamma}}{c_{n,\gamma}}
\int_{\R^{n+1}_{+}}y^{1-2\gamma}\vert \nabla \hat \delta_{a,\lambda} \vert^{2}dxdy
=
\frac{c_{n,2}^{\gamma}}{c_{n,\gamma}}.
\end{split}
\end{equation*}
The standard bubbles \;$\delta_{a,\lambda}$\; after stereographic projection also minimize 
\begin{equation}\label{eq:yamabehsphere}
\mathcal{Y}_\gamma(S^n)=\mathcal{Y}^\gamma(S^n, [g_{S^n}])
\end{equation}
with \;$\mathcal{Y}^\gamma(S^n, [g_{S^n}])$\; as in \eqref{eq:gaminv}, and therefore
\begin{equation}\label{eq:relationcy}
 c_{n, 2}^{\gamma}=c_{n, \gamma}c_{n, 1}^{\gamma}\;\;\;\;\;\text{and}\;\;\;\;\;\mathcal{Y}_{\gamma}(S^n)=\frac{c_{n, 2}^{\gamma}}{(c_{n, 1}^{\gamma})^{\frac{n-2\gamma}{n}}}.
\end{equation} 
See also Section \ref{bubinter} and we refer to our previous work \cite{martndia3} for details. Furthermore  we set
\begin{equation}\label{defc3}
c_{n, 3 }^{\gamma}=\int_{\R^{n}}\left(\frac{1}{1+|x|^2}\right)^{\frac{n+2\gamma}{2}}dx,\;\;\;
c_{n, 4}^{\gamma}=
c_{n, \gamma}c_{n, 3}^{\gamma},
\;\; \;
c_{n, \gamma}^*=d_\gamma^*c_{n, 3}^\gamma.
\end{equation}
We remark that by \cite{cafsyl}
\begin{equation}\label{dkconv}
\begin{split}
\hat \delta_{0,\lambda}(y,x)  
= & \; 
[K(\cdot , y)*\delta_{0, \l}(x)]
=
p_{n, \gamma} 
\int_{\R^n} \frac{y^{2\gamma}}{(y^{2}+\vert x-\xi\vert^{2})^{\frac{n+2\gamma}{2}}}
\left(\frac{\lambda}{1+\lambda^{2}\vert \xi\vert^{2}}\right)^{\frac{n-2\gamma}{2}}d\xi \\
= & \; 
p_{n, \gamma}
y^{-\frac{n-2\gamma}{2}}\int_{\R} \left(\frac{y^{-1}}{1+y^{-2}\vert x-\xi\vert^{2}}\right)^{\frac{n+2\gamma}{2}}
\left(\frac{\lambda}{1+\lambda^{2}\vert \xi\vert^{2}}\right)^{\frac{n-2\gamma}{2}} d\xi  \\
= & \;
 p_{n, \gamma}
y^{-\frac{n-2\gamma}{2}}
\int_{\R^n} \delta_{x,y^{-1}}^{\frac{n+2\gamma}{n-2\gamma}}(\xi) \;\delta_{0,\lambda}(\xi)d\xi
\end{split}\end{equation}
with \;$K$\; being the Poisson kernel at the origin \;$0\in \R^{n+1}$\; of the operator
\begin{equation}\label{Extension_Operator_flat}
D=-div(y^{1-2\gamma}\nabla (\,\cdot\,))
\end{equation}
and given by
\begin{equation}\label{pngamma}
K(y,x)=K^{\gamma}(y,x)=p_{n, \gamma}\frac{y^{2\gamma}}{(|x|^2+|y|^2)^{\frac{n+2\gamma}{2}}},
\;\;\;
p_{n,\gamma}=\frac{1}{c_{n,3}^{\gamma}},
\end{equation}
see \eqref{defc3}. We also have $$
\hat \delta_{a,\lambda}(\cdot)=\hat \delta_{0,\lambda}(\cdot-a)\;\;\;\forall a\in \R^n.
$$
\vspace{6pt}
 
\noindent
Following \cite{gq}, for an asymptotically hyperbolic manifold \;$(X, g^+)$\;  of dimension $n+1\geq 3$ and with conformal infinity \;$(M, [h])$ we introduce
\begin{equation}\label{eqdg}
\begin{split}
D_{g}U=-div_{g}(y^{1-2\gamma}\nabla_{g}U)+E_{g}U
\end{split}
\end{equation}
with \;$g=y^2g^+$, $y$\; the unique geodesic defining function associated to \;$h$ and 
\begin{equation*}
E_{g} :=y^{\frac{1-2\gamma}{2}}L_{g}y^{\frac{1-2\gamma}{2}} -(\frac{R_{g^{+}}}{c_{n}}+s(n-s))y^{(1-2\gamma)-2},
\end{equation*}
where 
$$L_g=-\Delta_{g}+\frac{R_{g}}{c_{n}},\; c_{n}=\frac{4n}{n-1}$$ denotes the conformal Laplacian of \;$(X, g)$.
This operator realizes via the Dirichlet to Neumann map
\begin{equation}\label{Dirichlet_to_Neumann_map}
u
\longrightarrow 
\left\{
\begin{matrix*}[c]
D_{g}U&=&0  &  \;\text{in}\;  &  X \\
U&=&u &  \;\text{on}\;  &   M\\
\end{matrix*}
\right\}
\longrightarrow
-d_{\gamma}^{*}\lim_{y\to 0}y^{1-2\gamma}\partial_{y}U=f
\end{equation}
the conformal fractional Laplacian \;$P_{h}^{\gamma}u=f$, provided
$$0<\gamma<\frac{1}{2} \; \text{ or, if } \; \frac{1}{2}<\gamma < 1, \; \text{ then additionally } \;  H_{g} \equiv 0.$$ 
We refer to our work \cite{martndia3} for details, where   
the existence and asymptotic behavior of 
the Poisson kernel \;$K_g:=K_{g}^{\gamma}$ \; of \;$D_g$, 
the Green's functions \;$\Gamma_g:=\Gamma_{g}^{\gamma}$\; of \;$D_g$\; under weighted normal boundary condition 
and \;$G_{h}:=G^{\gamma}_{h}$\; of the fractional conformal  Laplacian \;$P^{\gamma}_{h}$\; have been studied.
 
The fundamental solutions \;$K_g$, \;$\Gamma_g$\; and \;$G_h$ are defined via
\begin{equation}\label{kg}
\begin{split}
\begin{cases}
D_{g}K_g(\cdot, \xi)=0\;\; \text{ in }\;\; X \;\;\text{ and for all }\;\; \xi\in M\\
\lim_{y\rightarrow 0}K_g(y, x, \xi)=\d_{x}(\xi)\;\;\text{ and for all }\;\; x,\;\xi\in M,
\end{cases}
\end{split}\end{equation}
\begin{equation}\label{gamg}
\begin{split}
\begin{cases}
D_{g}\Gamma_g(\cdot, \xi)=0\; \text{ in } X \;\;\text{ and for all }\;\; \xi\in M\\
-d_\gamma^*\lim_{y\rightarrow 0}y^{1-2\gamma}\partial_y\Gamma_g(y, x, \xi)=\d_{x}(\xi)\;\;\text{ and for all }\;\; x,\xi\in M,
\end{cases}
\end{split}
\end{equation}
\begin{equation}\label{gh}
P^\gamma_h G_{h}(x, \xi)=\d_{x}(\xi),\;\;x,\xi\in M,
\end{equation}

where \;$d_\gamma^*$ is as in \eqref{dsgamma}, and are linked by  $$\Gamma_g(z, \xi)=\int_M K_g(z, x)G_h(x, \xi)dV_h(x), \;z\in \ov X, \xi \in M.$$ In particular
\begin{equation}\label{Gamma_is_G_on_M}
\lim_{y\to 0} \Gamma_{g}(y,\cdot)=G_{h}(\cdot).
\end{equation}

\bigskip

\noindent
\textbf{Notation:} \; We conclude this section with fixing some notations, used in this paper. 

$\N$\;denotes the set of non negative integers, $\N^*$\; the set of positive integers and  
for $k\in \N^*$, $\R^k$\;stands for the standard $k$-dimensional Euclidean space, 
$\R^k_+$ the open positive half-space of $\R^k$ 
and $\bar \R^k_+$ its closure in $\R^k$. 
For simplicity we let \;$\R_+=\R^1_+$ and $\bar \R_+=\bar \R^1_+$. 
For $r>0$ we denote respectively by
\begin{equation*}
B^{\R^k}_r(0) \;\; \text{ and }\;\; B^{\R^k_+}_r(0)=B_r^{\R^{k}}(0)\cap  \R^k_+ 
\end{equation*}
the open and open upper half ball of \;$\R^k$\; of center \;$0$\; and radius \;$r$, and set
$$B_r=B_r(0)=B_r^{\R^n}(0) \; \text{ and } \; B_r^+=B_r^+(0)=B_r^{\R^{n+1}_+}(0).$$ 

We recall that \;$X=X^{n+1}$ with $n\geq 2$ is a manifold of dimension \;$n+1$\; with boundary \;$M=M^{n}$\; and closure \;$\ov{X}$. 
For any Riemannian metric \;$\bar h$\; defined on \;$M$, $a\in M$\; and \;$r>0$ we use the notation \;$B^{\bar h}_{r}(a)$\; to denote 
the geodesic ball with respect to $\bar h$\; of radius \;$r$\;and center \;$a$.  
We also denote by \;$d_{\bar h}(x,y)$\; the geodesic distance with respect to \;$\bar h$\; 
between  two points \;$x,y\in M$.  
$dV_{\bar h}$\; denotes the Riemannian measure associated to the metric \;$\bar h$\; on \;$M$. 
For \;$a\in M$\; we use the notation \;$\exp^a_{\bar h}$\; to denote the exponential map with respect to \;$\bar h$\; on \;$M$. 

Similarly for any Riemannian metric \;$\bar g$\; defined on \;$\ov{X}$, $a\in M$\; and \;$r>0$ we use the notation \;$B^{\bar g, +}_{r}(a)$\; to denote the geodesic half ball with respect to \;$\bar g$\; of radius \;$r$\;and center \;$a$.  We also denote by \;$d_{\bar g}(x,y)$\; the geodesic distance with respect to \;$\bar g$\; between two points 
\;$x,y\in \ov {X}$.  
 $dV_{\bar g}$\; denotes the Riemannian measure associated to the metric \;$\bar g$\; on \;$\ov{X}$. 
 For \;$a\in M$\; we use the notation \;$\exp_a^{\bar g, +}$\; to denote the exponential map with respect to \;$\bar g$\; on \;$\ov{X}$. 

For \;$p\in \N^*$\; let \;$\sigma_p$\; denote the permutation group of \;$p$\; elements 
and let \;$M^p$\; denote the Cartesian product of \;$p$\; copies of \;$M$.
We define
\begin{equation}\label{M_squred_star}
(M^2)^*:=M^2\setminus Diag(M^2), 
\end{equation}
where \;$Diag(M^2)=\{(a, a): \,a\in M\}$\; is the diagonal of \;$M$. 

For \;$\bar g$, a Riemannian metric defined on \;$\ov{X}$, let 
 \;$C^{\infty}(\bar X, \bar g)$\; denote the space of infinitely differentiable functions on \;$\ov{X}$\; with respect to \;$\bar g$.

Large positive constants are usually denoted by \;$C$\; and the value of \;$C$\; is allowed to vary from formula to formula and also within the same line. Similarly small positive constants are denoted by \;$c$\; and their values may vary from formula to formula and also within the same line. 

$O(1)$ stands for quantities, which are bounded. 
For \;$\epsilon>0$\;  we denote by \;$o_{\epsilon}(1)$\; any quantity,  which tends to \;$0$, as \;$\epsilon\longrightarrow 0$.   
For \;$x\in \R$\; we denote by \;$O(x)$\; and 
\;$o_{\epsilon}(x)$\; respectively \;$|x|O(1)$\; and \; $|x|o_{\epsilon}(1)$.

For a topological space \;$Z$\; let \;$H_{*}(Z)$\; denote the singular homology of \;$Z$\; with \;$\Z_2$\; coefficients and 
for a subspace \;$Y$\; of \;$X$\; let \;$H_*(Z, Y)$\; denote the relative homology. 
For a map \;$f:Z\longrightarrow  Y$\; with \;$Z$\; and \;$Y$\; topological spaces we denote by \;$f_*$\; the induced map in homology.

\section{Bubbles and related interaction estimates}\label{bubinter}
In this section we recall the definition of the standard bubbles on  \;$\R^{n+1}_+$\; and their interpretation as a suitable interaction of standard bubbles on \;$\R^n$. 
Furthermore we establish some new sharp estimates of independent interest for the interaction of the standard bubbles on 
\;$\R^n$\; and use the latter to derive sharp estimates for the standard bubbles on \;$\R^{n+1}_+$. 
Moreover we recall the Schoen's  bubbles associated to the standard bubbles on \;$\R^{n+1}_+$\; and use them to define other bubbles,
called them Projective bubbles,
which \textit{talk} to the local formulation of the problem. Finally we derive sharp interaction estimates for the Projective bubbles.

\subsection{The bubbles on \;$\R^{n+1}_{+}$\; as an interaction of standard ones on \;$\R^{n}$}\label{standarrnp}
In this subsection we deal with the standard bubbles \;$\hat \delta_{a, \l} $\; on \;$\R^{n+1}_+$. They are the natural extension of the fractional 
bubbles \;$\delta_{a, \l}$\; on \;$\R^{n}$\; with respect to \;$D$, cf. \eqref{Extension_Operator_flat}, and are  given by the 
convolution of the Poisson kernel \;$K$\; of \;$D$, cf. \eqref{dkconv}, 
and can be treated as a standard bubble interaction on \;$\R^n$. Indeed \eqref{dkconv} is equivalent to
\begin{equation}\label{bubble_interpretation}
\hat \delta_{a,\lambda}(y,x)
= 
p_{n, \gamma}
y^{-\frac{n-2\gamma}{2}}
\int_{\R^{n}} \delta_{x,y^{-1}}^{\frac{n+2\gamma}{n-2\gamma}}\delta_{a,\lambda}.
\end{equation}

\noindent
This interpretation will be used to derive sharp estimates for \;$\hat \d_{a, \l}$\; relating it to the scale of the Green's function 
\;$\lambda^{\frac{2\gamma-n}{2}}\Gamma(a,\cdot)$, which is necessary for sharp energy estimates. 
Recalling that interaction always relates to suitable scales of the Green's function, see Lemma \ref{lem_standard interaction on Rn} for instance, 
then clearly \eqref{bubble_interpretation} provides a way to achieve that and we will follow this approach.

\bigskip

\noindent
We then need the following Lemma \ref{lem_standard interaction on Rn} on interaction estimates for standard  bubbles on \;$\R^{n}$, 
which are new in the fractional setting, sharp and of independent interest. 
While these estimates are completely analogous to the classical ones, 
e.g. in case of the conformal Laplacian and famously to those in \cite{bah}, 
they do not only serve the purpose of understanding the interaction of different bubbles on \;$\R^{n}$, 
but by virtue of \eqref{bubble_interpretation}, 
as explained above, 
give also rise to sharp estimates on their extensions on \;$\R^{n+1}_{+}$.  
We anticipate that the proof of Lemma \ref{lem_standard interaction on Rn} is postponed to the appendix given by Section \ref{appendix}.
\begin{lem} 
\label{lem_standard interaction on Rn}
For \;$c_{n, 3}^{\gamma}$\; given by \eqref{defc3}, \;$i\neq j$\; and \;$\varepsilon_{i,j}=(\frac{\lambda_{i}}{\lambda_{j}}+\frac{\lambda_{j}}{\lambda_{
i
}}+\lambda_{i}\lambda_{j}\vert a_{i}-a_{j}\vert^{2})^{\frac{2\gamma-n}{2}}
$\;  there holds
\begin{enumerate}[label=(\roman*)]
 \item \quad 
$\int_{\R^{n}}\delta_{a_{i}, \lambda_{i}}^{\frac{n+2\gamma}{n-2\gamma}}
\delta_{a_{j}, \lambda_{j}}
=
c_{n, 3}^{\gamma}\varepsilon_{i,j}
+
O(\min(\frac{\lambda_{i}}{\lambda_{j}}, \frac{\lambda_{j}}{\lambda_{i}})^{\gamma}
\varepsilon_{i,j}^{\frac{n}{n-2\gamma}})$
 \item \quad 
$
\lambda_{k}\partial_{\lambda_{k}}\int_{\R^{n}}\delta_{a_{i}, \lambda_{i}}^{\frac{
n+2\gamma}{n-2\gamma}}
\delta_{a_{j}, \lambda_{j}}
=
c_{n, 3}^{\gamma}\lambda_{k}\partial_{\lambda_{k}}\varepsilon_{i,j}
+
O(\min(\frac{\lambda_{i}}{\lambda_{j}}, \frac{\lambda_{j}}{\lambda_{i}})^{\gamma}
\varepsilon_{i,j}^{\frac{n}{n-2\gamma}})
$

\item \quad 
$
\nabla_{a_{k}}\int_{\R^{n}}\delta_{a_i, \l_i}^{\frac{n+2\gamma}{n-2\gamma}}\delta_{a_j, \l_j}
= 
c_{n, 3}^{\gamma}\lambda_{i}\lambda_{j}(a_{k}-a_{l})
\varepsilon_{i,j}^{\frac{n+2-2\gamma}{n-2\gamma}}
+
O(\min(\frac{\lambda_{i}}{\lambda_{j}}, \frac{\lambda_{j}}{\lambda_{i}})^{\gamma}
\sqrt{\lambda_{i}\lambda_{j}}\varepsilon_{i,j}^{\frac{n+1}{n-2\gamma}})
$

\item \quad 
$
\nabla^{2}_{a_{k},a_{l}}\int_{\R^{n}}\delta_{a_i, \l_i}^{\frac{n+2\gamma}{n-2\gamma}}
\delta_{a_j, \l_j}
= 
c_{n, 3}^{\gamma}\lambda_{i}\lambda_{j}
(
id-\frac{(n+2-2\gamma)\lambda_{i}\lambda_{j}(a_{i}-a_{j}).(a_{i}-a_{j})}{\frac{
\lambda_{i}}{\lambda_{j}}+\frac{\lambda_{j}}{\lambda_{i}}+\lambda_{i}\lambda_{j}
\vert a_{i}-a_{j} \vert^{2}}
)
\varepsilon_{i,j}^{\frac{n+2-2\gamma}{2}}
\\ _{}\quad\quad\quad\quad\quad\quad\quad\quad\quad\quad\,\quad \,\,  +
O(\min(\frac{\lambda_{i}}{\lambda_{j}}, \frac{\lambda_{j}}{\lambda_{i}})^{\gamma}
\lambda_{i}\lambda_{j}\varepsilon_{i,j}^{\frac{n+2}{n-2\gamma}})
$
\end{enumerate}
for \; $k,l\in \{ i,j \}$\; and \;$\{ k,l \} = \{ i,j \}$. 
\end{lem}

\noindent 
Lemma \ref{lem_standard interaction on Rn} has the following impacts on the standard bubbles \;$\hat \d_{a, \lambda}$\; on \;$\R^{n+1}_+$.
\begin{cor}
\label{cor_rough_estimates}
For $a\in \R^n$\; and  \;$r_{a}=\vert z-a \vert=\vert (y,x-a)\vert$\; we have on \;$\R^{n+1}_{+}$
\begin{enumerate}[label=(\roman*)]
 \item \quad
$
\hat \delta_{a,\lambda}(y,x) 
=
O((\frac{\lambda}{1+\lambda^{2}r_{a}^{2}})^{\frac{n-2\gamma}{2}})
$
 \item \quad
$
\partial_{y}\hat 
\delta_{a,\lambda}(y,x)
=
O(\lambda^{\gamma}y^{2\gamma-1}(\frac{\lambda}{1+\lambda^{2}r_{a}^
{ 2 } }
)^{\frac{n}{2}})
$
 \item \quad
$
\nabla_{x}\hat\delta_{a,\lambda}(y,x)
=
O(\sqrt{\lambda}(\frac{\lambda}{1+\lambda^{2}
r_{a}^{2}})^{\frac{n+1-2\gamma}{2}})
$
 \item \quad
$
\nabla^{2}_{x}\hat\delta_{a,\lambda}(y,x)
=
O(\lambda(\frac{\lambda}{1+\lambda^{2}r_{a}^{2}})^{\frac{n+2-2\gamma}{2}}).
$
\end{enumerate}
\end{cor}
\begin{pf} 
Let \;$a_{i}=x,\lambda_{i}=\frac{1}{y},a_{j}=a$ and $\lambda_{j}=\lambda$. Then \eqref{bubble_interpretation} and 
Lemma \ref{lem_standard interaction on Rn} (i) show
\begin{equation*}
\begin{split}
\hat \delta_{a,\lambda}(y,x)
\leq
Cy^{-\frac{n-2\gamma}{2}}\varepsilon_{i,j}
=
O(\frac{y^{-\frac{n-2\gamma}{2}}}{(\frac{1}{\lambda y}+\lambda 
y+\frac{\lambda}{y}\vert x-a\vert^{2})^{\frac{n-2\gamma}{2}}}),
\end{split}
\end{equation*}
whence (i) follows. From \eqref{pngamma}, \eqref{bubble_interpretation} and Lemma \ref{lem_standard interaction on Rn} 
(ii) we infer using \;$y\partial_{y}=-y^{-1}\partial_{y^{-1}}$\; that
\begin{equation*}
\begin{split}
y\partial_{y}&\hat \delta_{a,\lambda}(y,x)
= 
p_{n, \gamma}y\partial_{y}(y^{-\frac{n-2\gamma}{2}}\int_{\R^{n}}
\delta^{\frac{n+2\gamma}{n-2\gamma}}_{x,y^{-1}}\delta_{a,\lambda}) \\ 
= & \;
y\partial_{y}y^{-\frac{n-2\gamma}{2}}\frac{1}{(\frac{1}{\lambda y}+\lambda 
y+\frac{\lambda}{y}\vert x-a \vert^{2})^{\frac{n-2\gamma}{2}}} \\
& +
y^{-\frac{n-2\gamma}{2}}y\partial_{y}\frac{1}{(\frac{1}{\lambda y}+\lambda 
y+\frac{\lambda}{y}\vert x-a \vert^{2})^{\frac{n-2\gamma}{2}}}
+
O(\frac{\min(\frac{1}{\lambda y}, \lambda 
y)^{\gamma}y^{-\frac{n-2\gamma}{2}}}{(\frac{1}{\lambda y}+\lambda 
y+\frac{\lambda}{y}\vert x-a \vert^{2})^{\frac{n}{2}}}) \\
& = \;
-(n-2\gamma)(\frac{\lambda}{1+\lambda^{2}r_{a}^{2}})^{\frac{n-2\gamma}{2}}
[
\frac{\lambda^{2}y^{2}}{1+\lambda^{2}r_{a}^{2}}
+
O((\frac{\min(1,\lambda^{2} y^{2})}{1+\lambda^{2}r_{a}^{2}})^{\gamma})
],
\end{split}
\end{equation*}  
whence (ii) follows. Likewise (iii),(iv) follow from \eqref{bubble_interpretation} and Lemma \ref{lem_standard interaction on Rn} (iii),(iv)
respectively.
\end{pf}

\noindent
\begin{cor}
\label{cor_sharp_estimates}
For \;$a\in \R^n$\; and \;$r_{a}=\vert z-a \vert=\vert (y,x-a)\vert$\; we have on \;$B_{\varepsilon}(a)^{c}\cap \R^{n+1}_{+}$\
\begin{enumerate}[label=(\roman*)]
 \item \quad 
$\hat \delta_{a,\lambda}(y,x)
= 
\frac{1}{\lambda^{\frac{n-2\gamma}{2}}r_{a}^{n-2\gamma}}
+
o_{\frac{1}{1+\lambda^{2\gamma}\varepsilon^{n}}}(\lambda^{\frac{2\gamma-n}{2}})$
 \item \quad
$y\partial_{y}\hat \delta_{a,\lambda}(y,x)
= 
y\partial_{y}\frac{1}{\lambda^{\frac{n-2\gamma}{2}}r_{a}^{n-2\gamma}}
+
o_{\frac{1}{1+\lambda^{2\gamma}\varepsilon^{n}}}(\lambda^{\frac{2\gamma-n}{2}})$
 \item \quad
$x\nabla_{x}\hat \delta_{a,\lambda}(y,x)
= 
x\nabla_{x}\frac{1}{\lambda^{\frac{n-2\gamma}{2}}r_{a}^{n-2\gamma}}
+
o_{\frac{1}{1+\lambda^{2\gamma}\varepsilon^{n}}}(\lambda^{\frac{2\gamma-n}{2}})$
\end{enumerate}
\end{cor}

\begin{pf}
Let $a_{i}=x,\lambda_{i}=\frac{1}{y},a_{j}=0$ and $\lambda_{j}=\lambda$. Then \eqref{pngamma}, \eqref{bubble_interpretation} and 
Lemma \ref{lem_standard interaction on Rn} (i) give
\begin{equation*}
\begin{split}
\hat \delta_{a,\lambda}(y,x)
= &
y^{-\frac{n-2\gamma}{2}}(\varepsilon_{i,j}
+
O(\min(\frac{\lambda_{i}}{\lambda_{j}}, \frac{\lambda_{j}}{\lambda_{i}})^{\gamma}
\varepsilon_{i,j}^{\frac{n}{n-2\gamma}})) \\
= &
\frac{y^{-\frac{n-2\gamma}{2}}}{(\frac{1}{\lambda y}+\lambda 
y+\frac{\lambda}{y}\vert x-a\vert^{2})^{\frac{n-2\gamma}{2}}}
+
O(\frac{1}{(\lambda 
y)^{\gamma}}\frac{y^{-\frac{n-2\gamma}{2}}}{(\frac{1}{\lambda y}+\lambda 
y+\frac{\lambda}{y}\vert x-a\vert^{2})^{\frac{n}{2}}}) \\
= & 
(\frac{\lambda}{1+\lambda^{2}r_{a}^{2})})^{\frac{n-2\gamma}{2}}
+
O(\frac{1}{\lambda^{\frac{n+2\gamma}{2}}r_{a}^{n}})
,
\end{split}
\end{equation*} 
whence (i) follows.
Moreover we find from \eqref{bubble_interpretation} and Lemma \ref{lem_standard interaction on Rn} (ii) as before
\begin{equation*}
\begin{split}
y\partial_{y}\hat \delta_{a,\lambda}(y,x)
& = 
-(n-2\gamma)(\frac{\lambda}{1+\lambda^{2}r_{a}^{2}})^{\frac{n-2\gamma}{2}}
[
\frac{\lambda^{2}y^{2}}{1+\lambda^{2}r_{a}^{2}}
+
O((\frac{\min(1,\lambda^{2} y^{2})}{1+\lambda^{2} r_{a}^{2} 
})^{\gamma})
],
\end{split}
\end{equation*}  
whence 
\begin{equation*}
y\partial_{y}\hat \delta_{a,\lambda}(y,x)
=
y\partial_{y}\frac{1}{\lambda^{\frac{n-2\gamma}{2}}r_{a}^{n-2\gamma}}
+
\frac
{1}
{\lambda^ {\frac {n-2\gamma}{2}} 
}O
(
\frac
{1}
{\lambda^{2}r_{a}^{n-2\gamma+2}}
+
\frac
{1}
{\lambda^{2\gamma}r_{a}^{n}
} 
)
\end{equation*}
and (ii) follows. (iii) follows analogously from \eqref{bubble_interpretation} and Lemma \ref{lem_standard 
interaction on Rn} (iii).
\end{pf}

\subsection{Schoen's bubbles and Projective bubbles}
In this subsection we recall the Schoen's  bubbles associated to the standard bubbles on \;$\R^{n+1}_{+}$\; 
for an asymptotically hyperbolic manifold \;$(X, g^+)$ of dimension \;$n+1$\; with  \;$n\geq 2$\; and minimal conformal infinity \;$(M, [h])$.
Moreover we use them to define another type of bubbles, called them {\em Projective bubbles}, whose construction is motivated by the local interpretation of the problem under study. 
Furthermore  we establish some interaction estimates for our Projective bubbles, which by the way coincide with the  Schoen's  bubbles on the conformal infinity.
\vspace{6pt}

\noindent
First of all, because of \eqref{eq:uniqdef} and minimality of the conformal infinity, we can consider a geodesic defining function \;$y$\; splitting the metric 
\begin{equation*}
g=y^{2}g^{+}, \;\; g=dy^{2}+h_{y}\;\;\text{near}\;\;M\;\;\text{ and } \;\;h=h_{y}\lfloor_{M}
\end{equation*} 
in such a way, that \;$H_{g}=0$.
Moreover,  using the existence of conformal normal coordinates, cf.  \cite{gun}, there exists for every \;$a\in M$\; a conformal factor
\begin{equation}\label{conformal_factor_properties}
0<u_{a}\in C^{\infty}(M)\;\;\text{ satisfying }\;\;\frac{1}{C}\leq u_{a}\leq 
C, \;\;u_{a}(a)=1\;\;\text{ and }\;\;\nabla u_{a}(a)=0,
\end{equation}
inducing a conformal normal coordinate system close to \;$a$\; on \;$M$, in 
particular in normal coordinates with respect to  
\begin{equation}\label{h_a_metric}
h_{a}=u_{a}^{\frac{4}{n-2\gamma}}h
\;\; \text{ with } \;\;
dV_{h_{a}}=u_{a}^{\frac{2n}{n-2\gamma}}dV_{h}
\end{equation}
we have for some small \;$\epsilon>0$, that
\;$
h_{a}=\delta+O(\vert x \vert^{2}), \;\det h_{a}\equiv 1\;\text{ on } \;
B_{\epsilon}^{h_a}(a)
$.
As clarified in Subsection  3.2 of \cite{martndia3} the conformal factor \;$u_{a}\;$ then 
naturally extends onto \;$X$ via 
$$u_{a}=(\frac{y_{a}}{y})^{\frac{n-2\gamma}{2}},$$
where \;$y_{a}$\; close to the 
boundary \;$M$ is the unique geodesic defining function, for which
\begin{equation*}
g_{a}=y_{a}^{2}g^{+}, \;\; g_{a}=dy_{a}^{2}+h_{a, y_a}\;\;\text{near}\;\;M\;\;\text{ with } 
\;\;h_{a}=h_{a, y_a}\lfloor_{M}
\end{equation*} 
and there still holds \;$H_{g_{a}}=0$. Consequently
\begin{equation*}
g_{a}=\delta+O(y+\vert x \vert^{2})\;\;\text{ and }\;\; \det g_{a}=1+O(y^{2}) 
\;\;\text{ in } \;\; B_{\epsilon}^{g_a, +}(a).
\end{equation*}
With respect to \;$h_{a}$-normal coordinates\; $x=x_{a}$\; centered at \;$a$\; we define the standard bubble 
\begin{equation}\begin{split}\label{delta_boundary}
\delta_{a,\lambda}
=
(\frac{\lambda}{1+\lambda^{2}\vert x\vert^{2}})^{\frac{n-2\gamma}{2}}\,\;\;
\text{ on } \;\;B_{4\rho_{0}}^{h_a}(a)
\end{split}\end{equation}
for some small \;$\rho_0>0$, identifying \;$M\ni a \sim 0\in \R^{n}$, and the Schoen's  bubble associated to \;$\d_{a, \l}$\; by 
\begin{equation}\begin{split}\label{bubble_boundary}
\varphi_{a,\lambda,\varepsilon}
=
\eta_{a,\varepsilon}\delta_{a,\lambda}
+
(1-\eta_{a,\varepsilon})\lambda^{\frac{2\gamma-n}{2}}\frac{G_{a}}{g_{n,\gamma}}
 \;\;\text{ on } \;\;M
 ,\;\;
G_{a}=G_{h_{a}}(\cdot,a)
\end{split}\end{equation}
with \;$G_{h_a}$\; as \eqref{gh},  $\eta_{a,\varepsilon}$\; a cut-off function defined in \;$g_{a}$-normal Fermi-coordinates by
\begin{equation}\label{cut_off}
\begin{split}
\eta_{a,\varepsilon}=\eta(\frac{y^{2}+\vert x \vert^{2}}{\varepsilon^{2}}), \;\;
\eta\equiv 1 \,\,\;\text{ on } [0,1],\;\;
\eta\equiv 0\;\; \text{ on } \;\;[2,\infty)\;\;
\text{ and }\;\;
0<\varepsilon\ll 1
\end{split}
\end{equation} 
and \;$g_{n,\gamma}>0$\; 
such that 
\begin{equation}\label{motivation_g_n_gamma}
G_{a}/g_{n,\gamma}=(1+o_{r_{a}}(1))r_{a}^{2\gamma-n}
\end{equation}
cf. \eqref{delta_boundary} and see \cite{martndia3} for an expansion of \;$G_{a}$. 
 
\bigskip

In what follows we will always choose 
$
\varepsilon \sim \lambda^{-\frac{1}{k}} \text{ with }1\ll k\in \N,
$
for instance $\varepsilon =\lambda^{-\frac{1}{k}}$,
and hence may write
$\varphi_{a,\lambda,\varepsilon}=\varphi_{a,\lambda}$ 
by abusing the  notation.  
We define our Projective bubbles on \;$X$\; by
$$\overline{\varphi_{a,\lambda}}^{a}=K_{g_{a}}*\varphi_{a,\lambda},$$ 
cf. \eqref{kg},  i.e. $\varphi_{a,\lambda}$ uniquely solves
\begin{equation}\begin{split}\label{canonical_bubble}
\begin{cases}
D_{g_{a}}\overline{\varphi_{a,\lambda}}^{a}=0 \;\;\text{ in }\;\; X \\
\overline{\varphi_{a,\lambda}}^{a}=\varphi_{a,\lambda} \;\;\text{ on }\;\; M
\end{cases}
\end{split}\end{equation}
and is the canonical extension of \;$\varphi_{a,\lambda}$\; to \;$X$\; with respect to \;$D_{g_a}$\; as in \eqref{eqdg}. 
Using $g_{a}$-normal Fermi-coordinates 
\;$
z=z_{a}=(y ,x_{a})=(y,x)
$,
we consider as a 
second extension \begin{equation}\begin{split}\label{constructed_bubble}
\widetilde{\varphi_{a,\lambda}}^{a}
=
\eta_{a,\varepsilon}\hat 
\delta_{a,\lambda}+(1-\eta_{a,\varepsilon})\lambda^{\frac{2\gamma-n}{2}}\frac{\Gamma_{a}}{g_{n,\gamma}},
\;\;\;\;
\Gamma_{a}=\Gamma_{g_{a}}(\cdot,a),
\end{split}\end{equation}
with \;$\Gamma_{g_a}$\; as in \eqref{gamg},  namely the Schoen's bubble associated to the standard bubble \;$\hat\delta_{a,\lambda}$\; on \;$\R^{n+1}_+$, 
cf. \eqref{eq:yambaf} and \eqref{equation_for_delta_hat}, solving
\begin{equation}\begin{split}\label{delta_hat}
\begin{cases}
D\hat \delta_{a,\lambda}=0\;\;\text{ in }\;\;\R^{n+1}_{+} \\
\hat \delta_{a,\lambda}=\delta_{a,\lambda} \;\; \text{ on } \;\;\R^{n}\\
-d_\gamma^*\lim_{y\rightarrow 0}y^{1-2\gamma}\partial_{y}\hat \delta_{a,\lambda}= 
c_{n, \gamma}\delta_{a,\lambda}^{\frac{n+2\gamma}{n-2\gamma}} \;\;\text{ on } \;\;\R^{n},
\end{cases}
\end{split}\end{equation}
where \;$D$\; is as in \eqref{Extension_Operator_flat}.
Readily \;$\widetilde{\varphi_{a,\lambda}}^{a}$\; enjoys analogous, but weaker 
identities, namely
\begin{equation}\label{properties_of_tilde_phi}
\begin{cases}
D_{g_{a}}\widetilde{\varphi_{a,\lambda}}^{a}=0 \;\;\text{ in }\;\; 
 X \setminus B_{2\varepsilon}^{g_a, +}(a)  \\
\widetilde{\varphi_{a,\lambda}}^{a}=\varphi_{a,\lambda} \;\;\text{ on }\;\;M \\
-d_\gamma^*\lim_{y\rightarrow 0}y^{1-2\gamma}\partial_{y}\widetilde{\varphi_{a,\lambda}}^{a}=c_{n, \gamma}\eta_{a,
\varepsilon}\delta_{a,\lambda}^{\frac{n+2\gamma}{n-2\gamma}}
\;\;\text{ on }\;\;M.
\end{cases}
\end{equation}
In fact, recalling \eqref{constructed_bubble}, 
the first property in \eqref{properties_of_tilde_phi} is due to \eqref{gamg} and \eqref{cut_off}, while the second one is by definition. 
Furthermore, due to 
\;$\partial_{y}\eta_{a,\varepsilon}=\varepsilon^{-2}O(y), \; \gamma\in (0,1)$\; and, 
since \;$\Gamma_{a}$\; extends smoothly to the boundary away from $a$, 
cf. \eqref{Gamma_is_G_on_M}, 
there holds  
\begin{equation*}
\begin{split}
-d_\gamma^*\lim_{y\rightarrow 0}y^{1-2\gamma}\partial_{y}\widetilde{\varphi_{a,\lambda}}^{a} 
= \; &
-d_\gamma^*\lim_{y\rightarrow 0}y^{1-2\gamma}\partial_{y}(\eta_{a,\varepsilon}\hat \delta_{a,\lambda} ) ,
\end{split}
\end{equation*}
whence by (i) of Corollary \ref{cor_rough_estimates} and again by \;$\partial_{y}\eta_{a,\varepsilon}=\varepsilon^{-2}O(y)$\; and \;$\gamma\in (0,1)$\;
\begin{equation*}
\begin{split}
-d_\gamma^*\lim_{y\rightarrow 0}y^{1-2\gamma}\partial_{y}\widetilde{\varphi_{a,\lambda}}^{a} 
= \; &
-d_\gamma^* \,\eta_{a,\varepsilon}\lfloor_{y=0}\,\cdot \lim_{y\rightarrow 0}y^{1-2\gamma}\partial_{y}\hat \delta_{a,\lambda}  
=
c_{n, \gamma}\eta_{a,\epsilon}\delta_{a,\lambda}^{\frac{n+2\gamma}{n-2\gamma}} \lfloor_{y=0},
\end{split}
\end{equation*}
where the last equality is due to \;$\hat \delta_{a,\lambda}$\; being written \;$g_{a}$-normal Fermi-coordinates
and the third property in \eqref{delta_hat}.
Finally we set
\begin{equation}\label{dfbua}
v_{a, \l}=u_a\varphi_{a, \l} \;\; \text{ on } \;\; M  \;\; \text{ and }\;\;\ov u^a=K_{g_a}*u 
\; \text{ on } \; X.
\end{equation}
\subsection{Interaction estimates for \;$v_{a, \l}$\; and related identities}
In this subsection we derive several interaction estimates for \;$v_{a, \l}$\; in \eqref{dfbua}. 
Indeed we give estimates for the higher exponent  interaction as well as for the linear and nonlinear interaction and relate the latter two. 
Recalling \eqref{motivation_g_n_gamma}, we start with the higher exponent interaction estimates.
\begin{lem}
\label{lem_higher_interaction_estimates}
For \;$ i \neq j$\; there holds
\begin{enumerate}[label=(\roman*)]
 \item \quad 
$
\int_{M} v_{a_{i}, \lambda_{i}}^{\alpha}v_{a_{j}, \lambda_{j}}^{\beta}dV_{h}
=
O(\eps_{i,j}^{\beta})
$
\;\;for all \;\;$\alpha +\beta=\frac{2n}{n-2\gamma}, \;\;\alpha>\frac{n}{n-2\gamma}>\beta>0 $
\item \quad
$
\int_{M} v_{a_{i}, \lambda_{i}}^{\frac{n}{n-2\gamma}}v_{a_{j}, \lambda_{j}}^{\frac{n}{n-2\gamma}} dV_{h}
=
O(\eps^{\frac{n}{n-2\gamma}}_{i,j}\ln \eps_{i,j}),
$ 
\end{enumerate}
where 
\begin{equation*}
\varepsilon_{i,j}
=
(
\frac{\lambda_{i}}{\lambda_{j}}
+
\frac{\lambda_{j}}{\lambda_{i}}
+
\lambda_{i}\lambda_{j}\gamma_{n,\gamma}G_{h}^{\frac{2}{2\gamma-n}}(a_{i},a_{j})
)^{\frac{2\gamma-n}{2}},
\;\;\;\;\;\;
\gamma_{n,\gamma}=g_{n,\gamma}^{-\frac{2}{2\gamma-n}}.
\end{equation*}
\end{lem}
\noindent
The type of estimates stated in the latter lemma are standard, 
so we delay the proof to the appendix and pass to establishing the linear and nonlinear interaction estimates mentioned above.
\begin{lem}
\label{interactest}
For \;$i\neq j$\; and letting
\begin{enumerate}[label=(\roman*)]
 \item \quad 
$
e_{i,j}
=
\langle 
v_{a_{i}, \lambda_{i}},v_{a_{j}, \lambda_{j}}
\rangle_{P_{h}^\gamma}
$
 \item \quad
$
\epsilon_{i,j}
=
\int_{M}v_{a_{i}, \lambda_{i}}^{\frac{n+2\gamma}{
n-2\gamma } }v_{a_{j}, \lambda_{j}}dV_{h}
$
 \item \quad 
$
\varepsilon_{i,j}
=
(
\frac{\lambda_{i}}{\lambda_{j}}
+
\frac{\lambda_{j}}{\lambda_{i}}
+
\lambda_{i}\lambda_{j}\gamma_{n,\gamma}G_{h}^{\frac{2}{2\gamma-n}}(a_{i},a_{j})
)^{\frac{2\gamma-n}{2}}
$,\;\;\;\;\;\;$\gamma_{n,\gamma}=g_{n,\gamma}^{-\frac{2}{2\gamma-n}}$
\end{enumerate}
there holds
\begin{enumerate}[label=(\roman*)]
 \item \quad
$
c \leq 
\frac{e_{i,j}}{\varepsilon_{i,j}}, \frac{\epsilon_{i,j}}{\varepsilon_{i,j}}
\leq C
$
 \item \quad 
$
e_{i,j}=c_{n,\gamma}\epsilon_{i,j}(1+o_{\max(\frac{1}{\lambda_{i}}, \frac{1}{
\lambda_{j}} ) }(1))
$
 \item \quad 
$e_{i,j}=c_{n, 4}^\gamma\varepsilon_{i,j}(1+o_{\varepsilon_{i,j}}(1))$\;\; and  \;$\epsilon_{i,j}=c_{n, 3}^\gamma\varepsilon_{i,j}(1+o_{\varepsilon_{i,j}
} (1))$,
\end{enumerate}
where \;$c_{n,3}^{\gamma}, \;c_{n,4}^{\gamma}$\; and \;$c_{n,\gamma}$\; are given by \eqref{eq:yambaf} and \eqref{defc3}.  
\end{lem}
\begin{pf}
Writing abbreviatively \;$\varphi_{k}=\varphi_{a_{k}, \lambda_{k}}$\; we first remark
\begin{equation}\label{epsilon_i_j_symmetry}
\epsilon_{i,j}=\epsilon_{j,i}+o_{\max(\frac{1}{\lambda_{i}}, \frac{1}{\lambda_{j}
} )} (\varepsilon_{i,j} ).
\end{equation} 
We will prove \eqref{epsilon_i_j_symmetry} in the appendix.
Then, since \;$e_{i,j}$\; and \;$\varepsilon_{i,j}$\; are symmetric in \;$i$\; and \;$j$, we may assume 
\;$\frac{1}{\lambda_{i}}\leq \frac{1}{\lambda_{j}}$\; for the rest of the proof. Consider
\begin{equation*}\begin{split}
e_{i,j}
= &
\int_{M}P_{h}^\gamma(u_{a_{i}}  \varphi_{i}) (u_{a_{j}}\varphi_{j}) dV_{h} 
= 
\int_{M}P_{h_{a_{i}}}^\gamma\varphi_{i}\frac{u_{a_{j}}}{u_{a_{i}}}\varphi_{j}dV_{h_{a_
{i}}}
= 
\langle \varphi_{i}  , 
\frac{u_{a_{j}}}{u_{a_{i}}}\varphi_{j}\rangle_{P^\gamma_{h_{a_{i}}}},
\end{split}\end{equation*}
where for the second equality we have applied the covariance property \eqref{eq:confinv} 
to the conformal metric
$$
h_{a_{i}}=u_{a_{i}}^{\frac{4}{n-2\gamma}}h,
$$
cf. \eqref{conformal_factor_properties} and \eqref{h_a_metric}.  
Denoting by \;$y=y_{a_{i}}$\; the geodesic defining function attached to $h_{a_{i}}$  
and choosing close to \;$a_{i}$\; a 
corresponding \;$g_{a_{i}}$-normal Fermi-coordinate system 
$
z=z_{a_{i}}=(y_{a_{i}},x_{a_{i}})=(y,x),\; r=\vert z \vert
$, 
we find from \eqref{Dirichlet_to_Neumann_map}, \eqref{canonical_bubble} and the divergence theorem, that
\begin{equation}\begin{split}\label{quadratic_interaction}
\frac{1}{d_\gamma^*}e_{i,j}
= &
-\lim_{y\rightarrow 0}\int_{M}y^{1-2\gamma}\partial_{y}
\widetilde{\varphi_{i}}^{a_{i } }
\frac{u_{a_{j}}}{u_{a_{i}}}\varphi_{j}dV_{h_{a_{i}}} \\
& +
\int_{X}div_{g_{a_{i}}}(y^{1-2\gamma}\nabla_{g_{a_{i}}}(\overline{
\varphi_{i}}^{a_{i}}-\widetilde{\varphi_{i}}^{a_{i}})
\overline{\frac{u_{a_{j}}}{u_{a_{i}}}\varphi_{j}}^{a_{i}})
dV_{g_{a_{i}}} \\
= & \;
I_{1}+I_{2}.
\end{split}\end{equation}
We first show smallness of $I_{2}$, i.e. \eqref{I2_interaction_integral_estimate}. 
Recalling \eqref{eqdg}, since  \;$D_{g_{a_{i}}}\overline{\frac{u_{a_{j}}}{u_{a_{i}}}\varphi_{j}}^{a_{i}}=0$ due to \eqref{canonical_bubble}, we have
\begin{equation*}\begin{split}
I_{2}
= &
-\int_{X}D_{g_{a_{i}}}(\overline{\varphi_{i}}^{a_{i}}-\widetilde{\varphi_{i}}^{
a_{i}})
\overline{\frac{u_{a_{j}}}{u_{a_{i}}}\varphi_{j}}^{a_{i}}
dV_{g_{a_{i}}}
\end{split}\end{equation*}
and there holds, as we will prove in the appendix,
\begin{equation}\label{varphi_in_a_j_coordinates}
\begin{split}
\vert \overline{\frac{u_{a_{j}}}{u_{a_{i}}}\varphi_{j}}^{a_{i}} \vert
\leq 
C\widetilde{\varphi_{j}}^{a_{j}} 
 \;\;\text{ on }\;\,\; B_{\rho_{0}}^{g_{a_i}, +}(a_{i}).
\end{split}
\end{equation} 
Let us estimate
\begin{equation*}
D_{g_{a_{i}}}(\overline{\varphi_{i}}^{a_{i}}-\widetilde{\varphi_{i}}^{a_{i}})
=-
D_{g_{a_{i}}}\widetilde{\varphi_{i}}^{a_{i}}.
\end{equation*}
Recalling \eqref{properties_of_tilde_phi},  \;$D_{g_{a_{i}}}\widetilde{\varphi_{i}}^{a_{i}}=0$\; on
$B^{g_{a_i}, +}_{2\varepsilon_{i}}(a_{i})^{c}$, whereas on \;$B^{g_{a_i}, +}_{2\varepsilon}(a_{i})$\;  we have
\begin{equation*}
\begin{split}
D_{g_{a_{i}}}\widetilde{\varphi_{i}}^{a}
= & \;
D
(
\eta_{a_{i}, \varepsilon_{i}}
(
\hat  \delta_{a_{i},\lambda_{i}}
-
\frac{\lambda_{i}^{\frac{2\gamma-n}{2}}}{r^{n-2\gamma}}
)
) 
+
(D_{g_{a_{i}}}-D)
(\eta_{a_{i}, \varepsilon_{i}}\hat 
\delta_{a_{i}, \lambda_{i}} 
+
(1-\eta_{a_{i}, \varepsilon_{i}})
\frac{\lambda_{i}^{\frac{2\gamma-n}{2}}}{r^{n-2\gamma}}
) 
\\
&  \;
+
D_{g_{a_{i}}}((1-\eta_{a_{i}, \varepsilon_{i}})
\frac{H_{a_{i}}}{\lambda_{i}^{\frac{n-2\gamma}{2}}}) \\
= & \;
J_{1}+J_{2}+J_{3}
\end{split}
\end{equation*}
by writing \;$\Gamma_{a_{i}}=g_{n,\gamma}(r ^{2\gamma-n}+H_{a_{i}})$\; in \;$g_{a_{i}}$-normal 
Fermi-coordinates, see Subsection 4.2 of \cite{martndia3}. Then
\begin{equation*}\begin{split}
J_{1}
= &
-div (y^{1-2\gamma}\nabla
(
\eta_{a_{i}, \varepsilon_{i}}(\hat \delta_{a_{i}, \lambda_{i}}
-
\frac{\lambda_{i}^{\frac{2\gamma-n}{2}}}{r ^{n-2\gamma}})))
\\
= &
-
\partial_{y}y^{1-2\gamma}\partial_{y}\eta_{a_{i}, \varepsilon_{i}}(\hat 
\delta_{a_{i}, \lambda_{i}}-\frac{\lambda_{i}^{\frac{2\gamma-n}{2}}}{r 
^{n-2\gamma}}) 
+
2y^{1-2\gamma}\nabla \eta_{a_{i}, \varepsilon_{i}}\nabla (\hat 
\delta_{a_{i}, \lambda_{i}}-\frac{\lambda_{i}^{\frac{2\gamma-n}{2}}}{r 
^{n-2\gamma}}) \\
& +
y^{1-2\gamma}\Delta \eta_{a_{i}, \varepsilon_{i}}(\hat 
\delta_{a_{i}, \lambda_{i}}-\frac{\lambda_{i}^{\frac{2\gamma-n}{2}}}{r 
^{n-2\gamma}})
\end{split}\end{equation*}
due to \;$D\hat \delta_{a_{i}, \lambda_{i}}=Dr^{2\gamma-n}=0$.
Then \eqref{cut_off} and Corollary \ref{cor_sharp_estimates} show 
\begin{equation*}
J_{1}
=
o_{\frac{1}{1+\lambda_{i}^{2\gamma}\varepsilon_{i}^{n} 
}}(\frac{y^{1-2\gamma}}{\varepsilon_{i}^{2}  
\lambda_{i}^{\frac{n-2\gamma}{2}}})
=
O(y^{1-2\gamma}\widetilde{\varphi_{i}}^{a_{i}}),
\end{equation*}
whenever \;$\varepsilon_{i}\sim \lambda_{i}^{-\frac{1}{k}}$\; and \;$k\gg 1$\; is 
chosen sufficiently large. Letting 
$$p, q=1,\ldots, n+1
\; \text{ and } \; 
k,l=1,\ldots, n,$$ we evaluate 
\begin{equation}\label{OperatorsDifference_interaction}
\begin{split}
D - D_{g_{a_{i}}}
= &
\frac{\partial_{p}}{\sqrt{g_{a_{i}}}}(\sqrt{g_{a_{i}}}g_{a_{i}}^{p,q}y^{
1-2\gamma } \partial_{
q}\cdot)
-
E_{g_{a_{i}}}
-
\partial_{p}(\delta^{p,q}y^{1-2\gamma}\partial_{q}\cdot) \\
= &
-
E_{g_{a_{i}}}
+
\frac{\partial_{y}\sqrt{g_{a_{i}}}}{\sqrt{g_{a_{i}}}}y^{1-2\gamma}\partial_{y}
+
\frac{\partial_{k}\sqrt{g_{a_{i}}}}{\sqrt{g_{a_{i}}}}g_{a_{i}}^{k,l}y^{1-2\gamma
} \partial_{l
}
+
\partial_{k}((g_{a_{i}}^{k,l}-\delta^{k,l})y^{1-2\gamma}\partial_{l}\cdot)
 \\
= &
O(y^{1-2\gamma})
+
O(y^{2-2\gamma})\partial_{y} 
+
O(y^{1-2\gamma}r)\nabla_{x}
+
O(y^{1-2\gamma}r^{2}+y^{2-2\gamma})\nabla^{2}_{x},
\end{split}
\end{equation}
where we made use of formula (14) of \cite{martndia3}, and
\begin{equation*}
g_{a_{i}}=\delta+O(y+\vert x \vert^{2})\;\;\text{ and } \;\;\det g_{a_{i}}=1+O(y^{2}) 
\;\;\text{ 
on } \;\;
B_{2\varepsilon_{i}}^{g_{a_i}, +}(a_{i}).
\end{equation*}
From \eqref{OperatorsDifference_interaction} and Corollary 
\ref{cor_rough_estimates} we 
then find 
\begin{equation*}
J_{2}
=
O(y^{1-2\gamma}\widetilde{\varphi_{i}}^{a_{i}})
+
O(y^{1-2\gamma}\lambda_{i}y(\widetilde{\varphi_{i}}^{a_{i}})^{\frac{n+2-2\gamma}
{ n-2\gamma } } )
\end{equation*}
for \;$\varepsilon_{i}\sim \lambda_{i}^{-\frac{1}{k}}$\; and \;$k\gg 1$\; 
chosen sufficiently large.
Moreover from \eqref{OperatorsDifference_interaction} we find 
\begin{equation*}
D_{g_{a_{i}}}H_{a_{i}}=-D_{g_{a_{i}}}r^{2\gamma-n}=(D-D_{g_{a_{i}}})r^{2\gamma-n
}
=
O
( 
y^ { 1-2\gamma } r^ { 2\gamma-n-1 } )
, 
\end{equation*}
whence \;$H_{a_{i}}=O(r^{1+2\gamma-n})$, cf. Subsection 4.2 of \cite{martndia3}.
Thus, using \eqref{OperatorsDifference_interaction}, we obtain
\begin{equation*}
\begin{split}
\vert J_{3} \vert
\leq 
\frac{Cy^{1-2\gamma}\chi_{B_{\varepsilon_{i}}(a_{i})^{c}}}
{\lambda_{i}^{\frac{n-2\gamma}{2}}r^{n-2\gamma}}
=
O(y^{1-2\gamma}\widetilde{\varphi_{i}}^{a_{i}})
\;\;\text{ on }\;\;
B_{2\varepsilon_{i}}^{+}(a_{i})
\end{split}
\end{equation*}
for any choice \;$\varepsilon_{i}\sim \lambda_{i}^{-\frac{1}{k}}, \, k\gg 1$.
Collecting terms we conclude
\begin{equation*}
\begin{split}
D_{g_{a_{i}}}\widetilde{\varphi_{i}}^{a_{i}}
= &
O(y^{1-2\gamma}\widetilde{\varphi_{i}}^{a_{i}})
+
O(y^{1-2\gamma}\lambda_{i}y(\widetilde{\varphi_{i}}^{a_{i}})^{\frac{n+2-2\gamma}
{ n-2\gamma } } )\;\;
\text{ on }\;\;
B_{2\varepsilon_{i}}^{g_{a_i}, +}(a_{i}),
\end{split}
\end{equation*}
which in conjunction with \eqref{varphi_in_a_j_coordinates} implies
\begin{equation*}
I_{2}
=
O
(
\int_{B^{g_{a_i}, +}_{2\varepsilon_{i}}(a_{i})}y^{1-2\gamma}\widetilde{\varphi_{i}}^{
a_{i }}\widetilde{\varphi_{j}}^{a_{j}}
)
+
O(\int_{B^{g_{a_i}, +}_{2\varepsilon_{i}}(a_{i})}y^{1-2\gamma}\lambda_{i}y(\widetilde{
\varphi_{i}}^{a_{i}})^{\frac{n+2-2\gamma}
{ n-2\gamma } }\widetilde{\varphi_{j}}^{a_{j}})
. 
\end{equation*}
As we will prove in the appendix, this gives
\begin{equation}\label{I2_interaction_integral_estimate}
I_{2}
=
o_{\frac{1}{\lambda_{i}}}(\varepsilon_{i,j}). 
\end{equation} 
Recalling \eqref{quadratic_interaction} we thus arrive at
\begin{equation*}\begin{split}
\frac{1}{d_\gamma^*}e_{i,j}
= &
-\lim_{y\rightarrow 0}\int_{M}y^{1-2\gamma}\partial_{y}\widetilde{\varphi_{i}}^{a_{i}
}
\frac{u_{a_{j}}}{u_{a_{i}}}\varphi_{j}dV_{h_{a_{i}}} 
+
o_{\frac{1}{\lambda_{i}}}(\varepsilon_{i,j}),
\end{split}\end{equation*}
whence by virtue of \eqref{properties_of_tilde_phi} 
\begin{equation*}\begin{split}
\frac{1}{d_\gamma^*}e_{i,j}
= 
\int_{M}
\frac{c_{n, \gamma}}{d_\gamma^*}\eta_{a_{i}, \varepsilon_{i}}\delta_{i}^{\frac{n+2\gamma}{n-2\gamma}}
\frac{u_{a_{j}}}{u_{a_{i}}}\varphi_{j}dV_{h_{a_{i}}}
+
o_{\frac{1}{\lambda_{i}}}(\varepsilon_{i,j}),
\end{split}\end{equation*} 
where \;$\d_{i}=\d_{a_{i}, \l_{i}}$. Moreover, using Corollary \ref{cor_sharp_estimates}, we have
\begin{equation}\label{interaction_expansion_intermediate_form}
\begin{split}
\int_{M} 
\eta_{a_{i}, \varepsilon_{i}}\delta_{i}^{\frac{n+2\gamma}{n-2\gamma}}\frac{u_{a_
{j}}}{u_{a_{i}}}\varphi_{j}dV_{h_{a_{i}}} 
= &
\int_{M} 
\eta_{a_{i}, \varepsilon_{i}}(u_{a_{i}}\delta_{i})^{\frac{n+2\gamma}{n-2\gamma}}
(u_{a_{j}}\varphi_{j})dV_{h}  \\
= &
\underset{B^{h_{a_i}}_{\varepsilon_{i}}(a_{i})}\int 
(u_{a_{i}}\delta_{i})^{\frac{n+2\gamma}{n-2\gamma}}
(u_{a_{j}}\varphi_{j})dV_{h} 
+
O(\frac{\varepsilon_{i}^{-n-2\gamma}}{\lambda_{i}^{\frac{n+2\gamma}{2}}\lambda_{
j}^{\frac{n-2\gamma}{2}}})\\
= &
\int_{M} (u_{a_{i}}\varphi_{i})^{\frac{n+2\gamma}{n-2\gamma}}
(u_{a_{j}}\varphi_{j})dV_{h} 
+
o_{\frac{1}{\lambda_{i}}}(\varepsilon_{i,j})
=
\epsilon_{i,j}+o_{\frac{1}{\lambda_{i}}}(\varepsilon_{i,j})
\end{split}
\end{equation}
for any choice \;$\varepsilon_{i}\sim \lambda_{i}^{-\frac{1}{k}}, \, k$\; large. 
We conclude
\begin{equation}\begin{split}\label{e_ij_epsilon_ij_relation}
e_{i,j}
= &
c_{n,\gamma}\epsilon_{i,j}
+
o_{\frac{1}{\lambda_{i}}}(\varepsilon_{i,j}).
\end{split}\end{equation}
We turn to analyse \;$\epsilon_{i,j}$.
From \eqref{interaction_expansion_intermediate_form} we have
\begin{equation*}\begin{split}
\epsilon_{i,j}
= &
\underset{B^{h_{a_i}}_{\varepsilon_{i}}(a_{i})}\int 
(u_{a_{i}}\delta_{i})^{\frac{n+2\gamma}{n-2\gamma}}
(u_{a_{j}}\varphi_{j})dV_{h}
+
o_{\frac{1}{\lambda_{i}}}(\varepsilon_{i,j})
= 
\underset{B^{h_{a_i}}_{\varepsilon_{i}}(a_{i})}\int 
\delta_{i}^{\frac{n+2\gamma}{n-2\gamma}}
\frac{u_{a_{j}}}{u_{a_{i}}}\varphi_{j}dV_{h_{a_{i}}}
+
o_{\frac{1}{\lambda_{i}}}(\varepsilon_{i,j})
.
\end{split}\end{equation*}
Changing coordinates and rescaling we obtain
\begin{equation*}\begin{split}
\epsilon_{i,j}
= &
\underset{B^{h_{a_i}}_{ \varepsilon_{i}\lambda_{i} }(0)}{\int}
\frac{\frac{u_{ a _{j}}}{u_{a_{i}}}( 
\exp^{h_{a_{i}}}_{a_i}\frac{x}{\lambda_{i}})}{(1+r^{2})^{\frac{n+2\gamma}{2}}} 
\left(\frac{1}{\frac{\lambda_{i} }{ \lambda_{j} }+\lambda_{i}\lambda_{j}\gamma_{n}G_{ a 
_{j}}^{\frac{2}{2\gamma-n}}(\exp^{h_{ a _{i}}}_{a_i} 
\frac{x}{\lambda_{i}})}\right)^{\frac{n-2\gamma}{2}}
dV_{h_{a_{i}}}
+
o_{\frac{1}{\lambda_{i}}}(\varepsilon_{i,j})
.
\end{split}\end{equation*}
In particular
\begin{equation}\label{smallness_if_and_only_if}
\epsilon_{i,j}\ll 1 \Longleftrightarrow \varepsilon_{i,j}\ll 1
\end{equation}
and, as we will prove in the appendix, 
\begin{equation}\label{epsilon_varepsilon_relation}
\begin{split}
\epsilon_{i,j}
= &
c_{n, 3}^{\gamma}\varepsilon_{i,j}(1+o_{\varepsilon_{i,j}}(1)),
\end{split}\end{equation}
Now
\eqref{smallness_if_and_only_if} and \eqref{epsilon_varepsilon_relation} show
$c\leq \rfrac{\epsilon_{i,j}}{\varepsilon_{i,j}}\leq C$, 
whence by virtue 
of \eqref{e_ij_epsilon_ij_relation} for \;$\frac{1}{\lambda_{i}}\ll 1$\; sufficiently small
\begin{equation*}
c\leq \rfrac{e_{i,j}}{\varepsilon_{i,j}}\leq C.
\end{equation*}
This shows (i), i.e. that \;$e_{i,j}, \epsilon_{i,j}$\; and  \;$\varepsilon_{i,j}$\; are pairwise comparable. In particular (ii) follows from 
\eqref{e_ij_epsilon_ij_relation}. In view of 
\eqref{epsilon_varepsilon_relation} we are left with verifying 
\;$e_{i,j}=c_{n, 4}^\gamma\varepsilon_{i,j}(1+o_{\varepsilon_{i,j}}(1))$, but 
this follows easily from (ii) and \eqref{epsilon_varepsilon_relation} combined with \eqref{defc3}.
The proof of the lemma is thereby complete. 
\end{pf}

\section{Locally flat conformal infinities of Poincar\'e-Einstein manifolds}\label{pceinsten}
In this section we discuss Fermi-coordinates in case of a Poincar\'e-Einstein manifold \;$(X, g^+)$\; with  locally flat conformal infinity \;$(M, [h])$. 
Furthermore, in this particular case, we establish sharp \;$L^{\infty}$-estimates for 
\;$D_{g}(\ov{\varphi_{a, \l}}^a-\wtilde{\varphi_{a, \l}}^a)$ and \;$\ov{\varphi_{a, \l}}^a-\wtilde{\varphi_{a, \l}}^a$
as well as selfaction estimates for \;$v_{a, \l}$. 
\subsection{Fermi-coordinates in the particular case}
By our assumptions we have  
\begin{enumerate}[label=(\roman*)]
 \item a geodesic defining function \;$y$\; splitting the metric
\begin{equation*}
g=y^{2}g^{+}, \; g=dy^{2}+h_{y}\;\;\text{ near }\;\;M\;\;\text{ and }\;\; h=h_{y}\lfloor_{M}
\end{equation*} 
and for every \;$a\in M$\;  a conformal factor \;$u_{a}>0$\; as in \eqref{conformal_factor_properties},
whose conformal metric \;$h_{a}=u_{a}^{\frac{4}{n-2\gamma}}h$\; close to \;$a$\; admits an Euclidean 
coordinate system, \;$h_{a}=\delta$\; on $B_{\epsilon}^{h_a}(a)$. As clarified in Subsection 3.2 in \cite{martndia3} and recalling Remark \ref{eq:minimal}, this gives rise to a geodesic defining 
function \;$y_{a}$, for which
\begin{equation*}
g_{a}=y_{a}^{2}g^{+}, \;\; g_{a}=dy_{a}^{2}+h_{a,y_{a}}\;\;\text{near}\;\;M\;\;\text{ with }\;\; 
h_{a}=h_{a, y_a}\lfloor_{M}\;\;\text{ and }\;\;\delta=h_{a}\lfloor_{ 
B_{\epsilon}^{h_a}(a)},
\end{equation*}
the boundary \;$(M, [h_a])$\; is totally geodesic and the extension operator \;$D_{g_{a}}$\; is positive. 

\item
in \;$g_{a}$-normal  Fermi-coordinates around $a$ for some small \;$\epsilon >0$
\begin{equation}\label{flatness}
g_{a}=\delta+O(|y_{a}|^{n}) \;\;\text{ on } \;\;B_{\epsilon}^{g_a, +}(a),
\end{equation} 
as observed by Kim-Musso-Wei in case \;$n\geq 3$, cf. Lemma 43 in \cite{kmw1}, and for \;$n=2$\; due to Remark \ref{eq:minimal} and the existence of isothermal coordinates.
\end{enumerate}

\subsection{Comparing  Schoen's bubbles and the Projective bubbles}
In this subsection we compare the  Schoen's bubbles \;$\wtilde{\varphi_{a, \l}}^a$ and our Projective bubbles \;$\ov{\varphi_{a, \l}}^a$. 
Indeed we establish sharp $L^{\infty}$-estimates for \;$D_{g_a}\left(\ov{\varphi_{a, \l}}^a-\wtilde{\varphi_{a, \l}}^a\right)$ and,  
using the maximum principle for \;$D_{g_a}$\; under Dirichlet boundary conditions, sharp $\; L^{\infty}$-estimates for $\ov{\varphi_{a, \l}}^a-\wtilde{\varphi_{a, \l}}^a$.  

\begin{lem}
\label{lem_bubble comparsion} 
Writing in \;$g_{a}$-normal Fermi-coordinates 
\begin{equation*}
\frac{\Gamma_{a}}{g_{n,\gamma}}=r ^{2\gamma-n}+H_{a}=r^{2\gamma-n}+M_{a}+O(\max(r ,r ^{
2\gamma})),
\end{equation*} 
where \;$M_{a}$\; depends on \;$a\in M$ only, cf. Theorem 1.4 in \cite{martndia3}, there holds
\begin{enumerate}[label=(\roman*)]
 \item \quad
$
D_{g_{a}}\left(\ov{\varphi_{a,\lambda}}^{a}-\widetilde{\varphi_{a,\lambda}}^{a}\right)
= 
-D_{g_{a}}(\widetilde{\varphi_{a,\lambda}}^{a})=
D_{g_{a}}(\eta_{a,\varepsilon}\frac{H_{a}}{\lambda^{\frac{n-2\gamma}{2}}})
+
O(\frac{\lambda^{\frac{2\gamma-n}{2}}}{y^{\gamma}r ^{1-\gamma}}) \\
_{}\quad\quad\quad\;\;\,\quad\quad\quad\quad\quad\;\,\, = 
\frac{M_{a}}{\lambda^{\frac{n-2\gamma}{2}}}div_{g_{a}}(y^{1-2\gamma}\nabla_{g_{
a}}\eta_{a,\varepsilon})
+
O(\frac{\max(1,\varepsilon^{1-2\gamma})}{\lambda^{\frac{n-2\gamma}{2}}y^{\gamma}
r ^{1-\gamma}}) 
$
 \item \quad
$
\overline{\varphi_{a,\lambda}}^{a}-\widetilde{\varphi_{a,\lambda}}^{a}
=
O(\frac{1}{\lambda^{\frac{n-2\gamma}{2}}}),
$
\end{enumerate}
provided \;$\varepsilon\sim \lambda^{-\frac{1}{k}}$\; for some \;$k=k(n,\gamma)$\; sufficiently large.
\end{lem}
\begin{pf}
Recalling \eqref{eqdg} and \eqref{constructed_bubble}, since \;$D_{g_{a}}\Gamma_{a}=0$, we have 
\begin{equation*}
\begin{split}
D_{g_{a}}\widetilde{\varphi_{a,\lambda}}^{a}
= & \;
D_{g_{a}}(\eta_{a,\varepsilon}(\hat 
\delta_{a,\lambda}-\lambda^{\frac{2\gamma-n}{2}}\frac{\Gamma_{a}}{g_{n,\gamma}}))
\\
= & \;
D(\eta_{a,\varepsilon}(\hat 
 \delta_{a,\lambda}-\frac{\lambda^{\frac{2\gamma-n}{2}}}{r ^{n-2\gamma}}))
+
(D_{g_{a}}-D)(\eta_{a,\varepsilon}(\hat 
\delta_{a,\lambda}-\frac{\lambda^{\frac{2\gamma-n}{2}}}{r ^{n-2\gamma}}))  - 
D_{g_{a}}(\eta_{a,\varepsilon} \frac{H_{a}}{\lambda^{\frac{n-2\gamma}{2}}})
\\ 
= & \;
I_{1}+I_{2}+I_{3}
\end{split}
\end{equation*}
and start evaluating 
\begin{equation*}\begin{split}
I_{1}
= &
-div (y^{1-2\gamma}\nabla
(
\eta_{a,\varepsilon}(\hat \delta_{a,\lambda}
-
\frac{\lambda^{\frac{2\gamma-n}{2}}}{r ^{n-2\gamma}})))
\\
= &
-
\partial_{y}y^{1-2\gamma}\partial_{y}\eta_{a,\varepsilon}(\hat 
\delta_{a,\lambda}-\frac{\lambda^{\frac{2\gamma-n}{2}}}{r ^{n-2\gamma}}) 
+
2y^{1-2\gamma}\nabla \eta_{a,\varepsilon}\nabla (\hat 
\delta_{a,\lambda}-\frac{\lambda^{\frac{2\gamma-n}{2}}}{r ^{n-2\gamma}})  +
y^{1-2\gamma}\Delta \eta_{a,\varepsilon}(\hat 
\delta_{a,\lambda}-\frac{\lambda^{\frac{2\gamma-n}{2}}}{r ^{n-2\gamma}}),
\end{split}\end{equation*}
where we made use of \;$D\hat \delta_{a,\lambda}=Dr^{2\gamma-n}=0$. Then
\eqref{cut_off} and Corollary \ref{cor_sharp_estimates}  show 
\begin{equation*}
I_{1}
=
o_{\frac{1}{1+\lambda^{2\gamma}\varepsilon^{n} 
}}(\frac{y^{1-2\gamma}}{\varepsilon^{2}  
\lambda^{\frac{n-2\gamma}{2}}})
=
O(\frac{\lambda^{\frac{2\gamma-n}{2}}}{y^{\gamma}r ^{1-\gamma}}),
\end{equation*}
whenever \;$\varepsilon\sim \lambda^{-\frac{1}{k}}$\; and \;$k\gg 1$\; is 
chosen sufficiently large. Secondly, letting 
$$p,q=1,\ldots,n+1 \; \text{ and } \; i,j=1,\ldots, n,$$
we use of formula (14) of \cite{martndia3}, \eqref{flatness} and the splitting of the 
metric to evaluate  
\begin{equation}\label{OperatorsDifference}
\begin{split}
D - D_{g_{a}}
= &
\frac{\partial_{p}}{\sqrt{g_{a}}}(\sqrt{g_{a}}g_{a}^{p,q}y^{1-2\gamma}\partial_{
q}\cdot)
-
E_{g_{a}}
-
\partial_{p}(\delta^{p,q}y^{1-2\gamma}\partial_{q}\cdot) \\
= &
\frac{\partial_{y}\sqrt{g_{a}}}{\sqrt{g_{a}}}y^{1-2\gamma}\partial_{y} 
+
\frac{\partial_{i}\sqrt{g_{a}}}{\sqrt{g_{a}}}g_{a}^{i,j}y^{1-2\gamma}\partial_{j
} +
\partial_{i}((g_{a}^{i,j}-\delta^{i,j})y^{1-2\gamma}\partial_{j}\cdot)
-
E_{g_{a}} \\
= &
O(y^{n-1-2\gamma})
+
O(y^{n-2\gamma})\partial_{y}
+
O(y^{n+1-2\gamma})\nabla_{x}
+
O(y^{n+1-2\gamma})\nabla^{2}_{x},
\end{split}
\end{equation}
cf \eqref{OperatorsDifference_interaction} for the non flat case.  From \eqref{OperatorsDifference} and Corollary \ref{cor_rough_estimates} we 
then find 
\begin{equation*}
I_{2}=O(\frac{\lambda^{\frac{2\gamma-n}{2}}}{y^{\gamma}r ^{1-\gamma}}).
\end{equation*}
Moreover \eqref{OperatorsDifference} and Theorem 1.4 in \cite{martndia3} imply 
\;$
D_{g_{a}}H_{a}=(D-D_{g_{a}})r ^{2\gamma-n}=O(\frac{1}{y^{\gamma}r ^{ 
1-\gamma}}),
$ 
,
whence
\begin{equation*}
\begin{split}
I_{3}
=
-
\frac{H_{a}}{\lambda^{\frac{n-2\gamma}{2}}}div_{g_{a}}(y^{1-2\gamma}\nabla_{g_{a
}}\eta_{a,\varepsilon})
-
2y^{1-2\gamma}
\langle \nabla_{g_{a}} \eta_{a,\varepsilon}, \frac{\nabla_{g_{a}} 
H_{a}}{\lambda^{\frac{n-2\gamma}{2}}}\rangle_{g_{a}}
+
O(\frac{\lambda^{\frac{2\gamma-n}{2}}}{y^{\gamma}r ^{1-\gamma}}).
\end{split}
\end{equation*}
Due to  the structure of \;$H_{a}$, cf. Theorem 1.4 in \cite{martndia3}, we have 
\begin{equation*}
H_{a}=M_{a}+O(\max( r ,r ^{2\gamma})) 
\;\;\text{ and }\;\; 
\nabla_{g_{a}}H_{a}=
O(\max(1,r ^{2\gamma-1})
\end{equation*}
with \;$M_{a}$\; denoting the constant in the expansion of the Green's 
function \;$\Gamma_a$.
Recalling \eqref{cut_off} we then get
\begin{equation*}
\begin{split}
I_{3}
=
-\frac{M_{a}}{\lambda^{\frac{n-2\gamma}{2}}}div_{g_{a}}(y^{1-2\gamma}\nabla_{g_{
a}}\eta_{a,\varepsilon})
+
O(\frac{\max(1,\varepsilon^{1-2\gamma})}{\lambda^{\frac{n-2\gamma}{2}}y^{\gamma}
r ^{1-\gamma}}).
\end{split}
\end{equation*}
Collecting terms we conclude
\begin{equation}\label{Dga_Tilde_Bubble_Estimate}
\begin{split}
D_{g_{a}}\widetilde{\varphi_{a,\lambda}}^{a}
= &
-D_{g_{a}}(\eta_{a,\varepsilon}\frac{H_{a}}{\lambda^{\frac{n-2\gamma}{2}}})
+
O(\frac{\lambda^{\frac{2\gamma-n}{2}}}{y^{\gamma}r ^{1-\gamma}}) \\
= &
-\frac{M_{a}}{\lambda^{\frac{n-2\gamma}{2}}}div_{g_{a}}(y^{1-2\gamma}\nabla_{g_{
a}}\eta_{a,\varepsilon})
+
O(\frac{\max(1,\varepsilon^{1-2\gamma})}{\lambda^{\frac{n-2\gamma}{2}}y^{\gamma}
r ^{1-\gamma}}) 
.
\end{split}
\end{equation}
This and, that  $D_{g_a}\overline{\varphi_{a,\lambda}}^{a}=0$ by definition, show (i).  
And (ii) follows from the maximum principle for $D_{g_a}$\; under Dirichlet boundary conditions. Indeed consider
\begin{equation*}
\chi_{a,\rho}=\chi(\rho^{-2}y^{2}), \;\;\chi=1\;\;\text{ on }\;\; [0,1] \;\;
\text{ and } \;\; \chi=0 \;\;\text{ on } \;\;(2,+\infty),
\end{equation*}
i.e. a cut-off function for the boundary, and 
\begin{equation}\label{appearance_of_the_log}
\psi
=
\psi_{a,\lambda,\varepsilon,\rho}
=
\overline{\varphi_{a,\lambda}}^{a}
-
\widetilde{\varphi_{a,\lambda}}^{a}
-
\eta_{a,\varepsilon}\frac{H_{a}}{\lambda^{\frac{n-2\gamma}{2}}}
\pm
C\eta_{a,\rho}\frac{y^{2\gamma}\ln y}{\lambda^{\frac{n-2\gamma}{2}}}
\end{equation}
with constants \;$\rho, \;C>0$.
We then find from \eqref{canonical_bubble}, \eqref{constructed_bubble},
\eqref{Dga_Tilde_Bubble_Estimate} and
$
D(y^{2\gamma}\ln y)=-\frac{2\gamma}{y},
$
that 
\begin{equation*}
\begin{cases}
D_{g_{a}}\psi
=
O(\lambda^{\frac{2\gamma-n}{2}}) \;\;\text{ in } \;\;X
\\
D_{g_{a}}
\psi
\;\substack{\leq\\\geq}\;
0\;\;
\text{ in } \;\;[d_{g_{a}}(\cdot,M)\leq \rho]
\\
\psi
=
O(\lambda^{\frac{2\gamma-n}{2}})\;\;
\text{ on } \;\;M
\end{cases}
\end{equation*}
for suitable \;$0<\rho, \;C^{-1}\ll 1$ and according to Proposition 3.1 in \cite{martndia3} we 
may solve
\begin{equation*}
\begin{cases}
D_{g_{a}}u=c_{1}\;\; \text{ in } \;\;X \\
u=c_{2} \;\,\;\text{ on }\,\;\; M
\end{cases}
\end{equation*}
for constants \;$c_{1}, c_{2}\in \R$\; with \;$u\in 
C^{\infty}(\overline X, g_a)+y^{2\gamma}C^{\infty}(\overline X, g_a)$.
\end{pf}
\subsection{Selfaction estimates for \;$v_{a, \l}$\; and the emergence of mass}
In this short subsection we derive two selfaction estimates for \;$v_{a, \l}$. In particular we identify the fractional analogue of the mass for the classical Yamabe problem in the locally conformally flat case. 

\begin{lem}
\label{selfactestv}
Writing in \;$g_{a}$-normal Fermi-coordinates \;$z=z_{a}=(y_{a},x_{a})=(y,x)$ 
\begin{equation*}
\frac{\Gamma_{a}}{g_{n,\gamma}}=r ^{2\gamma-n}+H_{a}=r^{2\gamma-n}+M_{a}+O(\max(r ,r ^{
2\gamma})),
\end{equation*}
where \;$M_{a}$\;  depends on \;$a\in M$\; only, cf. Theorem 1.4 in \cite{martndia3}, there holds
\begin{enumerate}[label=(\roman*)]
 \item \quad 
$
\langle v_{a,\lambda}  , v_{a,\lambda}\rangle_{P_{h}^\gamma}
= 
\langle \varphi_{a,\lambda}  , \varphi_{a,\lambda}\rangle_{P_{h_{a}}^\gamma} 
=  
c_{n, 2}^\gamma-c_{n, \gamma}^*\;\dfrac{M_{a}}{\lambda^{n-2\gamma}}+o_{\frac{1}{\lambda}}(\frac{
1}{\lambda^{n-2\gamma}}),
$
 \item \quad 
$
\int_{M}v_{a,\lambda}^{\frac{2n}{n-2\gamma}}dV_{h}
=
\int_{M}\varphi_{a,\lambda}^{\frac{2n}{n-2\gamma}}dV_{h_{a}}
=
c_{n, 1}^\gamma+o_{\frac{1}{\lambda}}(\lambda^{2\gamma-n})
$
\end{enumerate}
provided \;$\varepsilon\sim \lambda^{-\frac{1}{k}}$\; for \;$k\in \N$\; sufficiently 
large. 
\end{lem}
\begin{pf}
Recalling \eqref{canonical_bubble} and \eqref{properties_of_tilde_phi}, we calculate the quadratic form
\begin{equation*}\begin{split}
\frac{1}{d_\gamma^*}\langle \varphi_{a,\lambda}  ,\varphi_{a,\lambda}\rangle_{P_{h_{a}}^\gamma}
= & \;
\frac{1}{d_\gamma^*}
\int_{M}P_{g_{a}}\varphi_{a,\lambda}  \varphi_{a,\lambda}dV_{h_{a}} \\
= & \;
-\lim_{y\rightarrow 0}\int_{M}y^{1-2\gamma}\partial_{y}\widetilde{\varphi_{a,\lambda}}^{a}\widetilde
{\varphi_{a,\lambda}}^{a} dV_{h_{a}} 
-
\lim_{y\rightarrow 0}\int_{M}y^{1-2\gamma}\partial_{y}(\overline{\varphi_{a,\lambda}}^{a}
-\widetilde{\varphi_{a,\lambda}}^{a})\widetilde{\varphi_{a,\lambda}}^{a} 
dV_{h_{a}}  \\
= & \;
\frac{c_{n, \gamma}}{d_\gamma^*}
\int_{M} 
\eta_{a,\varepsilon}\delta_{a,\lambda}^{\frac{n+2\gamma}{n-2\gamma}}\varphi_{a,
\lambda}dV_{h_{a}} 
+
\int_{X}div_{g_{a}}(y^{1-2\gamma}\nabla_{g_{a}}(\overline{\varphi_{a,\lambda}}^{
a}-\widetilde{\varphi_{a,\lambda}}^{a})\widetilde{\varphi_{a,\lambda}}^{a})dV_
{g_{a}}
\end{split}\end{equation*}
and from \eqref{delta_boundary}, \eqref{bubble_boundary} and 
\;$h_{a}=g_{a}\lfloor_{M}=\delta$\; close to \;$a$\; we easily find
\begin{equation}\label{super_quadratic_bubble_energy_expansion}
\begin{split}
\int_{M}\eta_{a,\varepsilon}\delta_{a,\lambda}^{\frac{n+2\gamma}{n-2\gamma}}
\varphi_{a,\lambda}dV_{h_{a}}
= &
\int_{\R^{n}}\delta_{0,\lambda}^{\frac{2n}{n-2\gamma}}
+
O(\int_{\R^{n}\setminus 
B_{\varepsilon}(0)}\delta_{0,\lambda}^{\frac{2n}{n-2\gamma}})
=
c_{n,1}^{\gamma}+O(\frac{1}{(\varepsilon \lambda)^{n}}),
\end{split}
\end{equation}
whence by integration by parts and any choice \;$\varepsilon \sim 
\lambda^{-\frac{1}{k}}, \,k$\; sufficiently large
\begin{equation*}\begin{split}
\frac{1}{d_\gamma^*}\langle \varphi_{a,\lambda}  , \varphi_{a,\lambda}\rangle_{P_{h_{a}}^\gamma}
= & \;
\frac{c_{n, 2}^\gamma}{d_\gamma^*} 
+
\int_{X} y^{1-2\gamma}
\langle 
\nabla_{g_{a}}(\overline{\varphi_{a,\lambda}}^{a}-\widetilde{\varphi_{a,\lambda}
}^{a}),
\nabla_{g_{a}}\widetilde{\varphi_{a,\lambda}}^{a}\rangle_{g_{a}} 
+
E_{g_{a}}(\overline{\varphi_{a,\lambda}}^{a}-\widetilde{\varphi_{a,\lambda}}^{a}
)\widetilde{\varphi_{a,\lambda}}^{a} dV_{g_{a}}\\
&  -
\int_{X}D_{g_{a}}(\overline{\varphi_{a,\lambda}}^{a}-\widetilde{\varphi_{a,
\lambda}}^{a})\widetilde{\varphi_{a,\lambda}}^{a}dV_{g_{a}}
+
o_{\frac{1}{\lambda}}(\lambda^{2\gamma-n})
\\
= & \;
\frac{c_{n, 2}^\gamma}{d_\gamma^*}
+
\int_{X}div_{g_{a}}(y^{1-2\gamma}\nabla_{g_{a}}\widetilde{\varphi_{a,\lambda}}^{
a}(\overline{\varphi_{a,\lambda}}^{a}-\widetilde{\varphi_{a,\lambda}}^{a})) 
dV_{g_{a}} \\
& -
\int_{X}D_{g_{a}}(\overline{\varphi_{a,\lambda}}^{a}-\widetilde{\varphi_{a,
\lambda}}^{a})\widetilde{\varphi_{a,\lambda}}^{a} dV_{g_{a}}   
+
\int_{X} 
D_{g_{a}}\widetilde{\varphi_{a,\lambda}}^{a}(\overline{\varphi_{a,\lambda}}^{a}
-\widetilde{\varphi_{a,\lambda}}^{a})dV_{g_{a}}
+
o_{\frac{1}{\lambda}}(\lambda^{2\gamma-n})\\
= & \;
I_{1}+\ldots+I_{5}.
\end{split}\end{equation*}
Clearly \;$I_{2}=0$\;  and from Lemma \ref{lem_bubble comparsion} (i) we obtain
\begin{equation*}\begin{split}
I_{3}
= &
-\int_{X}D_{g_{a}}(\eta_{a,\varepsilon}\frac{H_{a}}{\lambda^{\frac{n-2\gamma}{2}
}})\widetilde{\varphi_{a,\lambda}}^{a}dV_{g_{a}}
+
o_{\frac{1}{\lambda}}(\frac{1}{\lambda^{n-2\gamma}})
,
\end{split}\end{equation*}
since, whenever we choose \;$\varepsilon\sim \lambda^{-\frac{1}{k}}$, then
\begin{equation*}
\begin{split}
\int_{B_{2\varepsilon}^{+}}
\frac{\widetilde{\varphi_{a,\lambda}}^{a}}{y^{\gamma} r^{1-\gamma}}
\leq & \;
\frac{C}{\lambda^{\frac{n+2\gamma}{2}}}\int_{B_{2\varepsilon \lambda}^{+}}
\frac{1}{y^{\gamma}(y^{2}+\vert x 
\vert)^{\frac{1-\gamma}{2}}}(\frac{1}{1+y^{2}+\vert x 
\vert^{2}})^{\frac{n-2\gamma}{2}} \\
\leq & \;
\frac{C}{\lambda^{\frac{n+2\gamma}{2}}}\int^{\varepsilon\lambda}_{0}\frac{dy}{y^
{ \gamma } }
\int_{Q^{n}_{\varepsilon\lambda}}\frac{1}{\vert x 
\vert^{1-\gamma}}(\frac{1}{1+\vert x \vert^{2}})^{\frac{n-2\gamma}{2}}
=
\frac{\varepsilon^{2\gamma}}{\lambda^{\frac{n-2\gamma}{2}}},
\end{split}
\end{equation*}
where \;$Q^{n}_{l}=[-l/2,l/2]^{n}$.
Moreover from Lemma \ref{lem_bubble comparsion} we easily find 
\begin{equation*}
D_{g_{a}}\widetilde{\varphi_{a,\lambda}}^{a}=O(\frac{\varepsilon^{-2\gamma}}{
\lambda^{\frac{n-2\gamma}{2}}y^{\gamma}r^{1-\gamma}} )
\;\;\text{ and }\;\;
\overline{\varphi_{a,\lambda}}^{a}-\widetilde{\varphi_{a,\lambda}}^{a}=O(\frac{1
} {\lambda^{\frac{n-2\gamma}{2}}}),
\end{equation*}
whence \;$I_{4}=O(\frac{\varepsilon^{n-2\gamma}}{\lambda^{n-2\gamma}})$. 
Collecting terms we obtain
\begin{equation*}\begin{split}
\frac{1}{d_\gamma^*}\langle \varphi_{a,\lambda}  , &\varphi_{a,\lambda}\rangle_{P_{h_{a}}^\gamma}
= 
\frac{c_{n, 2}^\gamma}{d_\gamma^*}
-
\int_{X}D_{g_{a}}(\eta_{a,\varepsilon}\frac{H_{a}}{\lambda^{\frac{n-2\gamma}{2}}
})\widetilde{\varphi_{a,\lambda}}^{a}dV_{g_{a}}
+
o_{\frac{1}{\lambda}}(\frac{1}{\lambda^{n-2\gamma}}). 
\end{split}
\end{equation*}
From Lemma \ref{lem_bubble comparsion} (i) we then find
\begin{equation*}
\begin{split}
\int_{X}D_{g_{a}}(\eta_{a,\varepsilon}\frac{H_{a}}{\lambda^{\frac{n-2\gamma}{2}}
})\widetilde{\varphi_{a,\lambda}}^{a}dV_{g_{a}}
= &
\frac{M_{a}}{\lambda^{\frac{n-2\gamma}{2}}}
\int_{X}div_{g_{a}}(y^{1-2\gamma}\nabla_{g_{a}}\eta_{a,\varepsilon})\widetilde{
\varphi_{a,\lambda}}^{a}dV_{g_{a}}  
+
O(\frac{\max(\varepsilon,\varepsilon^{2\gamma})}{\lambda^{n-2\gamma}})
,
\end{split}
\end{equation*}
where integrating by parts, using \eqref{cut_off}, we have   
\begin{equation*}
\begin{split}
\int_{X}div_{g_{a}}(y^{1-2\gamma}\nabla_{g_{a}}\eta_{a,\varepsilon})\widetilde{
\varphi_{a,\lambda}}^{a}dV_{g_{a}} 
= &
-\int_{X}div_{g_{a}}(y^{1-2\gamma}\eta_{a,\varepsilon}\nabla_{g_{a}}\widetilde{
\varphi_{a,\lambda}}^{a})dV_{g_{a}} \\
& +
\int_{X}\eta_{a,\varepsilon}div_{g_{a}}(y^{1-2\gamma}\nabla_{g_{a}}\widetilde{
\varphi_{a,\lambda}}^{a})dV_{g_{a}}
\end{split}
\end{equation*}
with the latter integral being of order  \;$O(\frac{\varepsilon^{n-2\gamma}}{\lambda^{n-2\gamma}})$, since
\begin{equation*}
\begin{split}
div_{g_{a}}(y^{1-2\gamma}\nabla_{g_{a}}\widetilde{\varphi_{a,\lambda}}^{a})
= &
-D_{g_{a}}\widetilde{\varphi_{a,\lambda}}^{a}-E_{g_{a}}\widetilde{\varphi_{a,
\lambda}}^{a} \\
= &
-
\frac{M_{a}}{\lambda^{\frac{n-2\gamma}{2}}}div_{g_{a}}(y^{1-2\gamma}\nabla_{g_{a
}}\eta_{a,\varepsilon})
+
O(\frac{\max(1,\varepsilon^{1-2\gamma})}{\lambda^{\frac{n-2\gamma}{2}}y^{\gamma}
r_{a}^{1-\gamma}}) 
= 
O(\frac{\varepsilon^{-2\gamma}\lambda^{\frac{2\gamma-n}{2}}}{y^{\gamma}r_{a}^{
1-\gamma}})
\end{split}
\end{equation*}
due to formula 14 in \cite{martndia3}, \eqref{flatness}, \eqref{Dga_Tilde_Bubble_Estimate}.
Collecting terms we derive
\begin{equation*}\begin{split}
\frac{1}{d_\gamma^*}\langle \varphi_{a,\lambda}  , &\varphi_{a,\lambda}\rangle_{P_{h_{a}}^\gamma}
= 
\frac{c_{n, 2}^\gamma}{d_\gamma^*}
-
\frac{M_{a}}{\lambda^{\frac{n-2\gamma}{2}}}\int_{X}div_{g_{a}}(y^{1-2\gamma}
\eta_{a,\varepsilon}\nabla_{g_{a}}\widetilde{\varphi_{a,\lambda}}^{a})dV_{g_{a
}}
+
o_{\frac{1}{\lambda}}(\frac{1}{\lambda^{n-2\gamma}}),
\end{split}
\end{equation*}
provided we choose \;$\varepsilon\sim \lambda^{-\frac{1}{k}}$\; and \;$k\gg1$. Recalling \eqref{eq:relationcy}, from 
\eqref{properties_of_tilde_phi} we then find 
\begin{equation*}\begin{split}
\langle \varphi_{a,\lambda}  , &\varphi_{a,\lambda}\rangle_{P_{h_{a}}^\gamma}
= 
c_{n, 2}^\gamma
-
d_\gamma^*\frac{M_{a}}{\lambda^{\frac{n-2\gamma}{2}}}\int_{M}\eta^{2}_{a,\varepsilon}
\delta_{a,\lambda}^{\frac{n+2\gamma}{n-2\gamma}}dV_{h_{a}}
+
o_{\frac{1}{\lambda}}(\frac{1}{\lambda^{n-2\gamma}})
\end{split}
\end{equation*}
and arguing as for \eqref{super_quadratic_bubble_energy_expansion} we conclude 
\begin{equation*}\begin{split}
\langle \varphi_{a,\lambda}  , &\varphi_{a,\lambda}\rangle_{P_{h_{a}}^\gamma}
= 
c_{n, 2}^\gamma
-c_{n, \gamma}^*
\frac{M_{a}}{\lambda^{n-2\gamma}}
+
o_{\frac{1}{\lambda}}(\frac{1}{\lambda^{n-2\gamma}}),
\end{split}
\end{equation*}
cf.  \eqref{defc3}. This shows (i) and to see (ii) we have, again arguing as for \eqref{super_quadratic_bubble_energy_expansion}, that
\begin{equation*}
\begin{split}
\int_{M}\varphi_{a,\lambda}^{\frac{2n}{n-2\gamma}}dV_{h_{a}}
= &
\int_{B_{\varepsilon}(a)}\delta_{a,\lambda}^{\frac{2n}{n-2\gamma}}dV_{h_{a}
}
+
\int_{M\setminus 
B_{\varepsilon}(a)}O((\frac{1}{\lambda^{\frac{n-2\gamma}{2}}\varepsilon^{
n-2\gamma}})^{\frac{2n}{n-2\gamma}})dV_{h_{a}} \\
= &
\int_{\R^{n}}\delta_{a,\lambda}^{\frac{2n}{n-2\gamma}}dV_{h_{a}}
+
\int_{\R^{n}\setminus 
B_{\varepsilon}(a)}\delta_{a,\lambda}^{\frac{2n}{n-2\gamma}}
+
O(\frac{\varepsilon^{-2n}}{\lambda^{n}}) 
= 
c_{n, 1}^\gamma+o_{\frac{1}{\lambda}}(\lambda^{2\gamma-n})
\end{split}
\end{equation*}
provided \;$\varepsilon\sim \lambda^{-\frac{1}{k}}$\; for $\;k\gg 1$\; sufficiently 
large. The assertion then follows immediately from 
conformal covariance properties \eqref{eq:confinv} and \eqref{h_a_metric}.
\end{pf}
\begin{df}\label{fracmass}
In analogy to the classical Yamabe problem we define the fractional mass  map\;$\mathcal{M}_\gamma$\; as 
$$
\mathcal{M}_{\gamma}(a)=M_a,
$$
where \;$M_a$\; is given by Lemma \ref{lem_bubble comparsion}, and the fractional minimal mass as 
\;$ \mathcal{M}^{\min}_\gamma=\min_{a\in M}\mathcal{M}_{\gamma}(a)$.
\end{df}

\noindent
\section{Variational and algebraic topological argument}\label{varagtop}
In this section we present the proof of Theorem \ref{eq:thm} and therefore assume that we are under the assumptions of Theorem \ref{eq:thm}. 
Furthermore,  because of Remark \ref{eq:half} and the work of Gonzalez-Qing\cite{gq}, 
we assume also that \;$\gamma\neq \frac{1}{2}$\; and \; $\mathcal{Y}^\gamma(M, [h])>0$, and hence \;$G_h>0$\; for the Green's function.  
With the above agreement we carry out the variational and algebraic topological argument for existence. 
We point out that the algebraic topological argument of Bahri-Coron\cite{bc} has been used \cite{gam1}, \cite{gam2}, \cite{ould1}, \cite{martndia2} and \cite{nss}. 
Hence we will omit some standard poofs and encourage readers to find them in \cite{martndia2}.
\subsection{Variational principle in a non-compact setting}\label{var}
In this subsection we extend the classical variational principle to this non-compact setting.  
Clearly our bubbles \;$v_{a, \l}$\; can be used to replace the standard bubbles in the analysis of diverging Palais-Smale sequences of 
\;$\mathcal{E}^\gamma_{h}$\; 
in \cite{fanggonz} with 
\;$\mathcal{E}^{\gamma}_{h}$\; 
as in \eqref{eq:fracfunc}. 
Relying on this fact, we derive a deformation lemma that takes into account 
the bubbling phenomena via our bubbles \;$v_{a, \l}$, cf \eqref{dfbua}. 

To do so, we first  define for \;$p\in \N$
\begin{equation}\label{eq:dfwp}
W_p=W_{p, h}^\gamma:=\{u\in W^{\gamma, 2}_+(M, h)\;:\;\mathcal{E}^{\gamma}_{h}(u)\leq (p+1)^{\frac{2\gamma}{n}}\mathcal{Y}_{\gamma}(S^n)\},
\end{equation} 

where \;$\mathcal{Y}_{\gamma}(S^n)$\; is as in \eqref{eq:yamabehsphere}.
\vspace{6pt}

\noindent
Next  we introduce the {\em neighborhood of potential critical points at infinity} of  \;$\mathcal{E}^\gamma_h$, 
which depends on a universal constant \;$\nu_0>1$\; determined by Proposition \ref{eq:baryest} below.
\begin{df}
For \;$p\in \N^*$\; and \;$0<\varepsilon\leq \varepsilon_0$ (for some small \;$\varepsilon_{0}>0$) we call  
\begin{equation*}
\begin{split}
V(p, \varepsilon)=\{u\in W^{
\gamma, 2}(M,  h): &\;\;\exists a_1, \cdots, a_{p}\in M,\;\;\alpha_1, \cdots, \alpha_{p}>0, \;\l_1, \cdots,\l_{p}\geq \frac{1}{\varepsilon}, \\ &\;\; \Vert u-\sum_{i=1}^{p}\alpha_i\varphi_{a_i, \l_i}\Vert\leq \varepsilon, \;\;\;\;\frac{\alpha_i}{\alpha_j}\leq \nu_0 \;\;\text{ and }\;\;\varepsilon_{i, j}\leq \varepsilon, \;\;i\neq j=1, \cdots, p\},
\end{split}
\end{equation*}
the \;$(p, \varepsilon)$-neighborhood of potential critical points at infinity of \;$\mathcal{E}^\gamma_{h}$,
where \;$||\cdot||$\; denotes the norm associated to the scalar product \;$\left<\cdot, \cdot\right>_{P^{\gamma}_h}$\; defined by \eqref{fracscal}.
\end{df}

\noindent
Concerning the sets \;$V(p, \varepsilon)$, we have
\begin{lem}\label{lem_minimizer}
For every \;$p\in \N^*$ there exist \;$0<\varepsilon_p\leq\varepsilon_0$\; and \;$C_{p}>1$\; such that for every \;$0<\varepsilon\leq \varepsilon_p$ 
\begin{equation}\label{eq:mini}
\begin{cases}
\forall u\in V(p, \varepsilon)\;\; \text{the minimization problem}\;\;\min_{B_{2C_{p}\varepsilon}^{p}}\Vert u-\sum_{i=1}^{p}\alpha_i\varphi_{a_i, \l_i}\Vert \\
\text{has a solution \;$(\bar \alpha,  A, \bar \l)\in B_{C_p\varepsilon}^{p}$, which is unique up to permutations,}
\end{cases}
\end{equation}
where \;$B^{p}_{\gamma}$\; for \;$\gamma>0$\;  is defined by
\begin{equation*}
\begin{split}
B_{\gamma}^{p}
:=
\{
(& \bar \alpha=(\alpha_1,\cdots, \alpha_p),  A=(a_1, \cdots, a_p), \bar \l=(\l_1, \cdots, \l_p))\in \R^{p}_+\times M^p\times (0, +\infty)^{p}\;  
\\ &
\text{such that}\;\;\l_i\geq \frac{1}{\gamma} \;\;\text{ for } \;\; i=1, \cdots, p, \;\;\frac{\alpha_i}{\alpha_j}\leq \nu_0\;\; \text{and}\;\ \varepsilon_{i, j}\leq \gamma\;\;\text{for}\;\; i\neq j=1, \cdots, p
\}.
\end{split}
\end{equation*} 
\end{lem}
\noindent
Furthermore, for every \;$p\in \N^*$\; and \;$0<\varepsilon\leq \varepsilon_p$ ,  we consider the selection map 
\begin{equation}\label{eq:select}
s_{p}: V(p, \varepsilon)\longrightarrow ( M)^p/\sigma_p
:
u\longrightarrow s_{p}(u)=[A].
\end{equation}
Here $\sigma_{p}$ denotes the permutation group acting on $M^{p}$ 
and $A=(a_{1},\ldots,a_{p})\in M^{p}$ is derived  from a minimizer  
$\sum_{i=1}^{p} \alpha_{i}\varphi_{a_{i},\lambda_{i}}$ associated to 
$u\in V(p,\varepsilon)$ by virtue of Lemma \ref{lem_minimizer}. Since Lemma \ref{lem_minimizer} 
not only assures existence of a minimizer, but also its uniqueness up to permutations of indices, 
the class $[A]=A/\sigma_{p}$ is uniquely determined by $u$ and $s_{p}$ therefore  well defined.  
\vspace{6pt}
 
\noindent
Now, having introduced the neighborhoods of potential critical points at infinity of  \;$\mathcal{E}^\gamma_{h}$, 
we are ready to state a deformation lemma, 
which follows from the same arguments as for its counterparts in classical application of the barycenter technique for existence of Bahri-Coron\cite{bc} 
and the fact that the \;$v_{a, \l}$\; can replace the standard bubbles in the analysis of diverging Palais-Smale sequences of \; $\mathcal{E}^\gamma_{h}$. 
Indeed we have the following result, which was the aim of this subsection.

\begin{lem}\label{eq:classicdeform}
Assuming that \;$\mathcal{E}^\gamma_{h}$\; has no critical points, 
then for every \;$p\in \N^*$, up to taking \;$\varepsilon_p$\; given by \eqref{eq:mini} smaller, we have that for every \;$0<\varepsilon\leq \varepsilon_p$\; 
the topological pair \;$(W_p,\; W_{p-1})$\; retracts by deformation onto \;$(W_{p-1}\cup A_p, \;W_{p-1})$\; 
with 
\;$V(p, \;\tilde \varepsilon)\subset A_p\subset V(p, \;\varepsilon)$, 
where \;$0<\tilde \varepsilon<\frac{\varepsilon}{4}$\; is a very small positive real number and depends on \;$\varepsilon$.
\end{lem}
\subsection{Algebraic topological argument of Bahri-Coron}\label{agtop}
In this subsection we present the algebraic topological argument or barycenter technique of Bahri-Coron\cite{bc}  for existence to prove Theorem \ref{eq:thm}. 
To that end we start establishing  sharp \;$\mathcal{E}_h^\gamma$-energy estimates for \;$v_{a, \l}$. 
As in \cite{nss}, we will follow the presentation of the barycenter technique in \cite{martndia2} 
and hence omit some standard proofs and direct the readers to \cite{martndia2} for details. Recalling \eqref{eq:fracfunc}, we have
\begin{lem}
\label{sharpenergy} 
With \;$c_{n, \gamma}^*$ as in \eqref{defc3} there holds
$$
\mathcal{E}_h^\gamma(v_{a, \l})=\mathcal{Y}_\gamma(S^n)\left(1- c_{n, \gamma}^*\frac{\mathcal{M}_{\gamma}(a)}{\l^{n-2\gamma}}\right)\left(1+o_{\frac{1}{\l}}(1)\right).
$$ 
\end{lem}
\begin{pf}
The proof is a direct application of Lemma \ref{selfactestv}.
\end{pf}

\noindent
Next  we  define for \;$p\in \N^*$ and $\lambda >0$
\begin{equation*}
f_p(\l)
\; : \;
B_p(M)\longrightarrow W^{\gamma, 2}_+(M)
\;:\;
\sigma=\sum_{i=1}^p\alpha_i\d_{a_i}
\longrightarrow
f_p(\l)(\sigma)
=
\sum_{i=1}^p\alpha_i v_{a_i, \l}
\end{equation*}
with \;$B_p(M)$\; as in \eqref{dfbp}. Using Lemmas \ref{interactest} and \ref{sharpenergy}, we will prove that at \textit{true infinity} the quantity 
\begin{equation*}
\D^{\l}_p\mathcal{E}^{h}_{\gamma}(\sigma)
:=
\l^{n-2\gamma}\left(\frac{\mathcal{E}_\gamma^{\hat h}(f_p(\sigma))}{p^{\frac{2\gamma}{n}}{}\mathcal{Y}_\gamma(S^n)}-1\right), 
\;\;
\sigma=\sum_{i=1}^p\alpha_i\d_{a_i}\in B_p(M)
\end{equation*}
has at most a linear growth in \;$p$\; with a negative slope and this uniformly in \;$\sigma \in B_p(M)$\; and \;$\l$\; large, namely at infinity
\begin{equation}\label{eq:ling1}
\sup_{\l>>1, \;\sigma\in B_p(M)}\D^{\l}_p\mathcal{E}^{\hat h}_{\gamma}(\sigma)\leq -cp+C.
\end{equation}
About estimate \eqref{eq:ling1} we precisely show the following proposition.
\begin{pro}\label{eq:baryest}
There exists  \;$\nu_0>1$\; such that for every \;$p\in \N^*,\;p\geq 2$\; and every \;$0<\varepsilon\leq \varepsilon_0$ 
there exists \;$\l_p=\l_p(\varepsilon)$\; such that for every \;$\l\geq \l_p$\; and for every 
\;$\sigma=\sum_{i=1}^p\alpha_i\delta_{a_i}\in B_p(M)$\; we have,
\begin{enumerate}[label=(\roman*)]
 \item
if \;$\sum_{i\neq j}\varepsilon_{i, j}> \varepsilon$\; or there exist \;$i_0\neq j_0$\; such that \;$\frac{\alpha_{i_0}}{\alpha_{j_0}}>\nu_0$, then
 $$
\mathcal{E}^\gamma_{h}(f_p(\l)(\sigma))\leq p^{\frac{2\gamma}{n}}\mathcal{Y}_{\gamma}(S^n).
$$ 
 \item
if \;$\sum_{i\neq j}\varepsilon_{i, j}\leq \varepsilon$\; and for every \;$i\neq j$\; we have \;$\frac{\alpha_{i}}{\alpha_j}\leq\nu_0$, then
$$
\mathcal{E}^\gamma_{h}(f_p(\l)(\sigma))\leq p^{\frac{2\gamma}{n}}\mathcal{Y}_{\gamma}(S^n)\left(1- c_{n, \gamma}^*\frac{\mathcal{M}^{\min}_{\gamma}}{\l^{n-2\gamma}}-c_{h}^\gamma\frac{(p-1)}{\l^{n-2\gamma}}\right),
$$
where \;$c_{h}^{\gamma}:=\frac{c^\gamma_{n, 3}}{4c^1_{n, 1}}\min_{(M^2)^*}G_{h}$,\; $\mathcal{M}^{\min}_{\gamma}$\; is as in Definition \ref{fracmass} 
and \;$(M^2)^*$\; defined in \eqref{M_squred_star}.
\end{enumerate}

\end{pro}
\begin{pf}
The proof is the same as the one of Proposition 3.1 in \cite{martndia2} using Lemmas \ref{interactest}, \ref{selfactestv},  \ref{sharpenergy}, \ref{lem_higher_interaction_estimates}.
\end{pf}

\noindent
Now we start transporting the topology of the manifold \;$M$\; into  sublevels of the Euler-Lagrange functional 
\;$\mathcal{E}_h^\gamma$\; by bubbling via \;$v_{a, \l}$. 
Recalling \eqref{orientation_classes} and \eqref{eq:dfwp}, we have
\begin{lem}\label{eq:nontrivialf1}
Assuming that \;$\mathcal{E}_{ h}^{\gamma}$\; has no critical points and \:$0<\varepsilon\leq  \varepsilon_1$, then up to taking \;$\varepsilon_1$\; smaller and \;$\l_1$\; larger, 
we have that for every \;$\l\geq \l_1$
$$
f_1(\l)\; : \;(B_1(M),\; B_0(M))\longrightarrow (W_1, \;W_0)
$$
is well defined and satisfies
$$
(f_1(\l))_*(w_1)\neq 0\;\;\text{in}\;\;H_{n}(W_1, \;W_0). 
$$
\end{lem}
\begin{pf}
The proof follows from the same arguments as the ones used in the proof of  Lemma 4.2  in \cite{martndia2} by using the selection map \;$s_1$\;(see \eqref{eq:select}), Lemma \ref{eq:classicdeform} and Lemma \ref{sharpenergy}.
\end{pf}

\noindent
Next  we use the previous lemma and {\em pile up masses}  by bubbling  via \;$v_{a, \l}$\, in a recursive way. Still recalling \eqref{orientation_classes} we have 
\begin{lem}\label{eq:nontrivialrecursive}
Assuming that \;$\mathcal{E}^{\gamma}_{ h}$\; has no critical points and \;$0<\varepsilon\leq  \varepsilon_{p+1}$, then up to taking \;$\varepsilon_{p+1}$\; smaller 
and \;$\l_p$\; and \;$\l_{p+1}$\; larger, we have that for every \;$\l\geq \max\{\l_p, \l_{p+1}\}$
$$
f_{p+1}(\l): (B_{p+1}(M),\; B_{p}(M))\longrightarrow (W_{p+1}, \;W_{p})
$$
and 
$$
f_p(\l): (B_p(M), \;B_{p-1}(M))\longrightarrow (W_p, \; W_{p-1})
$$
are well defined and
$$(f_p(\l))_*(w_p)\neq 0\;\; \text{in}\;\; \;\;H_{np+p-1}(W_p, \;W_{p-1})$$ implies
$$(f_{p+1}(\l))_*(w_{p+1})\neq 0\;\; \text{in} \;\;H_{n(p+1)+p}(W_{p+1}, \;W_{p}).$$
\end{lem}
\begin{pf}
The proof follows from the same arguments as the ones used in the proof of Lemma 4.3  in \cite{martndia2}, by using the selection map \;$s_p$\; (see \eqref{eq:select}), Lemma \ref{eq:classicdeform} and Proposition \ref{eq:baryest}.
\end{pf}

\noindent
Finally we use the strength of Proposition \ref{eq:baryest} - namely point (ii) - to give a criterion ensuring that the recursive process of {{\em piling up} masses via Lemma \ref{eq:nontrivialrecursive} must stop at least after a very large number of steps. 
\begin{lem}\label{eq:largep}
Setting \;$$p^{\gamma}_{h}:=[1-c_{n, \gamma}^*\frac{\mathcal{M}^{\min}_{\gamma}}{c^{\gamma}_{h}} ]+2$$ and recalling \eqref{eq:dfwp}, 
then \;$\forall\;\;0<\varepsilon\leq\varepsilon_0$\; and \;$\forall\l\geq\l_{p_{\gamma}^{h}}$ there holds  
\;$f_{p^{\gamma}_{h}}(\l)(B_{p^{\gamma}_{h}}(M))\subset W_{p_h^\gamma-1}$.
\end{lem}
\begin{pf}
The proof is a direct application of Proposition \ref{eq:baryest}.
\end{pf}

\noindent
\begin{pfn} {of Theorem \ref{eq:thm}} 
It follows by a contradiction argument from Lemma \ref{eq:nontrivialf1} - Lemma \ref{eq:largep}.
\end{pfn}

\section{Appendix}\label{appendix}
In this section we provide the proofs of Lemmas \ref{lem_standard interaction on Rn}, \ref{lem_higher_interaction_estimates} and of some technical estimates. 

\bigskip

\begin{pfn}{ of Lemma \ref{lem_standard interaction on Rn}.}
Let us write \;$\delta_{i}=\delta_{a_{i},\lambda_{i}},\; \delta_{j}=\delta_{a_{j}, \lambda_{j}}$\; for simplicity.
\begin{enumerate}[label=(\roman*)]
 \item 
 Due to
\begin{equation*}\begin{split}
\underset{\R^{n}}{\int} \delta_{i}^{\frac{n+2\gamma}{n-2\gamma}}\delta_{j}
= &
\frac{1}{c_{n, \gamma}}\underset{\R^{n}}{\int} (-\Delta)^{\gamma}\delta_{i}\delta_{j}=\frac{1}{c_{n, \gamma}}\underset{\R^{n}}{\int} 
\delta_{i}(-\Delta)^{\gamma}\delta_{j}
=
\underset{\R^{n}}{\int} \delta_{i}\delta_{j}^{\frac{n+2\gamma}{n-2\gamma}}
\end{split}
\end{equation*} 
we may assume \;$ \frac{1}{\lambda_{i}} \leq  \frac{1}{\lambda_{j}}$. We then have with \;$r=\vert x \vert$\;
\begin{equation*}\begin{split}
\underset{\R^{n}}{\int} \delta_{i}^{\frac{n+2\gamma}{n-2\gamma}} \delta_{j}
= &
\underset{\R^{n}}{\int}
(\frac{1}{1+r^{2}})^{\frac{n+2\gamma}{2}}
(\frac{1}{\frac{\lambda_{i} }{ \lambda_{j} }+\lambda_{i}\lambda_{j}\vert 
\frac{x}{\lambda_{i}}+a_{i}-a_{j}\vert^{2}})^{\frac{n-2\gamma}{2}}dx.
\end{split}\end{equation*}
Since \;$\eps_{i,j}^{\frac{2}{2\gamma-n}}
\sim 
\lambda_{i}\lambda_{j}\vert a_{i}-a_{j}\vert^{2}
\;\text{ or }\;
\eps_{i,j}^{\frac{2}{2\gamma-n}}
\sim
\frac{\lambda_{i}}{\lambda_{j}},
$\; we  may expand on
\begin{equation*}\begin{split}
\mathcal{A}
=
\left[
\vert  \frac{x}{\lambda_{i}}\vert \leq \epsilon\vert a_{i}-a_{j}\vert
\right]
\cup
\left[
\vert  \frac{x}{\lambda_{i}}\vert \leq \epsilon  \frac{1}{\lambda_{j}} 
\right]
\end{split}\end{equation*}
for \;$\epsilon>0$\; sufficiently small
\begin{equation*}\begin{split}
(
\frac{\lambda_{i} }{ \lambda_{j} }
& +
\lambda_{i}\lambda_{j}\vert \frac{x}{\lambda_{i}} 
+a_{i}-a_{j}\vert^{2})
)^{\frac{2\gamma-n}{2}} \\
= & 
(
\frac{\lambda_{i} }{ \lambda_{j} }
+
\lambda_{i}\lambda_{j}
\vert a_{i}-a_{j}\vert^{2}
)^{\frac{2\gamma-n}{2}} 
+
(2\gamma-n)
\frac
{
\langle a_{i}-a_{j}, \lambda_{j}x\rangle
+
O(\frac{ \lambda_{j} }{ \lambda_{i} }\vert x\vert^{2})
}
{
(
\frac{\lambda_{i} }{ \lambda_{j} }
+
\lambda_{i}\lambda_{j}\vert a_{i}-a_{j}\vert^{2}
)^{\frac{n+2-2\gamma}{2}}
}.
\end{split}\end{equation*}
Thus  \;$
\underset{\R^{n}}{\int} \delta_{i}^{\frac{n+2\gamma}{n-2\gamma}}\delta_{j}
= 
\sum_{k=1}^{4}I_{k}
$\; with 
\begin{enumerate}[label=(\roman*)]
 \item \quad
$
I_{1}
= 
\frac
{1}
{(\frac{\lambda_{i} }{ \lambda_{j} }
+
\lambda_{i}\lambda_{j}
\vert a_{i}-a_{j}\vert^{2}
)^{\frac{n-2\gamma}{2}}
}
\int_{\mathcal{A}}\frac{1}{(1+r^{2})^{\frac{n+2\gamma}{2}}}
$
 \item \quad
$
I_{2}
=
\frac
{2\gamma-n}
{(\frac{\lambda_{i} }{ \lambda_{j} }
+
\lambda_{i}\lambda_{j}
\vert a_{i}-a_{j}\vert^{2}
)^{\frac{n+2-2\gamma}{2}}
}
\int_{\mathcal{A}} \frac{\langle a_{i}-a_{j}, 
\lambda_{j}x\rangle}{(1+r^{2})^{\frac{n+2\gamma}{2}}}dx
$
 \item \quad
$
I_{3}
= 
\frac
{1}
{(\frac{\lambda_{i} }{ \lambda_{j} }
+
\lambda_{i}\lambda_{j}
\vert a_{i}-a_{j}\vert^{2}
)^{\frac{n+2-2\gamma}{2}}
}
\int_{\mathcal{A}}
\frac{O(\frac{ \lambda_{j} }{ \lambda_{i} }r^{2})}{(1+r^{2})^{\frac{n+2\gamma}{2}}}
$
 \item \quad
$
I_{4}
= 
\underset{\mathcal{A}^{c}}{\int} (\frac{ 1}{1+r^{2}})^{\frac{n+2\gamma}{2}}
(\frac{1}{\frac{\lambda_{i} }{ \lambda_{j} }+\lambda_{i}\lambda_{j}\vert 
\frac{x}{\lambda_{i}}+a_{i}-a_{j}\vert^{2}})^{\frac{n-2\gamma}{2}}dx.
$
\end{enumerate}
We get
\begin{equation*}\begin{split}
I_{1}
= &
\frac
{b_{1}}
{(\frac{\lambda_{i} }{ \lambda_{j} }
+
\lambda_{i}\lambda_{j}
\vert a_{i}-a_{j}\vert^{2}
)^{\frac{n-2\gamma}{2}}
}
+
O((\frac{\lambda_{j}}{\lambda_{i}})^{\gamma}\eps_{i,j}^{\frac{n}{n-2\gamma}})
\end{split}\end{equation*}
with \;$b_{1}=\underset{\R^{n}}{\int} (\frac{1}{1+r^{2}})^{\frac{n+2\gamma}{2}}=c_{n,3}^{\gamma}$, \; $I_{2}=0$\;  and \;$
I_{3}
= 
O((\frac{\lambda_{j}}{\lambda_{i}})^{\gamma}\eps_{i,j}^{\frac{n}{n-2\gamma}})
$. Moreover  
\begin{equation*}
\begin{split}
I_{4}
\leq 
\frac{C}{(\frac{\lambda_{i}}{\lambda_{j}})^{\frac{n+2\gamma}{2}}}
=
O((\frac{\lambda_{j}}{\lambda_{i}})^{\gamma}\eps_{i,j}^{\frac{n}{n-2\gamma}})\;\;
\text{ in case }\;\;\eps_{i,j}^{\frac{2}{2\gamma-n}}\sim 
\frac{\lambda_{i}}{\lambda_{j}}.
\end{split}
\end{equation*} 
Otherwise we may assume \;$\frac{\lambda_{j}}{\lambda_{i}}\leq 
\frac{\lambda_{i}}{\lambda_{j}}\ll \lambda_{i}\lambda_{j}\vert 
a_{i}-a_{j}\vert^{2}$\; and decompose 
\begin{equation*}
\begin{split}
\mathcal{A}^{c}
\subseteq &
\mathcal{B}_{1}\cup \mathcal{B}_{2},
\end{split}
\end{equation*}
where for a sufficiently large constant \;$E>0$
\begin{equation*}
\begin{split}
\mathcal{B}_{1}
= &
\left[\epsilon \vert a_{i}-a_{j}\vert
\leq 
\vert \frac{x}{\lambda_{i}}\vert 
\leq 
E\vert a_{i}-a_{j}\vert\right
] 
\;\text{ and }\;
\mathcal{B}_{2}
= 
\left[
E\vert a_{i}-a_{j}\vert
\leq 
\vert \frac{x}{\lambda_{i}}\vert 
<\infty
\right].
\end{split}
\end{equation*} 
We then may estimate
\begin{equation*}
\begin{split}
I_{4}^{1}
= &
\underset{\mathcal{B}_{1}}{\int} (\frac{1}{1+r^{2}})^{\frac{n+2\gamma}{2}}
(\frac{1}{\frac{\lambda_{i} }{ \lambda_{j} }+\lambda_{i}\lambda_{j}
\vert \frac{x}{\lambda_{i}}+a_{i}-a_{j}\vert^{2}})^{\frac{n-2\gamma}{2}} dx\\
\leq &
\frac{C(\frac{\lambda_{i}}{\lambda_{j}})^{\frac{n+2\gamma}{2}}}{(1+\lambda_{i}^{
2}\vert a_{i}-a_{j}\vert^{2})^{\frac{n+2\gamma}{2}}} 
\underset{[\vert \frac{x}{\lambda_{j}}\vert \leq E\vert 
a_{i}-a_{j}\vert]}{\int} 
(\frac{1}{1+\lambda_{j}^{2}\vert 
\frac{x}{\lambda_{j}}+a_{i}-a_{j}\vert^{2}})^{\frac{n-2\gamma}{2}}dx
\\
\leq &
\frac{C}{(\frac{\lambda_{j}}{\lambda_{i}}+\lambda_{i}\lambda_{j}\vert 
a_{i}-a_{j}\vert^{2})^{\frac{n+2\gamma}{2}}}
\underset{B_{C\lambda_{j}\vert a_{i}-a_{j}\vert}(0)}{\int} 
(\frac{1}{1+r^{2}})^{\frac{n-2\gamma}{2}}\\
\leq &
\frac{C(\lambda_{j}^{2}\vert 
a_{i}-a_{j}\vert^{2})^{\gamma}}{(\frac{\lambda_{j}}{\lambda_{i}}+\lambda_{i}
\lambda_{j}\vert a_{i}-a_{j}\vert^{2})^{\frac{n+2\gamma}{2}}}
=
O((\frac{\lambda_{j}}{\lambda_{i}})^{\gamma}\eps_{i,j}^{\frac{n}{n-2\gamma}})
\end{split}
\end{equation*} 
and 
\begin{equation*}
\begin{split}
I_{4,2}
= &
\underset{\mathcal{B}_{2}}{\int} (\frac{ 1}{1+r^{2}})^{\frac{n+2\gamma}{2}}
(\frac{1}{\frac{\lambda_{i} }{ \lambda_{j} }+\lambda_{i}\lambda_{j}\vert  
\frac{x}{\lambda_{i}}+a_{i}-a_{j}\vert^{2}})^{\frac{n-2\gamma}{2}}dx \\
\leq &
\frac{C}{(\frac{\lambda_{i}}{\lambda_{j}}+\lambda_{i}\lambda_{j}\vert 
a_{i}-a_{j}\vert^{2})^{\frac{n-2\gamma}{2}}} 
\underset{[\vert x \vert \geq \lambda_{i}\vert a_{i}-a_{j}\vert]}{\int}
(\frac{1}{1+r^{2}})^{\frac{n+2\gamma}{2}} 
= 
O((\frac{\lambda_{j}}{\lambda_{i}})^{\gamma}\eps_{i,j}^{\frac{n}{n-2\gamma}}),
\end{split}
\end{equation*} 
so  \;$
I_{4}\leq 
I_{4}^{1}+I_{4}^{2}=O((\frac{\lambda_{j}}{\lambda_{i}})^{\gamma}\eps_{i,j}^{
\frac{n}{n-2\gamma}}).
$\;
Collecting terms the claim follows.

\item It is sufficient to consider the case \;$k=j$, since 
\begin{equation*}
\begin{split}
\lambda_{j}\partial_{\lambda_{j}}\underset{\R^{n}}{\int} 
\delta_{i}\delta_{j}^{\frac{n+2\gamma}{n-2\gamma}}
=
\lambda_{j}\partial_{\lambda_{j}}\underset{\R^{n}}{\int} 
\delta_{i}\delta_{j}^{\frac{n+2\gamma}{n-2\gamma}}.
\end{split}
\end{equation*} 
First we deal with the case \;$\frac{1}{\lambda_{i}}\leq \frac{1}{\lambda_{j}}$.
We have
\begin{equation*}\begin{split}
-\lambda_{j}\underset{\R^{n}}{\int} 
\delta_{i}^{\frac{n+2\gamma}{n-2\gamma}}\partial_{\lambda_{j}}\delta_{j}
= &
\frac{n-2\gamma}{2}
\underset{\R^{n}}{\int}
(\frac{ \lambda_{i} }{ 1+\lambda_{i}^{2}\vert 
x-a_{i}\vert^{2}})^{\frac{n+2\gamma}{2}} 
(\frac{ \lambda_{j} }{ 1+\lambda_{j}^{2}\vert 
x-a_{j}\vert^{2}})^{\frac{n-2\gamma}{2}}
\frac{\lambda_{j}^{2}\vert x-a_{j}\vert^{2}-1}{\lambda_{j}^{2}\vert 
x-a_{j}\vert^{2}+1}dx,
\end{split}\end{equation*}
whence
\begin{equation*}\begin{split}
-\lambda_{j}\underset{\R^{n}}{\int}  
\delta_{i}^{\frac{n+2\gamma}{n-2\gamma}}\partial_{\lambda_{j}}\delta_{j} 
= &
\frac{n-2\gamma}{2}
\underset{\R^{n}}{\int}
(\frac{1}{1+r^{2}})^{\frac{n+2\gamma}{2}} 
\frac{\lambda_{j}^{2}\vert\frac{x}{\lambda_{i}}+a_{i}-a_{j}\vert^{2}-1}
{\lambda_{j}^{2}\vert\frac{x}{\lambda_{i}}+a_{i}-a_{j}\vert^{2}+1}
(\frac{1}{\frac{\lambda_{i} }{ \lambda_{j} 
}+\lambda_{i}\lambda_{j}\vert\frac{x}{\lambda_{i}}+a_{i}-a_{j}\vert^{2}})^{\frac
{n-2\gamma}{2}} dx
.
\end{split}\end{equation*}
Since \;$
\eps_{i,j}^{\frac{2}{2-n}}
\sim 
\lambda_{i}\lambda_{j}\vert a_{i}-a_{j}\vert^{2}
\;\;\text{ or }\;\;
\eps_{i,j}^{\frac{2}{2-n}}
\sim
\frac{\lambda_{i}}{\lambda_{j}},
$
\; we may expand on
\begin{equation*}\begin{split}
\mathcal{A}
= &
\left[
\vert  \frac{x}{\lambda_{i}}\vert \leq \epsilon\vert a_{i}-a_{j}\vert
\right]
\cup
\left[\vert  \frac{x}{\lambda_{i}}\vert \leq \epsilon  \frac{1}{\lambda_{j}} 
\right]
\end{split}\end{equation*}
for \;$\epsilon>0$\; sufficiently small
\begin{equation*}\begin{split}
& 
\frac{1}
{
(\frac{\lambda_{i} }{ \lambda_{j} }
+
\lambda_{i}\lambda_{j}\vert\frac{x}{\lambda_{i}}+a_{i}-a_{j}\vert^{2})^{\frac{
n-2\gamma}{2}} }
\frac{\lambda_{j}^{2}\vert\frac{x}{\lambda_{i}}+a_{i}-a_{j}\vert^{2}-1}
{\lambda_{j}^{2}\vert\frac{x}{\lambda_{i}}+a_{i}-a_{j}\vert^{2}+1}
\\
& \quad = 
(
\frac{\lambda_{i} }{ \lambda_{j} }
+
\lambda_{i}\lambda_{j}
\vert a_{i}-a_{j}\vert^{2}
)^{\frac{2\gamma-n}{2}} 
\frac
{\lambda_{j}^{2}\vert a_{i}-a_{j}\vert^{2}-1}
{\lambda_{j}^{2}\vert a_{i}-a_{j}\vert^{2}+1}
\\
& \quad \;\; 
+
(2\gamma-n)
\frac
{
\langle a_{i}-a_{j}, \lambda_{j}x\rangle
}
{
(
\frac{\lambda_{i} }{ \lambda_{j} }
+
\lambda_{i}\lambda_{j}\vert a_{i}-a_{j}\vert^{2}
)^{\frac{n+2-2\gamma}{2}}
}
\frac
{\lambda_{j}^{2}\vert a_{i}-a_{j}\vert^{2}-1}
{\lambda_{j}^{2}\vert a_{i}-a_{j}\vert^{2}+1}
\\
& \quad \;\;+
\frac
{
4
}
{
(
\frac{\lambda_{i} }{ \lambda_{j} }
+
\lambda_{i}\lambda_{j}\vert a_{i}-a_{j}\vert^{2}
)^{\frac{n+2-2\gamma}{2}}
}
\frac
{\langle a_{i}-a_{j}, \lambda_{j}x\rangle}
{1+\lambda_{j}^{2}\vert a_{i}-a_{j}\vert^{2}}
+
\frac
{
O(\frac{ \lambda_{j} }{ \lambda_{i} }\vert x\vert^{2})
}
{
(
\frac{\lambda_{i} }{ \lambda_{j} }
+
\lambda_{i}\lambda_{j}\vert a_{i}-a_{j}\vert^{2}
)^{\frac{n+2-2\gamma}{2}}
}.
\end{split}\end{equation*}
From this we derive as before with \;$b_{2}=\frac{n-2\gamma}{2}c_{n,3}^{\gamma}$
\begin{equation*}\begin{split}
-\lambda_{j}\underset{\R^{n}}{\int}  
\delta_{i}^{\frac{n+2\gamma}{n-2\gamma}}\partial_{\lambda_{j}}\delta_{j} 
= &
b_{2}
\frac
{
\rfrac
{\lambda_{j}^{2}\vert a_{i}-a_{j}\vert^{2}-1}
{\lambda_{j}^{2}\vert a_{i}-a_{j}\vert^{2}+1}
}
{(\frac{\lambda_{i} }{ \lambda_{j} }
+
\lambda_{i}\lambda_{j}
\vert a_{i}-a_{j}\vert^{2}
)^{\frac{n-2\gamma}{2}}
} 
+
O((\frac{\lambda_{j}}{\lambda_{i}})^{\gamma}\varepsilon_{i,j}^{\frac{n}{
n-2\gamma}}),
\end{split}\end{equation*}
whence
\begin{equation*}\begin{split}
-\lambda_{j}\underset{\R^{n}}{\int} 
\delta_{i}^{\frac{n+2\gamma}{n-2\gamma}}\partial_{\lambda_{j}}\delta_{j}
= &
b_{2}
\frac
{
\lambda_{i}\lambda_{j}\vert a_{i}-a_{j}\vert^{2}-\frac{\lambda_{i}}{\lambda_{j}}
}
{(\frac{\lambda_{i} }{ \lambda_{j} }
+
\lambda_{i}\lambda_{j}
\vert a_{i}-a_{j}\vert^{2}
)^{\frac{n+2-2\gamma}{2}}
} 
+
O((\frac{\lambda_{j}}{\lambda_{i}})^{\gamma}\varepsilon_{i,j}^{\frac{n}{
n-2\gamma}}).
\end{split}\end{equation*}
We turn to the case  
\;$\frac{1}{\lambda_{i}}\geq \frac{1}{\lambda_{j}}
$. We then have
$
-\lambda_{j}\underset{\R^{n}}{\int} \delta_{i}^{\frac{n+2\gamma}{n-2\gamma}}\partial_{ 
\lambda_{j}}\delta_{j}
= 
-\lambda_{j}\underset{\R^{n}}{\int} \delta_{i}\partial_{\lambda_{j}} 
\delta_{j}^{\frac{n+2\gamma}{n-2\gamma}}
$ 
and
\begin{equation*}\begin{split}
-\lambda_{j}\underset{\R^{n}}{\int} \delta_{i} 
\partial_{\lambda_{j}}\delta_{j}^{\frac{n+2\gamma}{n-2\gamma}}
= &
\frac{n+2\gamma}{2}\underset{\R^{n}}{\int}
(\frac{ \lambda_{i} }{ 1+\lambda_{i}^{2}\vert 
x-a_{i}\vert^{2}})^{\frac{n-2\gamma}{2}} 
(\frac{ \lambda_{j} }{ 1+\lambda_{j}^{2}\vert 
x-a_{j}\vert^{2}})^{\frac{n+2\gamma}{2}}
\frac{\lambda_{j}^{2}\vert x-a_{j}\vert^{2}-1}{\lambda_{j}^{2}\vert 
x-a_{j}\vert^{2}+1}dx,
\end{split}\end{equation*}
whence
\begin{equation*}\begin{split}
-\lambda_{j}\underset{\R^{n}}{\int}  
\delta_{i}\partial_{\lambda_{j}}\delta_{j}^{\frac{n+2\gamma}{n-2\gamma}} 
= &
\frac{n+2\gamma}{2}
\underset{\R^{n}}{\int}
\frac
{r^{2}-1}
{r^{2}+1}
(\frac{1}{1+r^{2}})^{\frac{n+2\gamma}{2}}
\frac
{
1
}
{
(\frac{\lambda_{j}}{\lambda_{i}}+\lambda_{i}\lambda_{j}\vert\frac{x}{\lambda_{j}
}+a_{j}-a_{i}\vert^{2})^{\frac{n-2\gamma}{2}}
} dx.
\end{split}\end{equation*}
We may expand on
\;$
\mathcal{A}
= 
\left[
\vert  \frac{x}{\lambda_{j}}\leq \epsilon \vert a_{j}-a_{i}\vert
\right]
\cup
\left[\vert  \frac{x}{\lambda_{j}}\vert \leq \epsilon  \frac{1}{\lambda_{i}} 
\right]
$ \;
for \;$\epsilon>0$\; sufficiently small
\begin{equation*}\begin{split}
(
\frac{\lambda_{j} }{ \lambda_{i} }
& +
\lambda_{i}\lambda_{j}\vert 
\frac{x}{\lambda_{j}}+a_{j}-a_{i}\vert^{2})^{\frac{2\gamma-n}{2}} 
\\
= & 
(
\frac{\lambda_{j} }{ \lambda_{i} }
+
\lambda_{i}\lambda_{j}
\vert a_{j}-a_{i}\vert^{2}
)^{\frac{2\gamma-n}{2}} 
+
(2\gamma-n)
\frac
{
\langle a_{j}-a_{i}, \lambda_{i}x\rangle
+
O(\frac{ \lambda_{i} }{ \lambda_{j} }\vert x\vert^{2})
}
{
(
\frac{\lambda_{j} }{ \lambda_{i} }
+
\lambda_{i}\lambda_{j}\vert a_{j}-a_{i}\vert^{2}
)^{\frac{n+2-2\gamma}{2}}
}.
\end{split}\end{equation*}
As before we find with \;$
b_{2}
=
\frac{n-2\gamma}{2}\underset{\R^{n}}{\int} (\frac{1}{1+r^{2}})^{\frac{n+2\gamma}{2}}
=
\frac{n+2\gamma}{2}\underset{\R^{n}}{\int} 
\frac{r^{2}-1}{r^{2}+1}(\frac{1}{1+r^{2}})^{\frac{n+2\gamma}{2}}
$
\begin{equation*}\begin{split}
-\lambda_{j}\underset{\R^{n}}{\int} 
\delta_{i}\partial_{\lambda_{j}}\delta_{j}^{\frac{n+2\gamma}{n-2\gamma}}
= &
b_{2}\frac
{1}
{(\frac{\lambda_{i} }{ \lambda_{j} }
+
\lambda_{i}\lambda_{j}
\vert a_{j}-a_{i}\vert^{2}
)^{\frac{n-2\gamma}{2}}
} 
+
O((\frac{\lambda_{i}}{\lambda_{j}})^{\gamma}\eps_{i,j}^{\frac{n}{n-2\gamma}}).
\end{split}\end{equation*}
Collecting the results in cases  \;$\frac{1}{\lambda_{i}}\leq \frac{1}{\lambda_{j}}$\;
and \;$\frac{1}{\lambda_{j}}\leq \frac{1}{\lambda_{i}}$\; the claim follows.


\item 
By translation invariance and symmetry we may assume 
\;$
a_{i}=0, \;\;a_{j}=a \;\;\text{ and }\;\;\frac{1}{\lambda_{i}}\leq \frac{1}{\lambda_{j}}
$.
Then
\begin{equation*}
\begin{split}
\nabla_{a}  \underset{\R^{n}}{\int} \delta_{i}^{\frac{n+2\gamma}{n-2\gamma}}\delta_{j}
= &
\nabla_{a}\underset{\R^{n}}{\int}(\frac{\lambda_{i}}{1+\lambda_{i}^{2}r^{2}})^{\frac{n+2\gamma}{2}}
(\frac{\lambda_{j}}{1+\lambda_{j}^{2}\vert x+a\vert^{2}})^{\frac{n-2\gamma}{2}}dx
\\ = &
(2\gamma-n)
\underset{\R^{n}}{\int}(\frac{ \lambda_{i} }{ 1+\lambda_{i}^{2}r^{2}})^{\frac{n+2\gamma}{2}} 
(\frac{ \lambda_{j} }{ 1+\lambda_{j}^{2}\vert 
x+a\vert^{2}})^{\frac{n-2\gamma}{2}}
\frac
{\lambda_{j}^{2}( x+a)}
{1+\lambda_{j}^{2}\vert x+a\vert^{2}}  dx\\
= &
(2\gamma-n)
\underset{\R^{n}}{\int}(\frac{ 1 }{ 1+r^{2}})^{\frac{n+2\gamma}{2}} 
\frac{ \lambda_{i}\lambda_{j}( \frac{x}{\lambda_{i}}+a) }{ 
(\frac{\lambda_{i}}{\lambda_{j}}+\lambda_{i}\lambda_{j}\vert 
\frac{x}{\lambda_{i}}+a\vert^{2}))^{\frac{n+2-2\gamma}{2}}}dx.
\end{split}\end{equation*}
Since \;$
\eps_{i,j}^{\frac{2}{2\gamma-n}}
\sim 
\lambda_{i}\lambda_{j}\vert a\vert^{2}
\;\;\text{ or }\;\;
\eps_{i,j}^{\frac{2}{2\gamma-n}}
\sim
\frac{\lambda_{i}}{\lambda_{j}}
$\; we may expand on 
\begin{equation*}\begin{split}
\mathcal{A}
= &
[\vert  \frac{x}{\lambda_{i}}\vert \leq \epsilon\vert a\vert]
\cup
[\vert  \frac{x}{\lambda_{i}}\vert \leq \epsilon  \frac{1}{\lambda_{j}} ]
\end{split}\end{equation*}
for \;$\epsilon>0$\; sufficiently small 
\begin{equation*}
\begin{split}
& \hspace{-20pt}\frac{ \lambda_{i}\lambda_{j}(\frac{x}{\lambda_{i}}+a) }{ 
(\frac{\lambda_{i}}{\lambda_{j}}+\lambda_{i}\lambda_{j}\vert 
\frac{x}{\lambda_{i}}+a\vert^{2}))^{\frac{n+2-2\gamma}{2}}} \\
= &
\frac{ \lambda_{i}\lambda_{j}a }{ 
(\frac{\lambda_{i}}{\lambda_{j}}+\lambda_{i}\lambda_{j}\vert 
a\vert^{2}))^{\frac{n+2-2\gamma}{2}}} 
-(n+2-2\gamma)
\frac
{
\lambda_{i}\lambda_{j}a\lambda_{j}\langle x,a\rangle
}
{
(\frac{\lambda_{i}}{\lambda_{j}}+\lambda_{i}\lambda_{j}\vert 
a\vert^{2})^{\frac{n+4-2\gamma}{2}}
}
\\
& +
\frac{\lambda_{j}x }{ 
(\frac{\lambda_{i}}{\lambda_{j}}+\lambda_{i}\lambda_{j}\vert 
a\vert^{2}))^{\frac{n+2-2\gamma}{2}}} 
+
\frac{ \sqrt{\lambda_{i}\lambda_{j}} 
O(\frac{\lambda_{j}}{\lambda_{i}}\vert x 
\vert^{2}) }{ (\frac{\lambda_{i}}{\lambda_{j}}+\lambda_{i}\lambda_{j}\vert 
a\vert^{2}))^{\frac{n+3-2\gamma}{2}}}.
\end{split}
\end{equation*} 
Using radial symmetry we obtain with \;$b_{3}=(2\gamma-n)c_{n,3}^{\gamma}$
\begin{equation*}
\begin{split}
\nabla_{a}  \underset{\R^{n}}{\int} \delta_{i}^{\frac{n+2\gamma}{n-2\gamma}}\delta_{j} 
= &
b_{3}\frac{ \lambda_{i}\lambda_{j}a }{ 
(\frac{\lambda_{i}}{\lambda_{j}}+\lambda_{i}\lambda_{j}\vert 
a\vert^{2}))^{\frac{n+2-2\gamma}{2}}}
+
O((\frac{\lambda_{j}}{\lambda_{i}})^{\gamma}\sqrt{\lambda_{i}\lambda_{j}}
\varepsilon_{i,j}^{\frac{n+1}{n-2\gamma}})
\\
& +
(2\gamma-n)
\underset{\mathcal{A}^{C}}\int(\frac{ 1 }{ 1+r^{2}})^{\frac{n+2\gamma}{2}} 
\frac{ \lambda_{i}\lambda_{j}( \frac{x}{\lambda_{i}}+a) }{ 
(\frac{\lambda_{i}}{\lambda_{j}}+\lambda_{i}\lambda_{j}\vert 
\frac{x}{\lambda_{i}}+a\vert^{2})^{\frac{n+2-2\gamma}{2}}}dx.
\end{split}
\end{equation*} 
In case \;$\varepsilon_{i,j}^{\frac{2}{2\gamma-n}}\sim 
\frac{\lambda_{i}}{\lambda_{j}}$\; the last summand above is of order
\begin{equation*}
\begin{split}
O((\frac{\lambda_{j}}{\lambda_{i}})^{\gamma}\sqrt{\lambda_{i}\lambda_{j}}
\varepsilon_{i,j}^{\frac{n+1}{n-2\gamma}}).
\end{split}
\end{equation*} 
Thus we may assume \;$\frac{\lambda_{i}}{\lambda_{j}}\ll 
\lambda_{i}\lambda_{j}\vert a \vert^{2}$.  Letting \;$0<\theta<\frac{1}{4}$\; and
\begin{equation*}
\begin{split}
\mathcal{B}_{1}=\left[\vert \frac{x}{\lambda_{i}}+a\vert \leq \theta \vert 
a\vert\right]\subset \subset B^{c}_{\epsilon \lambda_{i}\vert a 
\vert}(0)=\mathcal{A}^{c}
\end{split}
\end{equation*} 
we may expand on \;$\mathcal{B}_{1}$
\begin{equation*}
\begin{split}
(\frac{1}{1+r^{2}})^{\frac{n+2\gamma}{2}}
= &
(\frac{1}{1+\lambda_{i}^{2}\vert \frac{x}{\lambda_{i}} +a -a 
\vert^{2}})^{\frac{n+2\gamma}{2}} 
=
(\frac{1}{1+\lambda_{i}^{2}\vert a \vert^{2}})^{\frac{n+2\gamma}{2}}
+
\frac{O(\lambda_{i}^{2}\vert \frac{x}{\lambda_{i}}+a\vert\vert a 
\vert)}{(1+\lambda_{i}^{2}\vert a \vert^{2})^{\frac{n+2+2\gamma}{2}}}
\end{split}
\end{equation*} 
and find using radial symmetry  
\begin{equation*}
\begin{split}
\vert \underset{\mathcal{B}_{1}}\int(\frac{ 1 }{ 
1+r^{2}})^{\frac{n+2\gamma}{2}} 
&
 \frac{ \lambda_{i}\lambda_{j}( \frac{x}{\lambda_{i}}+a) }{ 
(\frac{\lambda_{i}}{\lambda_{j}}+\lambda_{i}\lambda_{j}\vert 
\frac{x}{\lambda_{i}}+a\vert^{2}))^{\frac{n+2-2\gamma}{2}}} dx
\vert \\
\leq &
\frac{C}{(1+\lambda_{i}^{2}\vert a \vert^{2})^{\frac{n+2+2\gamma}{2}}}
\underset{\mathcal{B}_{1}}\int
\frac{\lambda_{i}^{3}\lambda_{j}\vert \frac{x}{\lambda_{i}}+a\vert^{2}\vert a 
\vert}
{(\frac{\lambda_{i}}{\lambda_{j}}+\lambda_{i}\lambda_{j}\vert\frac{x}{\lambda_{i
}}+a\vert^{2})^{\frac{n+2-2\gamma}{2}}} dx \\
\leq &
\frac{C\lambda_{i}}{(1+\lambda_{i}^{2}\vert a 
\vert^{2})^{\frac{n+1+2\gamma}{2}}}
\underset{\mathcal{B}_{1}}\int
\frac{1}
{(\frac{\lambda_{i}}{\lambda_{j}}+\lambda_{i}\lambda_{j}\vert\frac{x}{\lambda_{i
}}+a\vert^{2})^{\frac{n-2\gamma}{2}}}dx.
\end{split}
\end{equation*} 
Rescaling this gives 
\begin{equation*} 
\begin{split}
\vert \underset{\mathcal{B}_{1}}\int(\frac{ 1 }{ 
1+r^{2}})^{\frac{n+2\gamma}{2}} 
& \frac{ \lambda_{i}\lambda_{j}( \frac{x}{\lambda_{i}}+a) }{ 
(\frac{\lambda_{i}}{\lambda_{j}}+\lambda_{i}\lambda_{j}\vert 
\frac{x}{\lambda_{i}}+a\vert^{2}))^{\frac{n+2-2\gamma}{2}}} dx
\vert \\
\leq  &
\frac{C\lambda_{i}(\frac{\lambda_{i}}{\lambda_{j}})^{\frac{n+2\gamma}{2}}}{
(1+\lambda_{i}^{2}\vert a \vert^{2})^{\frac{n+1+2\gamma}{2}}}
\underset{[\vert \frac{x}{\lambda_{j}}+a\vert<\theta \vert a\vert]}\int
\frac{1}
{(1+\lambda_{j}^{2}\vert\frac{x}{\lambda_{j}}+a\vert^{2})^{\frac{n-2\gamma}{2}}}dx \\
\leq &
\frac{C\sqrt{\lambda_{i}\lambda_{j}}(\lambda_{j}\vert a 
\vert)^{2\gamma}}{(\frac{\lambda_{j}}{\lambda_{i}}+\lambda_{i}\lambda_{j}\vert 
a \vert^{2})^{\frac{n+1+2\gamma}{2}}}
=
O((\frac{\lambda_{j}}{\lambda_{i}})^{\gamma}\sqrt{\lambda_{i}\lambda_{j}}
\varepsilon_{i,j}^{\frac{n+1}{n-2\gamma}})dx.
\end{split}
\end{equation*} 
Moreover letting \;$\mathcal{B}_{2}=\mathcal{A}^{c}\setminus \mathcal{B}_{1}$\; we 
may estimate
\begin{equation*}
\begin{split}
\vert \underset{\mathcal{B}_{2}}\int & (\frac{ 1 }{ 
1+r^{2}})^{\frac{n+2\gamma}{2}} 
\frac{ \lambda_{i}\lambda_{j}( \frac{x}{\lambda_{i}}+a) }{ 
(\frac{\lambda_{i}}{\lambda_{j}}+\lambda_{i}\lambda_{j}\vert 
\frac{x}{\lambda_{i}}+a\vert^{2}))^{\frac{n+2-2\gamma}{2}}} dx
\vert \\
& \leq 
\frac{C\sqrt{\lambda_{i}\lambda_{j}}}{(\frac{\lambda_{i}}{\lambda_{j}}+\lambda_{
i}\lambda_{j}\vert a \vert^{2})^{\frac{n+1-2\gamma}{2}}}
\int_{\mathcal{A}^{c}}(\frac{1}{1+r^{2}})^{\frac{n+2\gamma}{2}}
=
O((\frac{\lambda_{j}}{\lambda_{i}})^{\gamma}\sqrt{\lambda_{i}\lambda_{j}}
\varepsilon_{i,j}^{\frac{n+1}{n-2\gamma}}).
\end{split}
\end{equation*} 
Collecting terms the claim follows.

\item 
By translation invariance and symmetry we may assume 
\;$
a_{i}=0, \;\;a_{j}=a \;\;\text{and}\;\;\frac{1}{\lambda_{i}}\leq \frac{1}{\lambda_{j}}
$.
Then
\begin{equation*}
\begin{split}
\nabla^{2}_{a}  \underset{\R^{n}}{\int} \delta_{i}^{\frac{n+2\gamma}{n-2\gamma}}\delta_{j}
= &
(2\gamma-n)\underset{\R^{n}}{\int}(\frac{\lambda_{i}}{1+\lambda_{i}^{2}r^{2}})^{\frac{n+2\gamma}{2}
}
(\frac{\lambda_{j}}{1+\lambda_{j}^{2}\vert x+a\vert^{2}})^{\frac{n-2\gamma}{2}} 
\\
& \quad \quad\ \quad\quad\;\;\, 
[
\frac{\lambda_{j}^{2}\,id}{1+\lambda_{j}^{2}\vert x+a\vert^{2}}
-
(n+2-2\gamma)\frac{\lambda_{j}^{4} (x+a).(x+a)}{(1+\lambda_{j}^{2}\vert 
x+a\vert^{2})^{2}}
]dx
\\
= &
(2\gamma-n)\lambda_{i}\lambda_{j}
\underset{\R^{n}}{\int}\frac
{
(\frac{\lambda_{i}}{\lambda_{j}}+\lambda_{i}\lambda_{j}\vert 
\frac{x}{\lambda_{i}}+a\vert^{2})\,id
-
(n+2-2\gamma)\lambda_{i}\lambda_{j} 
(\frac{x}{\lambda_{i}}+a).(\frac{x}{\lambda_{i}}+a)
}
{
(1+r^{2})^{\frac{n+2\gamma}{2}}(\frac{\lambda_{i}}{\lambda_{j}}+\lambda_{i}
\lambda_{j}\vert \frac{x}{\lambda_{i}}+a\vert^{2})^{\frac{n+4-2\gamma}{2}}
} dx.
\end{split}\end{equation*}
Since \;$
\eps_{i,j}^{\frac{2}{2\gamma-n}}
\sim 
\lambda_{i}\lambda_{j}\vert a\vert^{2}
\;\;\text{ or }\;\;
\eps_{i,j}^{\frac{2}{2\gamma-n}}
\sim
\frac{\lambda_{i}}{\lambda_{j}}
$\; we may expand on 
\begin{equation*}\begin{split}
\mathcal{A}
= &
[\vert  \frac{x}{\lambda_{i}}\vert \leq \epsilon\vert a\vert]
\cup
[\vert  \frac{x}{\lambda_{i}}\vert \leq \epsilon  \frac{1}{\lambda_{j}} ]
\end{split}\end{equation*}
for \;$\epsilon>0$\; sufficiently small 
\begin{equation*}
\begin{split}
& \hspace{-20pt}\frac
{
(\frac{\lambda_{i}}{\lambda_{j}}+\lambda_{i}\lambda_{j}\vert 
\frac{x}{\lambda_{i}}+a\vert^{2})\,id
-
(n+2-2\gamma)\lambda_{i}\lambda_{j} 
(\frac{x}{\lambda_{i}}+a).(\frac{x}{\lambda_{i}}+a)
}
{
(\frac{\lambda_{i}}{\lambda_{j}}+\lambda_{i}\lambda_{j}\vert 
\frac{x}{\lambda_{i}}+a\vert^{2})^{\frac{n+4-2\gamma}{2}}
} 
\\
= &
\frac
{
(\frac{\lambda_{i}}{\lambda_{j}}+\lambda_{i}\lambda_{j}\vert a\vert^{2})\,id
-
(n+2-2\gamma)\lambda_{i}\lambda_{j} a.a
+
O(\frac{\lambda_{j}}{\lambda_{i}}\vert x 
\vert^{2})
}
{
(\frac{\lambda_{i}}{\lambda_{j}}+\lambda_{i}\lambda_{j}\vert 
a\vert^{2})^{\frac{n+4-2\gamma}{2}}
}
\end{split}
\end{equation*} 
up to some integrable odd terms and thus obtain
\begin{equation*}
\begin{split}
\nabla_{a}^{2} \underset{\R^{n}}{\int} \delta_{i}^{\frac{n+2\gamma}{n-2\gamma}}\delta_{j} 
= &
b_{4}\lambda _{i}\lambda_{j}\frac
{
(\frac{\lambda_{i}}{\lambda_{j}}+\lambda_{i}\lambda_{j}\vert a\vert^{2})\,id
-
(n+2-2\gamma)\lambda_{i}\lambda_{j} a.a
}
{
(\frac{\lambda_{i}}{\lambda_{j}}+\lambda_{i}\lambda_{j}\vert 
a\vert^{2})^{\frac{n+4-2\gamma}{2}}
} \\
& +
O((\frac{\lambda_{j}}{\lambda_{i}})^{\gamma}\lambda_{i}\lambda_{j}\varepsilon_{i
,j}^{\frac{n+2}{n-2\gamma}}) 
+
(2\gamma-n)\lambda_{i}\lambda_{j}
\int_{\mathcal{A}^{c}} \;\cdots \;dx
,
\end{split}
\end{equation*} 
where \;$b_{4}=(2\gamma-n)c_{n,3}^{\gamma}$\; 
and
\begin{equation*}
\cdots
=
\frac
{
(\frac{\lambda_{i}}{\lambda_{j}}+\lambda_{i}\lambda_{j}\vert 
\frac{x}{\lambda_{i}}+a\vert^{2})\,id
-
(n+2-2\gamma)\lambda_{i}\lambda_{j} 
(\frac{x}{\lambda_{i}}+a).(\frac{x}{\lambda_{i}}+a)
}
{
(1+r^{2})^{\frac{n+2\gamma}{2}}(\frac{\lambda_{i}}{\lambda_{j}}+\lambda_{i}
\lambda_{j}\vert \frac{x}{\lambda_{i}}+a\vert^{2})^{\frac{n+4-2\gamma}{2}}
}
.
\end{equation*}
In case \;$\varepsilon_{i,j}^{\frac{2}{2\gamma-n}}\sim 
\frac{\lambda_{i}}{\lambda_{j}}$\; the integral above is of order
\begin{equation*}
\begin{split}
O((\frac{\lambda_{j}}{\lambda_{i}})^{\gamma}\varepsilon_{i,j}^{\frac{n+2}{
n-2\gamma}}).
\end{split}
\end{equation*} 
Thus we may assume \;$\frac{\lambda_{j}}{\lambda_{i}}\leq \frac{\lambda_{i}}{\lambda_{j}}\ll 
\lambda_{i}\lambda_{j}\vert a \vert^{2}$, i.e. for \;$0<\theta<\frac{1}{4}$
\begin{equation*}
\begin{split}
\mathcal{B}_{1}=\left[\vert \frac{x}{\lambda_{i}}+a\vert \leq \theta \vert 
a\vert\right]\subset \subset B^{c}_{\epsilon \lambda_{i}\vert a 
\vert}(0)=\mathcal{A}^{c}.
\end{split}
\end{equation*} 
On \;$\mathcal{B}_{1}$\; we then may expand  
\begin{equation*}
\begin{split}
(\frac{1}{1+r^{2}}  )^{\frac{n+2\gamma}{2}}
= &
(\frac{1}{1+\lambda_{i}^{2}\vert\frac{x}{\lambda_{i}}+a-a 
\vert^{2}})^{\frac{n+2\gamma}{2}} \\
= &
(\frac{1}{1+\lambda_{i}^{2}\vert a \vert^{2}})^{\frac{n+2\gamma}{2}}
+
(n+2\gamma)\frac{\lambda_{i}^{2}\langle 
\frac{x}{\lambda_{i}}+a,a\rangle + O(\lambda_{i}^{2}\vert 
\frac{x}{\lambda_{i}}+a\vert\vert a \vert)}{(1+\lambda_{i}^{2}\vert a 
\vert^{2})^{\frac{n+2+2\gamma}{2}}} .
\end{split}
\end{equation*}
Using radial symmetry and 
\;$
\int_{\R^{n}}\frac{1+r^{2}-(n+2-2\gamma)x_{k}^{2}}{(1+r^{2})^{\frac{n+4-2\gamma}
{2}}}=0
$\; 
for
\;$k=1,\ldots, n$, we obtain
\begin{equation*}
\begin{split}
\vert 
\int_{\mathcal{B}_{1}} & \frac
{
(\frac{\lambda_{i}}{\lambda_{j}}+\lambda_{i}\lambda_{j}\vert 
\frac{x}{\lambda_{i}}+a\vert^{2})\,id
-
(n+2-2\gamma)\lambda_{i}\lambda_{j} 
(\frac{x}{\lambda_{i}}+a).(\frac{x}{\lambda_{i}}+a)
}
{
(1+r^{2})^{\frac{n+2\gamma}{2}}(\frac{\lambda_{i}}{\lambda_{j}}+\lambda_{i}
\lambda_{j}\vert \frac{x}{\lambda_{i}}+a\vert^{2})^{\frac{n+4-2\gamma}{2}}
} dx
\vert 
\\
\leq &
\frac{C}{(1+\lambda_{i}^{2}\vert a \vert^{2})^{\frac{n+2\gamma}{2}}}
\int_{\mathcal{B}_{1}^{c}}  \frac
{
1
}
{
(\frac{\lambda_{i}}{\lambda_{j}}+\lambda_{i}\lambda_{j}\vert 
\frac{x}{\lambda_{i}}+a\vert^{2})^{\frac{n+2-2\gamma}{2}}
}dx
\\
& +
\frac{C\lambda_{i}^{2}\vert a \vert}{(1+\lambda_{i}^{2}\vert a 
\vert^{2})^{\frac{n+2+2\gamma}{2}}}
\int_{\mathcal{B}_{1}}  \frac
{
\vert \frac{x}{\lambda_{i}}+a\vert
}
{
(\frac{\lambda_{i}}{\lambda_{j}}+\lambda_{i}\lambda_{j}\vert 
\frac{x}{\lambda_{i}}+a\vert^{2})^{\frac{n+2-2\gamma}{2}}
} dx
\\
\leq &
\frac{C(\frac{\lambda_{j}}{\lambda_{i}})^{\frac{-n+2-2\gamma}{2}}}{(1+\lambda_{i
}^{2}\vert a \vert^{2})^{\frac{n+2\gamma}{2}}}
\int_{[\vert x \vert >\theta \lambda_{j}\vert a \vert]}  \frac
{
1
}
{
(1+r^{2})^{\frac{n+2-2\gamma}{2}}
}
\\
& +
\frac{C(\frac{\lambda_{j}}{\lambda_{i}})^{-\frac{n+2+2\gamma}{2}}\vert a \vert}{(1+\lambda_{i
}^{2}\vert a \vert^{2})^{\frac{n+2+2\gamma}{2}}}
\int_{[\vert x \vert \leq \theta\lambda_{j}\vert a \vert]}  \frac
{
r
}
{
(1+r^{2})^{\frac{n+2-2\gamma}{2}}
}
=
O((\frac{\lambda_{j}}{\lambda_{i}})^{\gamma}\varepsilon_{i,j}^{\frac{n+2}{
n-2\gamma
}}).
\end{split}
\end{equation*} 
Moreover letting \;$\mathcal{B}_{2}=\mathcal{A}^{c}\setminus \mathcal{B}_{1}$\; we 
may estimate
\begin{equation*}
\begin{split}
\vert \int_{\mathcal{B}_{2}} &\frac
{
(\frac{\lambda_{i}}{\lambda_{j}}+\lambda_{i}\lambda_{j}\vert 
\frac{x}{\lambda_{i}}+a\vert^{2})\,id
-
(n+2-2\gamma)\lambda_{i}\lambda_{j} 
(\frac{x}{\lambda_{i}}+a).(\frac{x}{\lambda_{i}}+a)
}
{
(1+r^{2})^{\frac{n+2\gamma}{2}}(\frac{\lambda_{i}}{\lambda_{j}}+\lambda_{i}
\lambda_{j}\vert \frac{x}{\lambda_{i}}+a\vert^{2})^{\frac{n+4-2\gamma}{2}}
} dx
\vert
\\
\leq &
\frac{C}{(\frac{\lambda_{i}}{\lambda_{j}+\lambda_{i}\lambda_{j}\vert a \vert^{2}})^{\frac{n+2-2\gamma}{2}}}
\int_{\mathcal{A}^{c}}
\frac{1}{(1+r^{2})^{\frac{n+2\gamma}{2}}}
\leq 
\frac{C(\frac{\lambda_{j}}{\lambda_{i}}+\lambda_{i}\vert a 
\vert)^{-2\gamma}}{(\frac{\lambda_{i}}{\lambda_{j}}+\lambda_{i}\lambda_{j}\vert a \vert^{2})^{\frac{n+2-2\gamma}{2}}}
= 
O((\frac{\lambda_{j}}{\lambda_{i}})^{\gamma}\varepsilon_{i,j}^{\frac{n+2}{
n-2\gamma}}).
\end{split}
\end{equation*} 
Collecting terms the claim follows.
\end{enumerate}
Thereby \;(i)-(iv)\; are proven and so is the lemma.
\end{pfn}

\bigskip

\noindent
\begin{pfn}{ of Lemma \ref{lem_higher_interaction_estimates}}
\begin{enumerate}[label=(\roman*)]
 \item 
Let \;$\alpha'=\frac{n-2\gamma}{2}\alpha, \; \beta'=\frac{n-2\gamma}{2}\beta$, so \;$\alpha'+\beta'=n$.
We distinguish the cases
\begin{enumerate}
 \item [($\alpha$)] \quad 
"$\eps_{i,j}^{\frac{2}{2\gamma-n}}\sim \frac{\lambda_{i} }{ \lambda_{j} }
\; \vee \;
\eps_{i,j}^{\frac{2}{2\gamma-n}} \sim \lambda_{i}\lambda_{j}\gamma_{n,\gamma}G_{h}^{\frac{2}{2\gamma-n}}( a _{i}, a _{j})$"
\quad 
We estimate for \;$c>0$\; small
\begin{equation*}\begin{split}
\int_{M} v_{i}^{\alpha} & v_{j}^{\beta}dV_{h}
\leq 
C \underset{B_{c}(0)}{\int}(\frac{\lambda_{i} }{1+\lambda_{i}^2r^{2}})^{\alpha'}(\frac{ \lambda_{j} }{1 +\lambda_{j}^{2}\gamma_{n,\gamma}G_{ a _{j}}^{\frac{2}{2\gamma-n}}(\exp_{ a _{i}}x)})^{\beta'}
+
C  \frac{1}{\lambda_{i}^{\alpha'}}\underset{B_{c}(0)}{\int}(\frac{ \lambda_{j} }{1+\lambda_{j}^2r^{2}})^{\beta'}
 \\
\leq &
C \underset{B_{ c\lambda_{i} }(0)}{\int}(\frac{1}{1+r^{2}})^{\alpha'}
(\frac{1}{\frac{\lambda_{i} }{ \lambda_{j} }
+
\lambda_{i}\lambda_{j}\gamma_{n,\gamma}G_{ a _{j}}^{\frac{2}{2\gamma-n}}(\exp_{ a _{i}}\frac{x}{\lambda_{i}})})^{\beta'}
+
C  \frac{1}{\lambda_{i}^{\alpha'}} \frac{1}{\lambda_{j}^{n-\beta'}} 
\underset{B_{ c\lambda_{j} }(0)}{\int}(\frac{1}{1+r^{2}})^{\beta'}
\end{split}\end{equation*}
up to some \;$O(  \frac{1}{\lambda_{i}^{\alpha'}} \frac{1}{\lambda_{j}^{\beta'}} )$. Thus by 
\;$
\int_{B_{ c\lambda_{j} }(0)}(\frac{1}{1+r^{2}})^{\beta'}
\leq 
C \frac{1}{\lambda_{j} ^{2\beta'-n}}
$\; 
we get
\begin{equation*}\begin{split}
\int_{M}   v_{i}^{\alpha}  v_{j}^{\beta}  dV_{h}
\leq &
C\underset{B_{ c\lambda_{i} }(0)}{\int}(\frac{1}{1+r^{2}})^{\alpha'}(\frac{1}{\frac{\lambda_{i} }{ \lambda_{j} }+\lambda_{i}\lambda_{j}\gamma_{n,\gamma}G_{ a _{j}}^{\frac{2}{2\gamma-n}}(\exp_{ a _{i}} \frac{x}{\lambda_{i}})})^{\beta'} 
\end{split}\end{equation*}
up to some \;$O(  \frac{1}{\lambda_{i}^{\alpha'}} \frac{1}{\lambda_{j}^{\beta'}} )$, whence the claim follows in cases
\begin{equation*}\begin{split}
\frac{ \lambda_{j} }{ \lambda_{i} }+\frac{\lambda_{i} }{ \lambda_{j} }
+
\lambda_{i}\lambda_{j}\gamma_{n,\gamma}G_{h}^{\frac{2}{2\gamma-n}}( a _{i}, a _{j})
\sim \frac{\lambda_{i}}{\lambda_{j}}
\; \text{ or }\;
d_{g_{a_{j}}}( a _{j}, a _{i})>3c.
\end{split}\end{equation*}
Else we may assume \;$d_{g_{a_{j}}}( a _{j}, a _{i})<3c$\; and 
\begin{equation*}\begin{split}
\frac{ \lambda_{j} }{ \lambda_{i} }+\frac{\lambda_{i} }{ \lambda_{j} }+\lambda_{i}\lambda_{j}\gamma_{n,\gamma}G_{h}^{\frac{2}{2\gamma-n}}( a _{i}, a _{j}) 
& \sim 
\lambda_{i}\lambda_{j}d^{2}_{g_{a_{j}}}( a _{j}, a _{i}).
\end{split}\end{equation*}
We then get with \;$\mathcal{B}=[\frac{1}{2}d_{g_{a_{j}}}( a _{j}, a _{i})\leq \vert  \frac{x}{\lambda_{i}}\vert \leq 2 d_{g_{a_{j}}}( a _{j}, a _{i})]$
\begin{equation*}\begin{split}
\int_{M} v_{i}^{\alpha}  v_{j}^{\beta} dV_{h}
\leq &
C \underset{\mathcal{B}}{\int}
(\frac{1}{1+r^{2}})^{\alpha'}(\frac{1}{\frac{\lambda_{i} }{ \lambda_{j} }+\lambda_{i}\lambda_{j}d^{2}_{g_{a_{j}}}( a _{j}, \exp_{ a _{i}}( \frac{x}{\lambda_{i}}))})^{\beta'} +
O(\eps_{i,j}^{\beta})\\
\leq &
C
(\frac{1}{1+\vert\lambda_{i}d_{h}( a _{i}, a _{j})\vert^{2}})^{\alpha'}
\underset{ [\vert  \frac{x}{\lambda_{i}}\vert \leq 4 d_{h}( a _{i}, a _{j})]}{\int}
(\frac{1}{\frac{\lambda_{i} }{ \lambda_{j} }+\frac{ \lambda_{j} }{ \lambda_{i} }r^{2}})^{\beta'}  +
O(\eps_{i,j}^{\beta})\\
\leq &
C
\frac{(\frac{ \lambda_{j} }{ \lambda_{i} })^{\beta'-n}}{(1+\vert\lambda_{i}d_{h}( a _{i}, a _{j})\vert^{2})^{\alpha'}}
 \underset{ [r \leq 4 \lambda_{j}d_{h}( a _{i}, a _{j})]}{\int}
(\frac{1}{1+r ^{2}})^{\beta'}
+
O(\eps_{i,j}^{\beta}).
\end{split}\end{equation*}
Note, that in case 
\;$\lambda_{j}d_{g_{a_{j}}}( a _{j}, a _{i})$\; remains bounded, we are done. Else
\begin{equation*}\begin{split}
\int_{M} v_{i}^{\alpha}v_{j}^{\beta}dV_{h}
\leq &
C
\frac{(\frac{ \lambda_{j} }{ \lambda_{i} })^{\beta'-n}(\lambda_{j}d_{h}( a _{j}, a _{i}))^{n-2\beta'}}{(1+\vert\lambda_{i}d_{h}( a _{j}, a _{i})\vert^{2})^{\alpha'}}
+
O(\eps_{i,j}^{\beta}) \\
\leq &
C
(\frac{1}{1+\vert\lambda_{i}d_{h}( a _{i}, a _{j})\vert^{2}})^{\alpha'-\frac{n}{2}+\beta'}
(\frac{\lambda_{i} }{ \lambda_{j} })^{\beta'}+
O(\eps_{i,j}^{\beta}),
\end{split}\end{equation*}
whence due to \;$\alpha'>\frac{n}{2}$\; the claim follows.

 \item [$(\beta)$]\quad "$\eps_{i,j}^{\frac{2}{2\gamma-n}}\sim \frac{ \lambda_{j} }{ \lambda_{i} }$"
\quad 
We estimate for \;$c>0$\; small
\begin{equation*}\begin{split}
\int_{M} v_{i}^{\alpha}v_{j}^{\beta} dV_{h}
\leq &
C
\int_{B_{ c\lambda_{j} }(0)}
(\frac{1}{\frac{ \lambda_{j} }{ \lambda_{i} }+\lambda_{i}\lambda_{j}\gamma_{n,\gamma}G_{ a _{i}}^{\frac{2}{2\gamma-n}}( \exp_{ a _{j}}( \frac{x}{\lambda_{j}} ))})^{\alpha'} 
(\frac{1}{1+r^{2}})^{\beta'} \\
& +
C \frac{1}{\lambda_{j}^{\beta'}} \int_{B_{c}(0)}(\frac{\lambda_{i} }{1+\lambda_{i}^2r^{2}})^{\alpha'}
+
O(  \frac{1}{\lambda_{i}^{\alpha'}} \frac{1}{\lambda_{j}^{\beta'}} ),
\end{split}\end{equation*}
which by
\;$
\int_{B_{c}(0)}(\frac{\lambda_{i} }{1+\lambda_{i}^2r^{2}})^{\alpha'}
\leq
C\frac{1}{\lambda_{i}^{n-\alpha'}} =  C\frac{1}{\lambda_{i} ^{\beta'}}
$\;
gives

\begin{equation*}\begin{split}
\int  v_{i}^{\alpha}v_{j}^{\beta} 
\leq &
C
\int_{B_{ c\lambda_{j} }(0)}
(\frac{1}{\frac{ \lambda_{j} }{ \lambda_{i} }+\lambda_{i}\lambda_{j}\gamma_{n,\gamma}G_{a_{i}}^{\frac{2}{2\gamma-n}}(\exp_{ a _{j}}( \frac{x}{\lambda_{j}}))})^{\alpha'} 
(\frac{1}{1+r^{2}})^{\beta'} 
\end{split}\end{equation*}
up to some \;$O(\eps_{i,j}^{\beta})$. Since by assumption \;$d_{h}( a _{i}, a _{j})\leq  \frac{2}{\lambda_{i}}$, there holds
\begin{equation*}\begin{split}
\gamma_{n,\gamma}G_{ a _{i}}^{\frac{2}{2\gamma-n}}(\exp_{ a _{j}}( \frac{x}{\lambda_{j}}))\sim d_{h}^{2}( a _{i}, \exp_{ a _{j}}( \frac{x}{\lambda_{j}}))
\;\;\text{ on }\;\;
B_{ c\lambda_{j} }(0).
\end{split}\end{equation*}
Thus for \;$\gamma>3$\;
\begin{equation*}\begin{split}
\int_{M}  v_{i}^{\alpha}v_{j}^{\beta} dV_{h}
\leq &
C
\underset{[\gamma\frac{ \lambda_{j} }{ \lambda_{i} }\leq \vert x \vert\leq  c\lambda_{j} ]}{\int}
(\frac{1}{\frac{ \lambda_{j} }{ \lambda_{i} }+\lambda_{i}\lambda_{j}d_{h}^{2}( a _{i}, \exp_{ a _{j}}( \frac{x}{\lambda_{j}}))})^{\alpha'} 
(\frac{1}{1+r^{2}})^{\beta'} \\
& +C
(\frac{\lambda_{i} }{ \lambda_{j} })^{\alpha'}
\int_{[\vert x \vert<\gamma\frac{ \lambda_{j} }{ \lambda_{i} }]}
(\frac{1}{1+r^{2}})^{\beta'} 
+
o(\eps_{i,j}^{\beta}) \\
\leq &
C
\int_{ [\gamma\frac{ \lambda_{j} }{ \lambda_{i} }\leq \vert x \vert\leq  c\lambda_{j} ]}
(\frac{1}{\frac{ \lambda_{j} }{ \lambda_{i} }+\frac{\lambda_{i} }{ \lambda_{j} }r^{2}})^{\alpha'} 
(\frac{1}{1+r^{2}})^{\beta'} 
+
C
(\frac{\lambda_{i} }{ \lambda_{j} })^{\alpha'}
(\frac{ \lambda_{j} }{ \lambda_{i} })^{n-2\beta} 
+
o(\eps_{i,j}^{\beta}),
\end{split}\end{equation*}
since for \;$\vert x \vert \geq \gamma\frac{ \lambda_{j} }{ \lambda_{i} }$\; we may due to \;$d_{h}( a _{i}, a _{j})\leq \frac{2}{\lambda_{i}} $\; assume 
\;$
d_{h}( a _{i}, \exp_{ a _{j}}( \frac{x}{\lambda_{j}}))\geq  \frac{\vert x \vert}{\lambda_{j}}
$.  
Therefore
\begin{equation*}\begin{split}
\int_{M} v_{i}^{\alpha}v_{j}^{\beta}dV_{h}
\leq &
C
(\frac{ \lambda_{j} }{ \lambda_{i} })^{\alpha'}\int_{ [\vert x \vert \geq \gamma\frac{ \lambda_{j} }{ \lambda_{i} }]}
r^{-2n}
+
O(\eps_{i,j}^{\beta}) 
= 
O(\eps_{i,j}^{\beta}).
\end{split}\end{equation*}

\end{enumerate}
\item 
By symmetry we may assume \;$\frac{1}{\lambda_{i}}\leq \frac{1}{\lambda_{j}}$\; and thus
\begin{equation*}\begin{split}
\eps_{i,j}^{\frac{2}{2\gamma-n}}\sim \frac{\lambda_{i} }{ \lambda_{j} }
\; \vee \;
\eps_{i,j}^{\frac{2}{2\gamma-n}} \sim \lambda_{i}\lambda_{j}\gamma_{n,\gamma}G_{h}^{\frac{2}{2\gamma-n}}( a _{i}, a _{j})
\end{split}\end{equation*}
In case \;$d_{h_{a_{j}}}(a_{j},a_{i})>3\epsilon$\; for some \;$\epsilon>0$\;  we estimate
\begin{equation*}\begin{split}
\int_{M} v_{i}^{\frac{n}{n-2\gamma}}v_{j}^{\frac{n}{n-2\gamma}} dV_{h}
\leq &
C \int_{B_{\epsilon}(0)}(\frac{\lambda_{i} }{1+\lambda_{i}^2r^{2}})^{\frac{n}{2}}
(\frac{ \lambda_{j} }{1 +\lambda_{j}^{2}\gamma_{n,\gamma}G_{ a _{j}}^{\frac{2}{2\gamma-n}}(\exp_{ a _{i}}x)})^{\frac{n}{2}}\\
& +
\frac{C_{\epsilon}}{\lambda_{i}^{\frac{n}{2}}}
\int_{B_{\epsilon}(0)}(\frac{ \lambda_{j} }{1+\lambda_{j}^2r^{2}})^{\frac{n}{2}}
+
C_{\epsilon} \frac{1}{\lambda_{i}^{\frac{n}{2}}\lambda_{j}^{\frac{n}{2}}}   \\
\leq &
C \int_{B_{ \epsilon\lambda_{i} }(0)}(\frac{1}{1+r^{2}})^{\frac{n}{2}}
(\frac{1}{\frac{\lambda_{i} }{ \lambda_{j} }
+
\lambda_{i}\lambda_{j}\gamma_{n,\gamma}G_{ a _{j}}^{\frac{2}{2\gamma-n}}(\exp_{ a _{i}}( \frac{x}{\lambda_{i}}))})^{\frac{n}{2}}+
C_{\epsilon}\ln \lambda_{j}\eps_{i,j}^{\frac{n}{n-2\gamma}} \\
\leq &
C_{\epsilon}\ln(\lambda_{i}\lambda_{j})\epsilon_{i,j}^{\frac{n}{n-2\gamma}}
\leq C_{\epsilon}\ln \epsilon_{i,j}\epsilon_{i,j}^{\frac{n}{n-2\gamma}}.
\end{split}\end{equation*}
Thus we assume, that \;$d_{h_{a_{j}}}(a_{j},a_{i})$\; is arbitrarily small. Then for
\;$d_{h_{a_{j}}}(a_{j},a_{i})\ll c \ll 1 $\;
we estimate passing to normal coordinates around \;$a_{i}$\;
\begin{equation*}
\begin{split}
\int_{M} v_{i}^{\frac{n}{n-2\gamma}}v_{j}^{\frac{n}{n-2\gamma}} dV_{h}
\leq &
C\int_{B_{c\lambda_{i}}(0)}
(\frac{1}{1+r^{2}})^{\frac{n}{2}}
(
\frac{1}
{
\frac{\lambda_{i}}{\lambda_{j}}
+
\lambda_{i}\lambda_{j}\gamma_{n,\gamma}G_{a_{j}}^{\frac{2}{2\gamma-n}}(\exp_{a_{i}}(\frac{x}{\lambda_{i}}))
}
)^{\frac{n}{2}}  
=
C\int_{B_{c\lambda_{i}}(0)} I
\end{split}
\end{equation*}
up to some terms of order \;$\frac{1}{\lambda_{i}^{\frac{n}{2}}\lambda_{j}^{\frac{n}{2}}}=O(\epsilon_{i,j}^{\frac{n}{n-2\gamma}})$. Decompose \;$B_{c\lambda_{i}}(0)$\; into
\begin{enumerate}
 \item[$(\alpha)$] \quad 
$
\mathcal{A}
=
[
\vert \frac{x}{\lambda_{i}}\vert
\leq 
\epsilon \sqrt{G_{a_{j}}^{\frac{2}{2\gamma-n}}(a_{i})+\frac{1}{\lambda_{j}^{2}}}
]
$
 \item[$(\beta)$] \quad 
$
\mathcal{B}
=
[
\epsilon \sqrt{G_{a_{j}}^{\frac{2}{2\gamma-n}}(a_{i})+\frac{1}{\lambda_{j}^{2}}}
<
\vert \frac{x}{\lambda_{i}}\vert
<
E \sqrt{G_{a_{j}}^{\frac{2}{2\gamma-n}}(a_{i})+\frac{1}{\lambda_{j}^{2}}}
]
$
 \item[$(\gamma)$] \quad 
$
\mathcal{C}
=
[
E \sqrt{G_{a_{j}}^{\frac{2}{2\gamma-n}}(a_{i})+\frac{1}{\lambda_{j}^{2}}}
\leq
\vert \frac{x}{\lambda_{i}}\vert
<
c
]
$
\end{enumerate}
for some fixed \;$0<\epsilon \ll 1 \ll E <\infty$. We then find
\begin{equation*}
\begin{split}
\int_{\mathcal{A}}I
\leq 
C\varepsilon_{i,j}^{\frac{n}{n-2\gamma}}
\int_{\mathcal{A}}(\frac{1}{1+r^{2}})^{\frac{n}{2}} 
=
O(\ln \epsilon_{i,j}\epsilon_{i,j}^{\frac{n}{n-2\gamma}})
\end{split}
\end{equation*}
and 
\begin{equation*}
\begin{split}
\int_{\mathcal{B}} I
\leq &
\frac{C}
{
(\frac{\lambda_{i}^{2}}{\lambda_{j}^{2}}
+
\lambda_{i}^{2}G_{a_{j}}^{\frac{2}{2\gamma-n}}(a_{i})
)^{\frac{n}{2}}}
\int_{\mathcal{B}}
(
\frac{1}
{
\frac{\lambda_{i}}{\lambda_{j}}
+
\lambda_{i}\lambda_{j}\gamma_{n,\gamma}G_{a_{j}}^{\frac{2}{2\gamma-n}}(\exp_{a_{i}}(\frac{x}{\lambda_{i}}))
}
)^{\frac{n}{2}}  \\
\leq &
C\varepsilon_{i,j}^{\frac{n}{n-2\gamma}}(\frac{\lambda_{j}}{\lambda_{i}})^{\frac{n}{2}}
\underset
{
[
\vert \frac{x}{\lambda_{i}}\vert <
\tilde E \sqrt{G_{a_{j}}^{\frac{2}{2\gamma-n}}(a_{i})+\frac{1}{\lambda_{j}^{2}}}
]
}
{\int}
(
\frac{1}
{
\frac{\lambda_{i}}{\lambda_{j}}
+
\frac{\lambda_{j}}{\lambda_{i}}\vert x \vert^{2}
}
)^{\frac{n}{2}} 
\\
= &
C\varepsilon_{i,j}^{\frac{n}{n-2\gamma}}
\underset
{
[
\vert x\vert < \lambda_{j}
\tilde E \sqrt{G_{a_{j}}^{\frac{2}{2\gamma-n}}(a_{i})+\frac{1}{\lambda_{j}^{2}}}
]
}
{\int}
(
\frac{1}
{
1
+
\vert x \vert^{2}
}
)^{\frac{n}{2}} 
=
O(\ln \epsilon_{i,j}\epsilon_{i,j}^{\frac{n}{n-2\gamma}})
\end{split}
\end{equation*}
and
\begin{equation*}
\begin{split}
\int_{\mathcal{C}} I
\leq &
C\int_{\mathcal{C}}
(\frac{1}{1+r^{2}})^{\frac{n}{2}}
(
\frac{1}
{
\frac{\lambda_{i}}{\lambda_{j}}
+
\frac{\lambda_{j}}{\lambda_{i}}\vert x \vert^{2}
}
)^{\frac{n}{2}}  
\leq
C
(\frac{\lambda_{i}}{\lambda_{j}})^{\frac{n}{2}}
\int_{\mathcal{C}}
(\frac{1}{1+r^{2}})^{\frac{n}{2}}r^{-n}
\\
\\
\leq &
C(\frac{\lambda_{i}}{\lambda_{j}})^{\frac{n}{2}}
\underset
{
E\lambda^{i} \sqrt{G_{a_{j}}^{\frac{2}{2\gamma-n}}(a_{i})+\frac{1}{\lambda_{j}^{2}}}
}
{\int^{\infty}}
r^{-1-n}
=
O(\epsilon_{i,j}^{\frac{n}{n-2\gamma}}).
\end{split}
\end{equation*}
\end{enumerate}
Collecting terms the assertion follows.
\end{pfn}

\noindent

\smallskip

We turn now to the proofs of the estimates 
\eqref{epsilon_i_j_symmetry}, \eqref{varphi_in_a_j_coordinates}, \eqref{I2_interaction_integral_estimate} and
\eqref{epsilon_varepsilon_relation}.

\smallskip 

\noindent
\begin{pfn}{ of the estimate \eqref{epsilon_i_j_symmetry}}. 
First we notice, that due to Corollary 4.6 in \cite{martndia3} and 
\eqref{bubble_boundary} we have
\begin{equation*}
\varphi_{k}=(1+o_{\varepsilon_{k}}(1))(\frac{\lambda_{k}}{1+\lambda_{k}^{2}G_{a_
{ k } } ^{ \frac { 2 } { 2\gamma-n } } } )^{\frac{n-2\gamma}{2}}
\end{equation*}
for any choice \;$\varepsilon_{k}\sim \lambda_{k}^{-\frac{1}{l}}$, whence 
\begin{equation}\label{epsilon_i,j_expansion}
\begin{split}
\epsilon_{i,j}
= &
\int_{B_{\varepsilon_{i}}(a_{i})}(u_{a_{i}}\varphi_{a_{i}, \lambda_{i}})^{\frac{
n+2\gamma }{
n-2\gamma } }u_{a_{j}}\varphi_{a_{j}}dV_{h} \\
= &
(1+o_{\frac{1}{\lambda_{j}}}(1))
\underset{B_{\varepsilon_{i}\lambda_{i}}(0)}{\int}
(\frac{1}{1+r^{2}})^{\frac{n+2\gamma}{2}}
\frac
{\frac{u_{a_{j}}}{u_{a_{i}}}(\exp_{g_{a_{i}}}\frac{x}{\lambda_{i}})}
{
(\frac{\lambda_{i}}{\lambda_{j}}+\lambda_{i}\lambda_{j}G_{a_{j}}^{\frac{2}{
2\gamma-n} }(\exp_{g_{a_{i}}}\frac{x}{\lambda_{i}}) )^{\frac{n-2\gamma}{2}}
}
\end{split}
\end{equation}
up to some  \;$o_{\frac{1}{\lambda_{i}}}((\lambda_{i}\lambda_{j})^{\frac{2\gamma-n}{2}})$. In case \;$d_{g_{a_{j}}}(a_{i},a_{j})>\sqrt{\varepsilon_{i}}$\; we 
thus get
\begin{equation*}
\epsilon_{i,j}
=
(1+o_{\max(\frac{1}{\lambda_{i}}, \frac{1}{\lambda_{j}})}(1))
(\frac{1}{\lambda_{i
}\lambda_{j}G_{h}^{\frac{2}{2\gamma-n}}(a_{i},a_{j}) })^{\frac{n-2\gamma}{2}}
\end{equation*}
by conformal covariance of the Green's function, i.e.
\begin{equation}\begin{split}\label{covariance_greens_function}
G_{a_{j}}(a_{i})=G_{h_{a_{j}}}(a_{i},a_{j})=u_{a_{j}}^{-1}(a_{i})u_{a_{j}}^{-1}(a_{j})G_{h}(a_{i}
,a_{j})\;\;\text{ recalling }\;\;
h_{a_{j}}=u_{a_{j}}^{\frac{4}{n-2\gamma}}h.
\end{split}\end{equation}
Therefore \;\eqref{epsilon_i_j_symmetry}\; follows by symmetry in this case. 
Otherwise we may assume  \;$d_{g_{a_{j}}}(a_{i},a_{j})\leq 
\varepsilon_{i}^{\frac{1}{4}}$\;
and rewriting the Green's function in \eqref{epsilon_i,j_expansion} in \;$h_{a_{i}}$-normal coordinates via \eqref{covariance_greens_function} 
we obtain
\begin{equation*}
\begin{split}
\epsilon_{i,j}
= &
(1+o_{\max(\frac{1}{\lambda_{i}}, \frac{1}{\lambda_{j}})}(1))
\int_{B_{\varepsilon_{i}\lambda_{i}}(0)}
(\frac{1}{1+r^{2}})^{\frac{n+2\gamma}{2}}
(\frac
{1}
{
\frac{\lambda_{i}}{\lambda_{j}}+\lambda_{i}\lambda_{j}\vert 
\exp^{-}_{g_{a_{i}}}a_{j} -\frac{x}{\lambda_{i} 
}\vert^{2}
}
)^{\frac{n-2\gamma}{2}} \\
= &
(1+o_{\max(\frac{1}{\lambda_{i}}, \frac{1}{\lambda_{j}})}(1))
\int_{\R^{n}}
(\frac{\lambda_{i}}{1+\lambda_{i}^{2}r^{2}})^{\frac{n+2\gamma}{2}}
(\frac
{\lambda_{j}}
{
1+\lambda_{j}^{2}\vert 
\exp^{-}_{a_{i}}a_{j} -x\vert^{2}
}
)^{\frac{n-2\gamma}{2}}
\end{split}
\end{equation*}
up to some 
\;$o_{\frac{1}{\lambda_{i}}}((\lambda_{i}\lambda_{j})^{\frac{2\gamma-n}{2}})$.
Using \;$\exp^{-}_{g_{a_{i}}}a_{j}=(1+o_{\frac{1}{\lambda_{i}}}(1))\exp^{-}_{g_{a_{j}}}
a_ {i} $\; and 
\begin{equation*}
\int_{\R^{n}} \delta_{i}^{\frac{n+2\gamma}{n-2\gamma}}\delta_{j}
=\frac{1}{c_{n, \gamma}}
\int_{\R^{n}} (-\Delta)^\gamma\delta_{i}\delta_{j}
=\frac{1}{c_{n, \gamma}}
\int_{\R^{n}} \delta_{i}(-\Delta)^\gamma\delta_{j}
=
\int_{\R^{n}} \delta_{i}\delta_{j}^{\frac{n+2\gamma}{n-2\gamma}}
\end{equation*}
we thus find up to some 
\;$o_{\frac{1}{\lambda_{i}}}
((\lambda_{i}\lambda_{j})^{\frac{2\gamma-n}{2}})$
\begin{equation*}
\begin{split}
\epsilon_{i,j}
= &
(1+o_{\max(\frac{1}{\lambda_{i}}, \frac{1}{\lambda_{j}})}(1))
\int_{\R^{n}}
(\frac{\lambda_{j}}{1+\lambda_{j}^{2}r^{2}})^{\frac{n+2\gamma}{2}}
(\frac
{\lambda_{i}}
{
1+\lambda_{i}^{2}\vert 
\exp^{-}_{g_{a_{j}}}a_{i} -x\vert^{2}
}
)^{\frac{n-2\gamma}{2}}.
\end{split}
\end{equation*}
This shows \;$\epsilon_{i,j}=O(\varepsilon_{i,j})$, cf. Lemma 
\ref{lem_standard interaction on Rn} (i), and calculating back we find  
\begin{equation*}
\epsilon_{i,j}
=
(1+o_{\max(\frac{1}{\lambda_{i}}, \frac{1}{\lambda_{j}})} 
(1))\epsilon_{j,i}
+
o_{\max(\frac{1}{\lambda_{i}}, \frac{1}{\lambda_{j}})} 
((\lambda_{i}\lambda_{j})^{\frac{2\gamma-n}{2}})
=
(1+o_{\max(\frac{1}{\lambda_{i}}, \frac{1}{\lambda_{j}})} 
(1))\epsilon_{j,i}.
\end{equation*}
Thereby \;\eqref{epsilon_i_j_symmetry} follows.
\end{pfn}

\noindent
\begin{pfn}{ of the estimate \eqref{varphi_in_a_j_coordinates}}. 
Note, that due to Corollary 4.3 in \cite{martndia3} we have
\begin{equation*}
\begin{split}
K_{g_{a_{i}}}(y,x,\xi)
\leq &
Cy^{2\gamma}
(
(\frac{1}{y^{2}+d^{2}_{g_{a_{i}}}(x,\xi)})^{\frac{n+2\gamma}{2}}
+
1
)
\end{split}
\end{equation*} 
for \;$(y,x)\in B_{\rho_{0}}^{g_{a_i}, +}(a_{i})$, whence
\begin{equation*}
\begin{split}
\vert \overline{\frac{u_{a_{j}}}{u_{a_{i}}}\varphi_{j}}^{a_{i}} \vert(y,x)
\leq & 
C
\underset{M}{\int} 
(
\frac{y^{2\gamma}}{(y^{2}+d^{2}_{g_{a_{i}}}(x,\xi))^{\frac{n+2\gamma}{2}}}
+
y^{2\gamma}
) 
(\frac{\lambda_{j}}{1+\lambda_{j}^{2}d_{g_{a_{j}}}^{2}(a_{j}, \xi)})^{\frac{
n-2\gamma } {2}
}dV_{g_{a_{i}}}
\end{split}
\end{equation*} 
and thus 
\begin{equation*}\begin{split}
\vert \overline{\frac{u_{a_{j}}}{u_{a_{i}}}\varphi_{j}}^{a_{i}} & \vert(y,x)
\leq 
\frac{C}{y^{\frac{n-2\gamma}{2}}}
\int_{M} 
(\frac{y^{-1}}{1+y^{-2}d^{2}_{g_{a_{i}}}(x,\xi)})^{\frac{n+2\gamma}{2}}
(\frac{\lambda_{j}}{1+\lambda_{j}^{2}d_{g_{a_{j}}}^{2}(a_{j}, \xi)})^{\frac{
n-2\gamma}{2}}dV_{g_{a_{i}}}
+
\frac{Cy^{2\gamma}}{\lambda_{j}^{\frac{n-2\gamma}{2}}}
\\
\leq &
\frac{C}{y^{\frac{n-2\gamma}{2}}} 
\underset{[d_{g_{a_{i}}}(x,\xi)\leq \rho_{0}]\cap[d_{g_{a_{j}}}(a_{j}, \xi)\leq 
\rho_{0}]} {\int}
(\frac{y^{-1}}{1+y^{-2}d^{2}_{g_{a_{i}}}(x,\xi)})^{\frac{n+2\gamma}{2}}(\frac{\lambda_{j}}{1+\lambda_{j}^{2}d_{g_{a_{j}}}^{2}(a_{j}, \xi)})^{\frac{
n-2\gamma}{2}}dV_{g_{a_{i}}}
+
\frac{Cy^{2\gamma}}{\lambda_{j}^{\frac{n-2\gamma}{2}}}.
\end{split}
\end{equation*}
As all the distances involved are comparable to \;$d_{g}$, the first summand 
above vanishes for  \;$d_{g}(x,a_{j})\geq C\rho_{0}$\; and the claim follows.
Else we may assume \;$d_{g}(x,a_{j})\leq C\rho_{0}$\; and find passing to 
integration over \;$\R^{n}$\; by standard interaction estimates as given in Lemma 
\ref{lem_standard interaction on Rn} 
\begin{equation*}\begin{split}
\vert \overline{\frac{u_{a_{j}}}{u_{a_{i}}}\varphi_{j}}^{a_{i}} \vert(y,x)
\leq &
\frac{C}{y^{\frac{n-2\gamma}{2}}}
(\frac{1}
{
\lambda_{j} y+\frac{1}{\lambda_{j} y}+\frac{\lambda_{j}}{y}d^{2}_{g}(a_{j},x)
})^{\frac{n-2\gamma}{2}}
+
\frac{Cy^{2\gamma}}{\lambda_{j}^{\frac{n-2\gamma}{2}}} \\
\leq &
C(\frac{\lambda_{j}}{1+\lambda_{j}^{2}\vert 
y^{2}+d_{g}^{2}(a_{j},x)\vert})^{\frac{n-2\gamma}{2}}
+
\frac{C}{\lambda_{j}^{\frac{n-2\gamma}{2}}}.
\end{split}
\end{equation*} 
The proof is thereby complete.
\end{pfn}

\bigskip

\noindent
\begin{pfn}{ of the estimate \eqref{I2_interaction_integral_estimate}}. 
We let \;$B_{r}^{+}(a)=B_{r}^{g_{a}, +}(a)$\; and start showing 
\begin{equation}\label{I2_interaction_integral_estimate_proof}
\begin{split}
\int_{B^{ +}_{2\varepsilon_{i}}(a_{i})}y^{1-2\gamma}\widetilde{\varphi_{i}}^{a_{
i }} \widetilde{\varphi_{j}}^{a_{j}} 
=
o_{\frac{1}{\lambda_{i}}}(\varepsilon_{i,j}).
\end{split}
\end{equation}  
In case \;$d_{g}(a_{i},a_{j})>c\gg \varepsilon_{i}$\; we have
\begin{equation*}
\begin{split}
\int_{B^{+}_{2\varepsilon_{i}}(a_{i})}y^{1-2\gamma}\widetilde{\varphi_{i}}^{a_{
i }} \widetilde{\varphi_{j}}^{a_{j}}
\leq \frac{C}{\lambda_{i}^{\frac{n-2\gamma}{2}}}
\int_{B^{+}_{2\varepsilon_{i}}(a_{i})}y^{1-2\gamma}\widetilde{\varphi_{i}}^{a_{
i }}
=
O(\frac{\varepsilon_{i}^{2}}{(\lambda_{i}\lambda_{j})^{\frac{n-2\gamma}{2}}}),
\end{split}
\end{equation*}
so \;\eqref{I2_interaction_integral_estimate_proof}\; holds true in this case for 
any choice \;$\varepsilon_{i}\sim \lambda_{i}^{-\frac{1}{k}}$. Thus we may 
assume \;$d_{g}(a_{i},a_{j})\ll 1$\; for the rest of the proof and moreover, that
\begin{equation*}
\varepsilon_{i,j}^{\frac{2\gamma-n}{2}}\sim \frac{\lambda_{i}}{\lambda_{j}}\;\;
\text{ or }\,\;\;
\varepsilon_{i,j}^{\frac{2\gamma-n}{2}}\sim 
\lambda_{i}\lambda_{j}G_{h}^{\frac{2}{2\gamma-n}}(a_{i},a_{j}).
\end{equation*}
Passing to \;$g_{a_{i}}$- normal Fermi-coordinates and rescaling we then have
\begin{equation*}
\begin{split}
\int_{B^{ +}_{2\varepsilon_{i}}(a_{i})} 
y^{1-2\gamma}\widetilde{\varphi_{i}}^{a_{
i }} \widetilde{\varphi_{j}}^{a_{j}}  
\leq 
\frac{C}{\lambda_{i}^{2}}
\int_{B^{ +}_{2\varepsilon_{i}\lambda_{i}}(0)}
\frac{y^{1-2\gamma}}{(1+r^{2})^{\frac{n-2\gamma}{2}}}
(\frac{1}{\frac{\lambda_{i}}{\lambda_{j}}
+
\lambda_{i}\lambda_{j}\vert \frac{y}{\lambda_{i}},
a_{j}-a_{i}+\frac{x}{\lambda_{i}})\vert^{2}})^{\frac{n-2\gamma}{2}}.
\end{split}
\end{equation*}
In particular \;\eqref{I2_interaction_integral_estimate_proof} holds
in case \;$\varepsilon_{i,j}^{\frac{2\gamma-n}{2}}\sim 
\frac{\lambda_{i}}{\lambda_{j}}$, so we may assume 
\begin{equation*}
\varepsilon_{i,j}^{\frac{2\gamma-n}{2}}\sim 
\lambda_{i}\lambda_{j}G_{h}^{\frac{2}{2\gamma-n}}(a_{i},a_{j})
\sim 
\lambda_{i}\lambda_{j}\vert a_{i}-a_{j}\vert^{2}.
\end{equation*}
We subdivide the region of integration, 
i.e.  \;$B_{2\varepsilon_{i}\lambda_{i}}^{+}=B_{2\varepsilon_{i}\lambda_{i}}^{+}(0)$\; into
\begin{enumerate}[label=(\roman*)]
 \item \quad
$
\mathcal{B}_{1}
=
[\vert \frac{z}{\lambda_{i}}\vert \leq \epsilon \vert a_{i}-a_{j}\vert]
\cap B_{2\varepsilon_{i}\lambda_{i}}^{+}
$
 \item \quad 
$
\mathcal{B}_{2}
=
[\epsilon\vert a_{i}-a_{j}\vert <\vert \frac{z}{\lambda_{i}}\vert \leq E 
\vert a_{i}-a_{j}\vert]
\cap B_{2\varepsilon_{i}\lambda_{i}}^{+}
$
 \item \quad
 $
\mathcal{B}_{1}
=
[\vert \frac{z}{\lambda_{i}}\vert > E\vert a_{i}-a_{j}\vert]
\cap B_{2\varepsilon_{i}\lambda_{i}}^{+}
$
\end{enumerate}
for \;$z=(y,x)$\; and suitable \;$0<\epsilon,E^{-1}\ll 1$\; and obtain easily
\begin{equation*}
\begin{split}
\int_{B^{ +}_{2\varepsilon_{i}}(a_{i})} 
y^{1-2\gamma}\widetilde{\varphi_{i}}^{a_{
i }} \widetilde{\varphi_{j}}^{a_{j}}  
\leq 
\frac{C}{\lambda_{i}^{2}}
\int_{\mathcal{B}_{2}}
\frac{y^{1-2\gamma}}{(1+r^{2})^{\frac{n-2\gamma}{2}}}
(\frac{1}{\frac{\lambda_{i}}{\lambda_{j}}
+
\lambda_{i}\lambda_{j}\vert (\frac{y}{\lambda_{i}},
a_{j}-a_{i}+\frac{x}{\lambda_{i}})\vert^{2}})^{\frac{n-2\gamma}{2}}
\end{split}
\end{equation*}
up to some \;$o_{\frac{1}{\lambda_{i}}}(\varepsilon_{i,j})$\; for any choice \;
$\varepsilon_{i}\sim \lambda_{i}^{-\frac{1}{k}}$. Note, that on 
\;$\mathcal{B}_{2}$\; we have
\begin{equation*}
\epsilon \vert a_{i}-a_{j}\vert <\vert \frac{z}{\lambda_{i}}\vert \leq 
2\varepsilon_{i},
\end{equation*}
so  \;$\vert a_{i}-a_{j}\vert =O(\varepsilon_{i})$\; for 
\;$\mathcal{B}_{2}\neq \emptyset$.
We then find up to some \;$o_{\frac{1}{\lambda_{i}}}(\varepsilon_{i,j})$
\begin{equation*}
\begin{split}
\int_{B^{ +}_{2\varepsilon_{i}}(a_{i})} &
y^{1-2\gamma}\widetilde{\varphi_{i}}^{a_{
i }} \widetilde{\varphi_{j}}^{a_{j}}  
\leq 
\frac{C}{\lambda_{i}^{2}(1+\lambda_{i}^{2}\vert
a_{i}-a_{j}\vert^{2})^{\frac{n-2\gamma}{2}}}
\int_{[\vert z \vert \leq C\lambda_{i}\vert a_{i}-a_{j}\vert]}
\frac{y^{1-2\gamma}}{(\frac{\lambda_{i}}{\lambda_{j}}+\frac{\lambda_{j}}{\lambda
_{i}} \vert
z\vert^{2})^{\frac{n-2\gamma}{2}}} \\
\leq &
\frac{C(\frac{\lambda_{i}}{\lambda_{j}})^{\frac{n-2\gamma}{2}}}{\lambda_{i}^{2}
(1+\lambda_{i}^{2}\vert
a_{i}-a_{j}\vert^{2})^{\frac{n-2\gamma}{2}}}
\int_{[\vert z \vert \leq C\lambda_{i}\vert a_{i}-a_{j}\vert]}
y^{1-2\gamma}\vert z \vert^{2\gamma-n} 
\leq 
C\vert a_{i}-a_{j}\vert^{2}\varepsilon_{i,j}
=
o_{\frac{1}{\lambda_{i}}}(\varepsilon_{i,j}),
\end{split}
\end{equation*}
whence \;\eqref{I2_interaction_integral_estimate_proof}\; holds again and thus in 
any case. We are left with proving 
\begin{equation*}
\int_{B^{ +}_{2\varepsilon_{i}}(a_{i})}y^{1-2\gamma}\lambda_{i}y(\widetilde{
\varphi_{i}}^{a_{i}})^{\frac{n+2-2\gamma}
{ n-2\gamma } }\widetilde{\varphi_{j}}^{a_{j}}
=
o_{\frac{1}{\lambda_{i}}}(\varepsilon_{i,j}).
\end{equation*}
But this follows line by line as when showing  \;\eqref{I2_interaction_integral_estimate_proof}.
\end{pfn}

\smallskip

\noindent

\begin{pfn}{ of the estimate \eqref{epsilon_varepsilon_relation}}. 
We know
\begin{equation*}\begin{split}
\epsilon_{i,j}
= &
\underset{B_{ \varepsilon_{i}\lambda_{i} }(0)}{\int}
\frac{\frac{u_{ a _{j}}}{u_{a_{i}}}( 
\exp_{g_{a_{i}}}\frac{x}{\lambda_{i}})}{(1+r^{2})^{\frac{n+2\gamma}{2}}} 
(\frac{1}{\frac{\lambda_{i} }{ \lambda_{j} }+\lambda_{i}\lambda_{j}G_{ a 
_{j}}^{\frac{2}{2\gamma-n}}(\exp_{g_{ a _{i}}} 
\frac{x}{\lambda_{i}})})^{\frac{n-2\gamma}{2}}
dV_{h_{a_{i}}}.
\end{split}\end{equation*}
and have to show 
\begin{equation*}
\begin{split}
\epsilon_{i,j}
= &
c_{n,3}^{\gamma}\varepsilon_{i,j}(1+o_{\varepsilon_{i,j}}(1)).
\end{split}\end{equation*}
Since \;$
\eps_{i,j}^{\frac{2}{2\gamma-n}}
\sim 
\lambda_{i}\lambda_{j}G_{h}^{\frac{2}{2\gamma-n}}( a _{i}, a _{j})
\;\;\text{ or }\;\;
\eps_{i,j}^{\frac{2}{2\gamma-n}}
\sim
\frac{\lambda_{i}}{\lambda_{j}}
$, we may expand on 
\begin{equation*}\begin{split}
\mathcal{A}
=
\left(
\left[
\vert  \frac{x}{\lambda_{i}}\vert \leq \epsilon\sqrt{G^{\frac{2}{2-n}}_{ a 
_{j}}( a _{i})}
\right]
\cup
\left[
\vert  \frac{x}{\lambda_{i}}\vert \leq \epsilon  \frac{1}{\lambda_{j}} 
\right]
\right)
\cap
B_{\varepsilon_{i} \lambda_{i}}(0)
\subseteq
B_{\varepsilon_{i} \lambda_{i}}(0)
\end{split}\end{equation*}
for \;$\epsilon>0$\; small  \;$\frac{u_{a_{j}}}{u_{a_{i}}}(\exp_{g_{a_{i}}}\frac{x}{\lambda_{i}})=u_{a_{j}}(a_
{ i} )+\nabla 
u_{a_{j}}(a_{i})\frac{x}{\lambda_{i}}+O(\vert 
\frac{x}{\lambda_{i}} 
\vert^{2})$\; and 
\begin{equation*}\begin{split}
(
\frac{\lambda_{i} }{ \lambda_{j} }
& +
\lambda_{i}\lambda_{j}G_{a_{j}}^{\frac{2}{2\gamma-n}}(\exp_{g_{ a _{i}}} 
\frac{x}{\lambda_{i}}))^{\frac{2\gamma-n}{2}} \\
= &
(
\frac{\lambda_{i} }{ \lambda_{j} }
+
\lambda_{i}\lambda_{j}
G_{a _{j}}^{\frac{2}{2\gamma-n}}( a _{i})
)^{\frac{2\gamma-n}{2}}
+
\frac{2\gamma-n}{2}
\frac
{
\nabla G^{\frac{2}{2\gamma-n}}_{ a _{j}}( a _{i})\lambda_{j}x
+
O(\frac{ \lambda_{j} }{ \lambda_{i} }\vert x\vert^{2})
}
{
(
\frac{\lambda_{i} }{ \lambda_{j} }
+
\lambda_{i}\lambda_{j}G_{ a _{j}}^{\frac{2}{2\gamma-n}}( a _{i})
)^{\frac{n+2-2\gamma}{2}}
}.
\end{split}\end{equation*}
Thus  \;$\epsilon_{i,j}=\sum_{k=1}^{5}I_{k}+o_{\frac{1}{\lambda_{i}}}(\varepsilon_{i,j}
)$\; for
\begin{enumerate}[label=(\roman*)]
 \item \quad 
$
I_{1}
= 
\frac
{u_{ a _{j}}( a _{i})}
{(\frac{\lambda_{i} }{ \lambda_{j} }
+
\lambda_{i}\lambda_{j}
G_{ a _{j}}^{\frac{2}{2\gamma -n}}( a _{i})
)^{\frac{n-2\gamma}{2}}
}
\int_{\mathcal{A}}
\frac{1}{(1+r^{2})^{\frac{n+2\gamma}{2}}}
$
 \item \quad
$
I_{2}
= 
\frac
{1}
{(\frac{\lambda_{i}}{ \lambda_{j} }
+
\lambda_{i}\lambda_{j}
G_{ a _{j}}^{\frac{2}{2\gamma-n}}( a _{i})
)^{\frac{n-2\gamma}{2}}
}
\int_{\mathcal{A}} \frac{\nabla 
u_{a_{j}}(a_{i})\frac{x}{\lambda_{i}}}{(1+r^{2})^{\frac{n+2\gamma}{2}}}
$
 \item \quad
$
I_{3}
=
\frac
{\frac{2\gamma-n}{2}u_{ a _{j}}( a _{i})}
{(\frac{\lambda_{i} }{ \lambda_{j} }
+
\lambda_{i}\lambda_{j}
G_{a _{j}}^{\frac{2}{2\gamma-n}}( a _{i})
)^{\frac{n+2-2\gamma}{2}}
}
\int_{\mathcal{A}} \frac{\nabla G^{\frac{2}{2\gamma-n}}_{a _{j}}(a 
_{i})\lambda_{j}x}{(1+r^{2})^{\frac{n+2\gamma}{2}}}
$
 \item \quad
$
I_{4}
= 
\frac
{1}
{(\frac{\lambda_{i} }{ \lambda_{j} }
+
\lambda_{i}\lambda_{j}
G_{ a _{j}}^{\frac{2}{2\gamma-n}}( a _{i})
)^{\frac{n+2-2\gamma}{2}}
}
\int_{\mathcal{A}}
\frac{O(\frac{ \lambda_{j} }{ \lambda_{i} }\vert x 
\vert^{2})}{(1+r^{2})^{\frac{n+2\gamma}{2}}}
$
 \item \quad
$
I_{5}
 = 
\underset{\mathcal{A}^{c}}{\int} \frac{ \frac{u_{ a _{j}}}{u_{a_{i}}}( 
\exp_{g_{a_{i}}}\frac{x}{\lambda_{i}})}{(1+r^{2})^{\frac{n+2\gamma}{2}}}
 (\frac{1}{\frac{\lambda_{i} }{ \lambda_{j} }+\lambda_{i}\lambda_{j}G_{ a 
_{j}}^{\frac{2}{2\gamma-n}}(\exp_{g_{ a _{i}}} 
\frac{x}{\lambda_{i}})})^{\frac{n-2\gamma}{2}}dV_{h_{a_{i}}}.
$
\end{enumerate}
We then find with \;$c_{n,3}^{\gamma}=\int_{\R^{n}}(\frac{1}{1+r^{2}})^{\frac{n+2\gamma}{2}}
$
\begin{equation*}\begin{split}
I_{1}
= &
c_{n,3}^{\gamma}\frac
{u_{ a _{j}}( a _{i})}
{(\frac{\lambda_{i} }{ \lambda_{j} }
+
\lambda_{i}\lambda_{j}
G_{ a _{j}}^{\frac{2}{2\gamma-n}}(a _{i})
)^{\frac{n-2\gamma}{2}}
}
(1+o_{\varepsilon_{i,j}}(1)),
\end{split}\end{equation*}
whereas  \;$I_{2}=I_{3}=0$\; by radial symmetry and \;
$
I_{4}
= 
o_{\varepsilon_{i,j}}(\eps_{i,j})
$. Moreover  
\begin{equation*}
\begin{split}
I_{5}=o_{\varepsilon_{i,j}}(\eps_{i,j})
\end{split}
\end{equation*} 
in case \;$\eps_{i,j}^{\frac{2}{2\gamma-n}}\sim \frac{\lambda_{i}}{\lambda_{j}}$.  Else we have \;$\frac{\lambda_{j}}{\lambda_{i}}\leq 
\frac{\lambda_{i}}{\lambda_{j}}\leq 
\lambda_{i}\lambda_{j}G_{h}^{\frac{2}{2\gamma-n}}(a_{i},a_{j})$\; and decompose 
\begin{equation*}
\begin{split}
\mathcal{A}^{c}
=
B_{\varepsilon_{i}\lambda_{i}}\setminus \mathcal{A}
\subseteq &
\mathcal{B}_{1}\cup \mathcal{B}_{2},
\end{split}
\end{equation*}
where for a sufficiently large constant \;$E>0$ 
\begin{equation*}
\begin{split}
\mathcal{B}_{1}
= &
[\epsilon \sqrt{G^{\frac{2}{2\gamma-n}}_{a_{j}}(a_{i})}
\leq 
\vert \frac{x}{\lambda_{i}}\vert 
\leq 
E\sqrt{G^{\frac{2}{2\gamma-n}}_{a_{j}}(a_{i})}] 
\;\; \text{and} \;\;
\mathcal{B}_{2}
= 
[
E\sqrt{G^{\frac{2}{2\gamma-n}}_{a_{j}}(a_{i})}
\leq 
\vert \frac{x}{\lambda_{i}}\vert 
\leq \epsilon_{i}
].
\end{split}
\end{equation*} 
We may assume 
\;$
G_{a_{j}}(a_{i})^{\frac{2}{2\gamma-n}}\sim d^{2}_{g_{a_{j}}}(a_{j},a_{i})
\ll 1
$, 
since otherwise \;$I_{5}=o_{\varepsilon_{i,j}}(\varepsilon_{i,j})$, and estimate
\begin{equation*}
\begin{split}
I_{5}^{1}
= &
\underset{\mathcal{B}_{1}}{\int} \frac{\frac{ u_{a _{j}}}{u_{a_{i}}}(\exp_{g_{ 
a _{i}}} 
\frac{x}{\lambda_{i}})}{(1+r^{2})^{\frac{n+2\gamma}{2}}}dV_{h_{a_{i}}}
 (\frac{1}{\frac{\lambda_{i} }{ \lambda_{j} }+\lambda_{i}\lambda_{j}G_{ a 
_{j}}^{\frac{2}{2\gamma-n}}(\exp_{g_{ a _{i}}} 
\frac{x}{\lambda_{i}})})^{\frac{n-2\gamma}{2}} \\
\leq &
c\frac{(\frac{\lambda_{i}}{\lambda_{j}})^{\frac{n+2\gamma}{2}}}{(1+\lambda_{i}^{
2}G^{\frac{2}{2\gamma-n}}_{a_{j}}(a_{i}))^{\frac{n+2\gamma}{2}}} 
\underset{[\vert \frac{x}{\lambda_{j}}\vert \leq 
E\sqrt{G^{\frac{2}{2\gamma-n}}_{a_{j}}(a_{i})}]}{\int} 
 (\frac{1}{1+\lambda_{j}^{2}G_{ a _{j}}^{\frac{2}{2\gamma-n}}(\exp_{g_{ a 
_{i}}} \frac{x}{\lambda_{j}})})^{\frac{n-2\gamma}{2}}.
\end{split}
\end{equation*} 
Changing coordinates via \;$d_{i,j}=\exp_{g_{a_{i}}}^{-1}\exp_{g_{a_{j}}}$\; we get
\begin{equation*}
\begin{split}
I_{5}^{1}
\leq &
\frac{C}{(\frac{\lambda_{j}}{\lambda_{i}}+\lambda_{i}\lambda_{j}G^{\frac{2}{
2\gamma-n}}_{a_{j}}(a_{i}))^{\frac{n+2\gamma}{2}}}
\underset{[\vert \frac{x}{\lambda_{j}}\vert \leq  
C\sqrt{G_{a_{j}}^{\frac{2}{2\gamma-n}}}(a_{i})]}{\int} 
(\frac{1}{1+r^{2}})^{\frac{n-2\gamma}{2}}
\end{split}
\end{equation*} 
and thus \;$I_{5}^{1}=o_{\varepsilon_{i,j}}(\eps_{i,j})$. Moreover
\begin{equation*}
\begin{split}
I_{5,2}
= &
\underset{\mathcal{B}_{2}}{\int} \frac{ \frac{u_{ a _{j}}}{u_{a_{i}}}(\exp_{g_{ 
a _{i}}} 
\frac{x}{\lambda_{i}})}{(1+r^{2})^{\frac{n+2\gamma}{2}}}
 (\frac{1}{\frac{\lambda_{i} }{ \lambda_{j} }+\lambda_{i}\lambda_{j}G_{ a 
_{j}}^{\frac{2}{2\gamma-n}}(\exp_{g_{ a _{i}}} 
\frac{x}{\lambda_{i}})})^{\frac{n-2\gamma}{2}} \\
 \leq &
\frac{C}{(\frac{\lambda_{i} }{ \lambda_{j} }+\lambda_{i}\lambda_{j}G_{a 
_{j}}^{\frac{2}{2\gamma-n}}(a_{i}))^{\frac{n-2\gamma}{2}}} 
\underset{[\vert x \vert \geq 
c\sqrt{\lambda_{i}^{2}G^{\frac{2}{2\gamma-n}}_{a_{j}}(a_{i})}]}{\int}
\frac{1}{(1+r^{2})^{\frac{n+2\gamma}{2}}}
=
o_{\varepsilon_{i,j}}(\eps_{i,j}).
\end{split}
\end{equation*} 
Therefore \;$I_{5}=I_{5}^{1}+I_{5}^{2}=o_{\varepsilon_{i,j}}(\varepsilon_{i,j})$. 
Collecting terms we get 
\begin{equation*}
\begin{split}
\epsilon_{i,j}
= &
c_{n,3}^{\gamma}\frac
{u_{ a _{j}}( a _{i})}
{(\frac{\lambda_{i} }{ \lambda_{j} }
+
\lambda_{i}\lambda_{j}
G_{ a _{j}}^{\frac{2}{2\gamma-n}}(a _{i})
)^{\frac{n-2\gamma}{2}}
}
+
o_{\varepsilon_{i,j}}(\varepsilon_{i,j}).
\end{split}
\end{equation*}
By conformal covariance of the Green's function, cf. \eqref{covariance_greens_function}, we conclude
\;$
\epsilon_{i,j}
= 
c_{n,3}^{\gamma}\varepsilon_{i,j}(1+o_{\varepsilon_{i,j}}(1)).
$
\end{pfn}

\end{document}